\documentclass[10pt,reqno]{extarticle}
\pdfoutput=1
\usepackage{geometry}
\usepackage{graphicx}

\usepackage{tocloft}
\usepackage{titlesec}
\titleformat{\section}{\normalfont\rmfamily\large\bfseries\color{black}}{\thesection}{0.8em}{}
\titleformat{\subsection}{\normalfont\rmfamily\normalsize\bfseries\color{black}}{\thesubsection}{0.8em}{}
\titleformat{\subsubsection}{\normalfont\rmfamily\normalsize\bfseries\color{black}}{\thesubsubsection}{0.8em}{}

\usepackage[utf8]{inputenc}
\usepackage[english]{babel}

\usepackage{amssymb}
\usepackage{amsmath}
\usepackage{amsthm}

\usepackage{tikz}
\usepackage{pgfplots}

\usepackage{float}
\usepackage{caption}
\usepackage{subcaption}
\usepackage{scalerel}
\usepackage{enumitem}
\usepackage{hyperref}
\hypersetup{
  colorlinks   = true,
  urlcolor     = blue,
  linkcolor    = blue,
  citecolor   = red
}
\usepackage{showlabels}

\numberwithin{equation}{section}
\allowdisplaybreaks

\newtheorem{theorem}{Theorem}[section]

\newtheorem{definition}[theorem]{Definition}
\newtheorem{subdefinition}{Definition}[theorem]

\newtheorem{corollary}[theorem]{Corollary}
\newtheorem{remark}[theorem]{Remark}
\newtheorem{lemma}[theorem]{Lemma}

\newtheorem{proposition}[theorem]{Proposition}
\newtheorem{subproposition}{Proposition}[theorem]

\newtheorem{example}[theorem]{Example}
\newtheorem{subexample}{Example}[theorem]

\newtheorem{problem}[theorem]{Problem}

\newtheorem{construction}[theorem]{Construction}

\newcommand{\calD}{\mathcal{D} \makebox[0ex]{}}

\newcommand{\defeq}{\mathrel{\mathop:}=}
\newcommand{\slashg}{g  \mkern-8.2mu \scaleto{\boldsymbol{\slash}}{1.6ex} \mkern+1mu \makebox[0ex]{}}
\newcommand{\circg}{\mathring{g} \makebox[0ex]{}}
\newcommand{\dvol}{\mathrm{dvol} \makebox[0ex]{}}
\newcommand{\calC}{\mathcal{C} \makebox[0ex]{}}

\newcommand{\uC}{\underline{C} \makebox[0ex]{}}

\newcommand{\calV}{\mathcal{V} \makebox[0ex]{}}

\newcommand{\barf}{\bar{f} \makebox[0ex]{}}

\newcommand{\baru}{\bar{u} \makebox[0ex]{}}

\newcommand{\calE}{\mathcal{E} \makebox[0ex]{}}

\newcommand{\ed}{\mathrm{d} \makebox[0ex]{}}

\newcommand{\rmI}{\mathrm{I} \makebox[0ex]{}}

\newcommand{\barE}{\bar{E} \makebox[0ex]{}}

\newcommand{\bart}{\bar{t} \makebox[0ex]{}}
\newcommand{\barx}{\bar{x} \makebox[0ex]{}}

\title{\Large \textsc{Lorentz polarisation and isoperimetric inequality in Minkowski spacetime}}
\author{Pengyu Le}
\newcommand{\Address}{{
  \bigskip
  \footnotesize
  \textsc{Beijing Institute of Mathematical Sciences and Applications, Beijing, China}
  
  \textit{E-mail address}: \texttt{pengyu.le@bimsa.cn}
}}
\date{}

\begin{document}

\maketitle

\begin{abstract}
In this paper, we prove an isoperimetric inequality for the domain of dependence of a finite lightcone in the Minkowski spacetime of dimension greater than or equal to 3. The inequality involves two quantities: the volume of the domain of dependence, and the perimeter of the finite lightcone. It states that among all finite lightcones with the same perimeter, the maximal volume of the domain of dependence is achieved by the spacelike hyperplane truncated finite lightcone. A novelty of this isoperimetric inequality is the codimension 2 comparison feature.

We introduce the Lorentz polarisation to prove the isoperimetric inequality by studying the corresponding variational problem. A key observation is the monotonicity of the domain of dependence of a finite lightcone under the Lorentz polarisation. We show that any finite lightcone can be transformed by Lorentz polarisations to approximate a spacelike hyperplane truncated finite lightcone with an equal or less perimeter.

As further applications of the method of Lorentz polarisation, we prove the following isoperimetric type inequalities: a) For a set with the given perimeter in the hyperboloid in the Minkowski spacetime, the geodesic ball in the hyperboloid has the maximal volume of the domain of dependence of the set; b) For an achronal hypersurface with boundary in the lightcone (or the hyperboloid), given the perimeter of the boundary fixed, the spacelike hyperplane disk has the maximal area.
\end{abstract}

\tableofcontents

\section{Introduction}\label{sec 1}
The well-known classical isoperimetric inequality states that among all the shapes with the same perimeter in the Euclidean space $\mathbb{E}^n$, the ball has the maximal volume. The proofs and various applications of the isoperimetric inequality in Euclidean space can be found in the literature (\cite{PS51} \cite{Pa67} \cite{O78} \cite{B80} \cite{BZ88} \cite{C01}).

\subsection{Isoperimetric inequality for domain of dependence of finite lightcone}
In this paper, we investigate the isoperimetric inequality problem in the Minkowski spacetime $\mathbb{M}^{n,1}$.  The naive formal generalisation of the Euclidean isoperimetric inequality problem to the case of the Minkowski spacetime proves immediately problematic, since the area of a null hypersurface is zero. Therefore to formulate a meaningful isoperimetric inequality problem in the Minkowski spacetime is already an interesting and nontrivial question.

One direction to formulate the isoperimetric inequality in the Minkowski spacetime or general Lorentz manifolds is to investigate the isoperimetric problem on a spacelike hypersurface, for example \cite{CGGK07} 
\cite{TW22}.

Another direction to formulate the inequality is to investigate the relation between the area of a spacelike hypersurface in the lightcone and the volume of the cone bounded by the spacelike hypersurface, for example \cite{Ba99} \cite{BaE99} \cite{BaH01} \cite{ACKW09}.

In this paper, we formulate the following new type of the isoperimetric inequality problem and give its solution.
\begin{problem}\label{prob 1.1}
Let $S$ be a spacelike surface (curve if $n=2$) of the lightcone in the Minkowski spacetime $\mathbb{M}^{n,1}$ and $C_S$ be the finite lightcone bounded by $S$. Let $\calD(C_S)$ be the domain of dependence of $C_S$. Find the supremum of the volume of $\calD(C_S)$ among all $S$ with the same perimeter. See figure \ref{fig 1}.

Moreover determine whether the supremum can be achieved, and if it can be achieved, then find out the shape of $S$ achieving the supremum.

\begin{figure}[h]
\begin{center}
\begin{tikzpicture}[scale=0.7]
\draw[fill] (0,0) circle (1pt);
\draw (0,0) node[below right] {$o$}  -- (4,4) (0,0) -- (-4,4);
\draw[thick, dashed] (3,3) to [out=45,in=15] (2,3.5)
to [out=-165,in=-40] (-2,3)
to [out=140,in=135] (-3.5,3.5); 
\draw[dashed] (0,0) -- (2,3.5)  (0,0) -- (-2,3);
\draw[->] (1,1.5) to [out=-45,in=180] (2,1) node[right] {$C_S$};
\path[fill=blue,opacity=0.1] 
(-3.5,3.5) to [out=-45,in=-160] (-1.8,3.2)
to [out=20,in=-170] (-1,3.5)
-- (-0.4,5.8) to [out=120,in=-70] (-0.5,6.3)
to [out=-150,in=10] (-1.5,5.8)
-- (-3.6,3.6) -- cycle; 
\path[fill=blue,opacity=0.1] 
(-1,3.5)
to [out=10,in=150] (1,2.5)
to [out=-30,in=-135] (3,3)
--  (3.12,3.1) -- (1,6)
to [out=-160,in=-20] (-0.5,6.3)
to [out=-70,in=120] (-0.4,5.8)
-- cycle; 
\draw (-3.6,3.6) -- (-1.5,5.8)
(-1,3.5) -- (-0.4,5.8)
(3.12,3.1) -- (1,6); 
\draw[thick] (-3.5,3.5) to [out=-45,in=-160] (-1.8,3.2)
to [out=20,in=-170] (-1,3.5)
to [out=10,in=150] (1,2.5)
to [out=-30,in=-135] (3,3); 
\path[fill=gray,opacity=0.3] 
(0,0) -- 
(-3.5,3.5) to [out=-45,in=-160] (-1.8,3.2)
to [out=20,in=-170] (-1,3.5)
to [out=10,in=150] (1,2.5)
to [out=-30,in=-135] (3,3)
-- (0,0); 
\draw (0,0) -- (-1,3.5);
\draw[->] (0,5.5) to [out=90,in=180] (1,6.5) node[right] {$\mathcal{D}(C_S)$};
\draw[thick,dashed] (-0.5,6.3) to [out=-150,in=10] (-1.5,5.8)
(-0.5,6.3) to [out=-70,in=120] (-0.4,5.8)
(-0.5,6.3) to [out=-20,in=-160] (1,6); 
\end{tikzpicture}
\end{center}
\caption{Finite lightcone $C_S$ and its domain of dependence $\calD(C_S)$.}
\label{fig 1}
\end{figure}
\end{problem}
The following isoperimetric inequality answers the above problem.
\begin{theorem}[Rough version]\label{thm 1.2}
Among all the spacelike surfaces with the same area (curves with the same length if $n=2$) in the lightcone of the Minkowski spacetime $\mathbb{M}^{n,1}$, the spacelike hyperplane section bounds the finite lightcone whose domain of dependence has the maximal volume.

Let $S$ be the spacelike surface in the lightcone, $C_S$ be the open finite lightcone bounded by $S$, and $\calD(C_S)$ be the domain of dependence of $C_S$, the following isoperimetric inequality holds
\begin{align}
	\frac{|\calD(C_S)|}{2\omega_n/(n+1)} 
	\leq 
	\Big(\frac{|S|}{n \omega_n} \Big)^{\frac{n+1}{n-1}}.
	\label{eqn 1.1}
\end{align}
where $\omega_n$ is the volume of the unit $n$-dimensional ball and the inequality is achieved when $S$ is a spacelike hyperplane section of the lightcone. See figure \ref{fig 2}.
\begin{figure}[h]
\begin{center}
\begin{tikzpicture}[scale=0.7]
\node[right] at (1.5,1.5) {$C_S$};
\path[fill=gray!30] (-2.8,2.8)
to [out=-45,in=-135] (2.8,2.8)
-- (0,0)
-- cycle;
\path[fill=blue!10] (-2.8,3.2)
to [out=-135,in=135] (-2.8,2.8)
to [out=-45,in=-135] (2.8,2.8)
to [out=45, in=-45] (2.8,3.2)
-- (0,6)
-- cycle;
\draw 
(0,0) -- (-4,4)
(0,0) -- (4,4);
\draw
(-3,3) -- (0,6)
(3,3) -- (0,6);
\draw[thick] (-2.8,3.2)
to [out=-135,in=135] (-2.8,2.8)
to [out=-45,in=-135] (2.8,2.8)
to [out=45, in=-45] (2.8,3.2);
\draw[thick,dashed] (-2.8,3.2)
to [out=45,in=135] (2.8,3.2);
\node[right] at (2.8,2.8) {$S$};
\draw 
(0,0) -- (-0.8,1.7)
(0,0) -- (1.4,1.9);
\draw[dashed] 
(0,0) -- (-0.6,4.4);
\draw 
(0,6) -- (-0.8,1.7)
(0,6) -- (1.4,1.9);
\draw[dashed] 
(0,6) -- (-0.6,4.4);
\end{tikzpicture}
\end{center}
\caption{$S$ a spacelike hyperplane section.}
\label{fig 2}
\end{figure}
\end{theorem}

In fact, the above isoperimetric inequality is proved for a more general class of finite lightcones which includes the ones with much rougher boundaries (definitions \ref{def 2.1.a}, \ref{def 2.1.b}). See theorem \ref{thm 6.1} and \ref{thm 6.2} for the precise statements. We also show that the above isoperimetric inequality in the Minkowski spacetime has an implication in the Euclidean space, see theorem \ref{thm 6.3}.

\subsection{Lorentz reflection and polarisation}
The Euclidean isoperimetric inequality can be proved using the Euclidean reflection and polarisation. See \cite{Bu09} \cite{Bae19} for expositions on the Euclidean polarisation and its application to the Euclidean isoperimetric inequality. The polarisation has many applications to other problems, for example \cite{W52} \cite{BT76} \cite{Be84} \cite{BS00} \cite{AF04}.

In this paper, we introduce the Lorentz polarisation to prove the isoperimetric inequality for the domain of dependence of a finite lightcone. The Lorentz polarisation is a generalisation of the Euclidean polarisation, which we briefly describe here.

\begin{definition}
Let $v$ be a timelike vector at the origin $o$ and $H$ be a timelike hyperplane passing through $o$, that $v\not\in H$. Let $\gamma$ be the Lorentz reflection about $H$. Then the Lorentz polarisation of a set $E$ about the pair $(H,v)$ is the following set
\begin{align*}
	(E \cap H)
	\ \cup\ 
	[(E \cup \gamma(E) ) \cap H_+]
	\ \cup\ 
	[(E \cap \gamma(E) ) \cap H_-],
\end{align*}
where $H_+$ is the open half space containing $v$ and $H_-$ is the other one not containing $v$, see figure \ref{fig 3}. Roughly speaking, the Lorentz polarisation moves more part to the one side of $H$ by Lorentz reflection. See definition \ref{def 3.2} for more detailed description.
\begin{figure}[h]
\hfill
\begin{subfigure}[h]{0.45\textwidth}
\centering
\begin{tikzpicture}[scale=0.5]
\draw (-8 *.7,8 *.7) -- (0,0) node[below] {$o$} -- (8 *.7,8 *.7);
\draw (-1 *.7,1 *.7) -- (-2 *.7,2 *.7) -- (4 *.7,8 *.7) -- (5 *.7,7 *.7) -- (-1 *.7,1 *.7);
\path[fill=gray,opacity=0.3]
(-1 *.7,1 *.7) -- (-2 *.7,2 *.7) -- (4 *.7,8 *.7) -- (5 *.7,7 *.7) circle;
\node at (3 *.7,6 *.7) {$E$};
\draw (1.5 *.7,1.5 *.7) -- (-2.5 *.7,5.5 *.7) -- (-1 *.7,7 *.7) -- (3 *.7,3 *.7) -- (1.5 *.7,1.5 *.7);
\path[fill=blue,opacity=0.2]
(1.5 *.7,1.5 *.7) -- (-2.5 *.7,5.5 *.7) -- (-1 *.7,7 *.7) -- (3 *.7,3 *.7) circle;
\node at (-1 *.7, 5.5 *.7) {$\gamma(E)$};
\draw[->] (0,0) -- (-.1,1) node[above] {$v$};
\node[above] at (4 *.7+1,8 *.7) {$H_{-}$};
\node[above] at (-1 *.7-2,7 *.7) {$H_{+}$};
\draw (0,0) -- (1 *.7,5 *.7) -- (1.6 *.7,8 *.7) node[right] {$H$};
\end{tikzpicture}
\end{subfigure}
\hfill
\begin{subfigure}[h]{0.45\textwidth}
\centering
\begin{tikzpicture}[scale=0.5]
\draw (-8 *.7,8 *.7) -- (0,0) node[below] {$o$} -- (8 *.7,8 *.7);
\draw (-1 *.7,1 *.7) -- (-2 *.7,2 *.7) -- (4 *.7,8 *.7) -- (5 *.7,7 *.7) -- (-1 *.7,1 *.7);
\draw (1.5 *.7,1.5 *.7) -- (-2.5 *.7,5.5 *.7) -- (-1 *.7,7 *.7) -- (3 *.7,3 *.7) -- (1.5 *.7,1.5 *.7);
\path[fill=blue,opacity=0.2]
(-1 *.7,1 *.7) -- (-2 *.7,2 *.7) -- (-2 *.7+1.5*.7,2 *.7+1.5*.7) -- (-2.5 *.7,5.5 *.7) -- (-1 *.7,7 *.7)  -- (-1 *.7+3*.7, 1 *.7+3*.7) circle;
\draw[->] (0,0) -- (-.1,1) node[above] {$v$};
\node[above] at (4 *.7+1,8 *.7) {$H_{-}$};
\node[above] at (-1 *.7-2,7 *.7) {$H_{+}$};
\draw (0,0) -- (1 *.7,5 *.7) -- (1.6 *.7,8 *.7) node[right] {$H$};
\end{tikzpicture}
\end{subfigure}
\caption{Lorentz polarisation of $E$ about $(H,v)$}
\label{fig 3}
\end{figure}
\end{definition}

A key property of the Lorentz polarisation for the domain of dependence of a finite lightcone is that the Lorentz polarisation does not change the perimeter of the finite lightcone, while perserves or increases the volume of the domain of dependence. See proposition \ref{prop 3.10} and corollary \ref{coro 3.11}.

An important step in the proof is to show that starting from an arbitrary finite lightcone, one can apply Lorentz polarisations to obtain a sequence of finite lightcones converging to a spacelike hyperplane truncated finite lightcone with an equal or less perimeter. See proposition \ref{prop 4.16} and lemma \ref{lem 4.22}. Then the isoperimetric inequality follows.

We identify the case of equality in the isoperimetric inequality by introducing the equal perimeter separation hyperplane, see definition \ref{def 5.3}. We can show that for a finite lightcone which achieves the equality, no null generator in the strict future boundary of the domain of dependence can cross the equal perimeter separation hyperplane, see proposition \ref{prop 5.8}. Applying this proposition, we can show that the equality is achieved if and only if the finite lightcone is spacelike hyperplane truncated, see proposition \ref{prop 5.11} and \ref{prop 5.13}.

\subsection{Other isoperimetric type inequalities}

The method of Lorentz polarisation can be employed to investigate other isoperimetric type inequalities.

\subsubsection{Domain of dependence of set in hyperboloid}
In section \ref{sec 7}, we study the isoperimetric inequality for the domain of dependence of the set in the hyperboloid, see theorem \ref{thm 7.7}. We briefly explain the inequality.
\begin{figure}[h]
\hfill
\begin{subfigure}[h]{0.45\textwidth}
\centering
\begin{tikzpicture}[scale=.8, 
	declare function={
			h1(\x,\y)=sqrt(\x^2 + \y^2 + 1);
			h2(\x,\y)=sqrt(\x^2 + \y^2);
			rho(\x)=.8*(cos(\x))^2+.8;
		}
	]
	\begin{axis}[view = {50}{15},
		xmin=-2, xmax=2,
		ymin=-2, ymax=2,
		zmin=0, zmax=3,
		xlabel = $x_1$,
		ylabel = $x_2$,
		ticklabel style = {font = \scriptsize},
		unit vector ratio = 1 1 1,
		grid
	]
	\addplot3 [domain=-2:2, 
		domain y=-2:2, 
		samples=20, 
		samples y=20, 
		surf,
		shader = interp, 
		opacity=1,
	] 
		{h1(\x,\y)} ;
	\addplot3 [domain=-1.5:1.5, 
		domain y=-1.5:1.5, 
		samples=20, 
		samples y=20, 
		surf,
		shader = interp,
		opacity=0.4,
		fill=green!30!red!20,
	]
		{h2(\x,\y)};
	\end{axis}
\end{tikzpicture}
\end{subfigure}
\hfill
\begin{subfigure}[h]{0.45\textwidth}
\centering
\begin{tikzpicture}[scale=.8, 
	declare function={
			h1(\x,\y)=sqrt(\x^2 + \y^2 + 1);
			h2(\x,\y)=sqrt(\x^2 + \y^2);
			rho(\x)=(cos(\x))^2+1.2;
		}
	]
	\begin{axis}[view = {50}{15},
		xmin=-2, xmax=2,
		ymin=-2, ymax=2,
		zmin=0, zmax=3,
		xlabel = $x_1$,
		ylabel = $x_2$,
		ticklabel style = {font = \scriptsize},
		unit vector ratio = 1 1 1,
		grid
	]
	\addplot3 [
	domain=0:360,
	domain y=0:1,
	samples=30, 
	samples y=10,
	surf, 
	shader = interp, 
	colormap/bluered,
	] 
		( {y*rho(x)*cos(x)}, 
			{y*rho(x)*sin(x)}, 
			{h1({y*rho(x)*cos(x)}, {y*rho(x)*sin(x)})}
		);
	\addplot3 [blue,
	domain=0:360,
	samples=20, 
	] 
		( {rho(x)*cos(x)}, 
			{rho(x)*sin(x)}, 
			{h1({rho(x)*cos(x)}, {rho(x)*sin(x)})}
		);
	\end{axis}
\end{tikzpicture}
\end{subfigure}
\caption{Hyperboloid $S_{-1}$ and the set in $S_{-1}$.}
\label{fig 4}
\end{figure}

\begin{theorem}[Rough version]\label{thm 1.4}
Let $S_{-1}$ be the hyperboloid $\{ -t^2 + x_1^2 + \cdots + x_n^2=-1 \}$ in $\mathbb{M}^{n,1}$. Let $E$ be a set in $S_{-1}$ and $\partial E$ be the boundary of $E$ in $S_{-1}$, see figure \ref{fig 4}. Consider the domain of dependence $\calD(E)$ of $E$ in $\mathbb{M}^{n,1}$. Then we have the following isoperimetric inequality
\begin{align}
	\frac{|\calD(E)|}{2\omega_n/(n+1)} 
	\leq
	\Big(\frac{|\partial E|}{n \omega_n} \Big)^{\frac{n+1}{n-1}},
	\label{eqn 1.2}
\end{align}
where the equality is achieved by the set truncated by the spacelike hyperplane, which is a geodesic ball in $S_{-1}$.
\end{theorem}

In the precise formulation of theorem \ref{thm 7.7}, we assume that $E$ is a set of finite perimeter and use its perimeter $P(E)$ in the inequality, see figure \ref{fig 5}.
\begin{figure}[h]
\centering
\begin{tikzpicture}[scale=.7, 
	declare function={
			h1(\x,\y)=sqrt(\x^2 + \y^2 + 1);
			rho=1.5;
			b1(\x,\y)= rho+sqrt((rho)^2+1)-sqrt(\x^2 + \y^2);
			b2(\x,\y)= -rho+sqrt((rho)^2+1)+sqrt(\x^2 + \y^2);
		}
	]
	\begin{axis}[view = {50}{15},
		xmin=-2, xmax=2,
		ymin=-2, ymax=2,
		zmin=0, zmax=3,
		xlabel = $x_1$,
		ylabel = $x_2$,
		ticklabel style = {font = \scriptsize},
		unit vector ratio = 1 1 1,
		grid
	]
	\addplot3 [domain=-2.2:2.2, 
		domain y=-2.2:2.2, 
		samples=20, 
		samples y=20, 
		surf,
		shader = interp, 
		opacity=0.1,
	] 
		{h1(\x,\y)} ;
	\addplot3 [
	domain=0:360,
	domain y=0:1,
	samples=30, 
	samples y=10,
	surf, 
	shader = interp, 
	colormap/violet,
	] 
		( {y*rho*cos(x)}, 
			{y*rho*sin(x)}, 
			{h1({y*rho*cos(x)}, {y*rho*sin(x)})}
		);
	\addplot3 [blue,
	domain=0:360,
	samples=20, 
	] 
		( {rho*cos(x)}, 
			{rho*sin(x)}, 
			{h1({rho*cos(x)}, {rho*sin(x)})}
		);

	\addplot3 [
	domain=0:360,
	domain y=0:1,
	samples=30, 
	samples y=10,
	surf, 
	shader = interp, 
	colormap/bluered,
	opacity=0.4,
	] 
		( {y*rho*cos(x)}, 
			{y*rho*sin(x)}, 
			{b1({y*rho*cos(x)}, {y*rho*sin(x)})}
		);

	\addplot3 [
	domain=0:360,
	domain y=0:1,
	samples=30, 
	samples y=10,
	surf, 
	shader = interp, 
	colormap/bluered,
	opacity=0.4,
	] 
		( {y*rho*cos(x)}, 
			{y*rho*sin(x)}, 
			{b2({y*rho*cos(x)}, {y*rho*sin(x)})}
		);
	\end{axis}
\end{tikzpicture}
\caption{Case of equality.}
\label{fig 5}
\end{figure}

\subsubsection{Achronal hypersurface with boundary in lightcone}
In section \ref{sec 8}, we study the isoperimetric inequality for the achronal hypersurface with boundary in the lightcone, see theorem \ref{thm 8.1}. We briefly explain the inequality below.

\begin{figure}[h]
\centering
\begin{tikzpicture}[scale=.8,
	declare function={
			c1(\x,\y)=sqrt(\x^2 + \y^2);
			c2(\r)=max(\r, \r - ( 7.5*((\r-.7)^4) + 7*((\r-.7)^3) ) );
			rho(\x)=.6(cos(\x))^2+.8;
			h1(\y)= 1.1-.1*\y^2;
		}
	]
	\begin{axis}[view = {35}{15},
		xmin=-2, xmax=2,
		ymin=-2, ymax=2,
		zmin=0, zmax=3,
		xlabel = $x_1$,
		ylabel = $x_2$,
		ticklabel style = {font = \scriptsize},
		unit vector ratio = 1 1 1,
		grid
	]
	\addplot3 [domain=-2.2:2.2,
		domain y=-2.2:2.2,
		samples=20,
		samples y=20,
		surf,
		shader = interp,
		opacity=0.3,
		colormap/bluered,
	] 
		{c1(\x,\y)} ;
	\addplot3 [domain=0:360,
		domain y=0:1,
		samples=20,
		samples y=10,
		surf,
		shader = interp,
		opacity=.8,
		colormap/violet,
	] 
		(	{y*rho(x)*cos(x)}, 
			{y*rho(x)*sin(x)} , 
			{h1(y)*c2( y*rho(x))}
		);
	\addplot3 [domain=0:360,
		samples=40,
		opacity=.8,
		violet,
	] 
		(	{rho(x)*cos(x)}, 
			{rho(x)*sin(x)} , 
			{c2( rho(x))}
		);
	\end{axis}
\end{tikzpicture}
\caption{Achronal hypersurface with boundary in the lightcone}
\label{fig 6}
\end{figure}

\begin{theorem}[Rough version]\label{thm 1.5}
Let $\Sigma$ be a closed achronal hypersurface with $\partial \Sigma \subset C_0$, see figure \ref{fig 6}. Then we have the isoperimetric inequality similar as the one in the Euclidean space
\begin{align}
	\frac{|\Sigma|}{\omega_n}
	\leq
	\Big( \frac{|\partial \Sigma|}{n \omega_n} \Big)^{\frac{n}{n-1}},
	\label{eqn 1.3}
\end{align}
where the equality is achieved by the spacelike hyperplane disk.
\end{theorem}

In the precise formulation of theorem \ref{thm 8.1}, in order to include more general cases, instead of considering the boundary of $\Sigma$, we consider the closed achronal hypersurface $\Sigma$ in $I^+(o)$ which is contained in the domain of dependence of an open finite lightcone $C_f$ and use the perimeter of $C_f$ in the inequality.

The isoperimetric inequality \eqref{eqn 1.3} leads to a functional inequality on the hyperbolic space as a corollary, see corollary \ref{coro 8.10} and inequality \eqref{eqn 8.3}.

\subsubsection{Achronal hypersurface with boundary in hyperboloid}

In section \ref{sec 9}, we study the isoperimetric inequality for the archronal hypersurface with boundary in the hyperboloid $S_{-1}$, see theorem \ref{thm 9.1}. We briefly explain the inequality.

\begin{figure}[h]
\centering
\begin{tikzpicture}[scale=.8, 
	declare function={
			h1(\x,\y)=sqrt(\x^2 + \y^2 + 1);
			rho(\x)=.6(cos(\x))^2+1;
			h2(\y)= 1.8-.8*\y;
		}
	]
	\begin{axis}[view = {25}{10},
		xmin=-2, xmax=2,
		ymin=-2, ymax=2,
		zmin=0, zmax=3,
		xlabel = $x_1$,
		ylabel = $x_2$,
		ticklabel style = {font = \scriptsize},
		unit vector ratio = 1 1 1,
		grid
	]
	\addplot3 [domain=-2.2:2.2, 
		domain y=-2.2:2.2, 
		samples=20, 
		samples y=20, 
		surf,
		shader = interp, 
		opacity=0.2,
	] 
		{h1(\x,\y)} ;
	\addplot3 [
	domain=0:360,
	domain y=0:1,
	samples=30, 
	samples y=10,
	surf, 
	shader = interp, 
	colormap/violet,
	] 
		( {y*rho(x)*cos(x)}, 
			{y*rho(x)*sin(x)}, 
			{h2(y)*h1({y*rho(x)*cos(x)}, {y*rho(x)*sin(x)})}
		);
	\addplot3 [blue,
	domain=0:360,
	samples=20, 
	] 
		( {rho(x)*cos(x)}, 
			{rho(x)*sin(x)}, 
			{h1({rho(x)*cos(x)}, {rho(x)*sin(x)})}
		);
	\end{axis}
\end{tikzpicture}
\caption{Achronal hypersurface with boundary in the hyperboloid $S_{-1}$.}
\label{fig 7}
\end{figure}

\begin{theorem}[Rough version]\label{thm 1.6}
Let $\Sigma$ be a closed achronal hypersurface with $\partial \Sigma \subset S_{-1}$, see figure \ref{fig 7}. Then we have the isoperimetric inequality similar as the one in the Euclidean space
\begin{align}
	\frac{|\Sigma|}{\omega_n}
	\leq
	\Big( \frac{|\partial \Sigma|}{n \omega_n} \Big)^{\frac{n}{n-1}},
	\label{eqn 1.4}
\end{align}
where the equality is achieved by the spacelike hyperplane disk.
\end{theorem}

In the precise formulation of theorem \ref{thm 9.1}, we consider the closed achronal hypersurface $\Sigma$ contained in the domain of dependence of a set $E\subset S_{-1}$. We assume that $E$ is a set of finite perimeter in the hyperboloid $S_{-1}$ and use the perimeter $P(E)$ in the inequality.

\section{Notions in isoperimetric inequality for domain of dependence of finite lightcone}\label{sec 2}
In this section, we introduce the basic notions in the isoperimetric inequality for domain of dependence of finite lightcone. See the references \cite{P72}\cite{HE73} for the essential background knowledge.

\subsection{Finite lightcone and its boundary}\label{sec 2.1}
Let $S$ be a spacelike surface in the lightcone $C_0$ emanating from the origin $o$ of the Minkowski spacetime $(\mathbb{M}^{n,1},\eta)$. Denote the open finite lightcone bounded by $S$ as $C_S$. We call $S$ the boundary surface of the finite lightcone $C_S$.

We introduce the parameterisation of $S$ in the following. Let $\vartheta$ be a coordinate system on the round sphere $(\mathbb{S}^{n-1}, \circg)$ of radius $1$, which we abuse the notation to denote $\vartheta$ as the point of $\mathbb{S}^{n-1}$. Let $\{t, r, \vartheta\}$ be the spatial polar coordinate system of $\mathbb{M}^{n,1}$, where
\begin{align*}
	\eta
	=
	-\ed t^2 + \ed r^2 + r^2 \circg
\end{align*}
 then $\{r, \vartheta\}$ is a coordinate system of $C_0$. In this $\{r,\vartheta\}$ coordinate system, $S$ can be parameterised by a function $f$ as its graph of $r$ over the $\vartheta$ domain, i.e.
\begin{align}
S = \{ (r,\vartheta): r=f(\vartheta)\}.
\label{eqn 2.1}
\end{align}
In order to emphasize the parameterisation function $f$ of $S$, we sometimes use $S_f$ to denote the surface $S$, and $C_f$ to denote the open finite lightcone $C_S$, which is the set
\begin{align}
C_f = \{ (r,\vartheta): r < f(\vartheta)\} \cup \{ o \}
\label{eqn 2.2}
\end{align}
in $\{r,\vartheta\}$ coordinate system. See figure \ref{fig 8}.
\begin{figure}[h]
\begin{center}
\begin{tikzpicture}[scale=0.7]
\draw[fill] (0,0) circle (1pt);
\draw (0,0) node[below right] {$o$}  -- (4,4) (0,0) -- (-4,4);
\draw[thick, dashed] (3,3) to [out=45,in=15] (2,3.5)
to [out=-165,in=-40] (-2,3)
to [out=140,in=135] (-3.5,3.5); 
\draw[thick] (-3.5,3.5) to [out=-45,in=-160] (-1.8,3.2)
to [out=20,in=-170] (-1,3.5)
to [out=10,in=150] (1,2.5)
to [out=-30,in=-135] (3,3); 
\path[fill=gray,opacity=0.1] 
(0,0) -- (3,3)
to [out=45,in=15] (2,3.5)
to [out=-165,in=-40] (-2,3)
to [out=140,in=135] (-3.5,3.5)
-- (0,0); 
\path[fill=gray,opacity=0.3] 
(0,0) -- 
(-3.5,3.5) to [out=-45,in=-160] (-1.8,3.2)
to [out=20,in=-170] (-1,3.5)
to [out=10,in=150] (1,2.5)
to [out=-30,in=-135] (3,3)
-- (0,0); 
\draw (0,0) --(-1.8,3.2)
(0,0) -- (-1,3.5) node[above] {$f$} 
(0,0) -- (1,2.5);
\draw[dashed] (0,0) -- (2,3.5)  (0,0) -- (-2,3);
\draw[->] (1,1.5) to [out=-45,in=180] (2,1) node[right] {$C_f$};
\draw[->] (2,3.6) to [out=60,in=180] (3,4) node[right] {$S_f$};
\end{tikzpicture}
\end{center}
\caption{The closed finite future lightcone bounded by $f$.}
\label{fig 8}
\end{figure}
Note that the parameterisation formulae \eqref{eqn 2.1}\eqref{eqn 2.2} of $S_f$ and $C_f$ can be generalised easily to any function $f$, for example $f$ could be discontinuous. And by identifying the coordinate systems $(r,\vartheta)$ in $C_0$ and $\mathbb{R}^n$, one can view a finite lightcone $C_f$ as simply a star-shaped domain.  Such generalisation is necessary when considering the limit of a sequence of finite lightcones and their boundaries. Thus we introduce the following definition for a finite lightcone.

\addtocounter{theorem}{1}
\begin{subdefinition}\label{def 2.1.a}
An open finite lightcone in $C_0$ is a bounded open star-shaped domain at the origin $o$ by identifying $C_0$ with $\mathbb{R}^n$ using the $(r,\vartheta)$ coordinate system. Equivalently, an open finite lightcone in $C_0$ is a set $C_f$ of the following form in the $(r,\vartheta)$ coordinate system
\begin{align*}
C_f = \{ (r,\vartheta): r < f(\vartheta)\} \cup \{ o \},
\end{align*}
where  $f$ is positive, bounded and lower semicontinuous. $f$ is called the parameterisation function of the open finite lightcone $C_f$. Define $S_f$ be the graph of $f$ in the $(r,\vartheta)$ coordinate system of $C_0$
\begin{align*}
S_f = \{ (r,\vartheta): r = f(\vartheta) \},
\end{align*}
$S_f$ is called the lower envelope of the open finite lightcone $C_f$.

We define the closure and boundary of an open finite lightcone as its closure and boundary in $C_0$.
\end{subdefinition}

Similarly, we define the notion of a closed finite lightcone as follows.
\begin{subdefinition}\label{def 2.1.b}
A closed finite lightcone in $C_0$ is a bounded closed star-shaped domain at the origin $o$  with $o$ being an interior, by identifying $C_0$ with $\mathbb{R}^n$ using the $(r,\vartheta)$ coordinate system. Equivalently, a closed finite lightcone in $C_0$ is a set $\calC_h$ of the following form in the $(r,\vartheta)$ coordinate system\footnote{We use the calligraphic letter $\calC$ to denote the closed finite lightcone.}
\begin{align*}
\calC_h = \{ (r,\vartheta): r \leq h(\vartheta)\} \cup \{ o \},
\end{align*}
where  $h$ is positive, bounded and upper semicontinuous. $h$ is called the parameterisation function of the closed finite lightcone $\calC_h$. Define $S_h$ be the graph of $f$ in the $(r,\vartheta)$ coordinate system of $C_0$
\begin{align*}
S_h = \{ (r,\vartheta): r = h(\vartheta) \},
\end{align*}
$S_h$ is called the upper envelope of the closed finite lightcone $\calC_h$.

We define the interior and boundary of an open finite lightcone as its interior and boundary in $C_0$.
\end{subdefinition}

\begin{remark}
The lower and upper semicontinuity of $f$ and $h$ are equivalent to the openness of $C_f$ and the closedness of $\calC_h$ respectively.
\end{remark}

\subsection{Some point-set topological properties of finite lightcone}\label{sec 2.2}
When $f$ is continuous, then the boundary $\partial C_f$ of the open finite lightcone $C_f$ is simply the lower envelope $S_f$, while for discontinuous $f$,  $\partial C_f$ is more complicated. We have the following proposition for the boundary and closure of an open finite lightcone, whose proof is straightforward.

\addtocounter{theorem}{1}
\begin{subproposition}\label{prop 2.3.a}
Let $C_f$ be an open finite lightcone. The closure $\overline{C_f}$ of $C_f$ is given by
\begin{align*}
\overline{C_f} =  \{ (r,\vartheta): r \leq f^{\sup}(\vartheta)\} \cup \{ o \}.
\end{align*}
where
\begin{align*}
f^{\sup}(\vartheta) = \lim_{\delta \rightarrow 0^+} \sup_{B_{\vartheta}(\delta)}\{f(x)\}.
\end{align*}
$f^{\sup}$ is positive, bounded and upper semicontinuous, thus the closure $\overline{C_f}$ is the closed finite lightcone $\calC_{f^{\sup}}$. Let $Z$ be the set of discontinuous points of $f$, then the boundary $\partial C_f$ of $C_f$ is given by
\begin{align*}
\partial C_f = S_f \cup B_Z,
\quad
B_Z = \{ ( r, \vartheta): \vartheta \in Z, r\in [ f(\vartheta), f^{\sup} (\vartheta)] \}.
\end{align*}
\end{subproposition}

Similarly for the boundary and interior of a closed finite ligthcone, we have the following analogous proposition.
\begin{subproposition}\label{prop 2.3.b}
Let $\calC_h$ be a closed finite lightcone. The interior $(\calC_h)^{\circ}$ of $\calC_h$ is given by
\begin{align*}
(\calC_h)^{\circ} =  \{ (r,\vartheta): r < h_{\inf}(\vartheta)\} \cup \{ o \}.
\end{align*}
where 
\begin{align*}
h_{\inf} (\vartheta) = \lim_{\delta \rightarrow 0^+} \inf_{B_{\vartheta}(\delta)} \{ h(x) \}.
\end{align*}
$h_{\inf}$ is positive bounded and lower semicontinuous, thus the interior $(\calC_h)^{\circ}$ is the open finite lightcone $C_{h_{\inf}}$. Let $Z$ be the set of discontinuous points of $h$, then the boundary $\partial \calC_h$ of $\calC_h$ is given by
\begin{align*}
\partial \calC_h = S_h \cup B_Z,
\quad
B_Z = \{ ( r, \vartheta): \vartheta \in Z, r\in [ h_{\inf}(\vartheta), h(\vartheta)] \}.
\end{align*}
\end{subproposition}

\begin{example}\label{ex 2.4}
A simple example to illustrate the above is a function $f$ with one jump point on the circle which defines an open finite lightcone in $\mathbb{M}^{2,1}$, whose boundary contains the light ray at the jump point. See figure \ref{fig 9}.
\begin{figure}[h]
\hfill
\begin{subfigure}[b]{0.45\textwidth}
\centering
\begin{tikzpicture}[scale=.8]
\draw[dashed,color=black,domain=1:1+2*pi,samples=200,smooth] plot (canvas polar cs:angle=\x r,radius={3*(6 + \x)});
\draw[dashed] (canvas polar cs:angle=1 r,radius={3*(6 + 1)}) -- (canvas polar cs:angle=1 r+ 2*pi r,radius={3*(6 + (1+2*pi))});
\fill (0,0) node[below] {$o$} circle (1pt);
\end{tikzpicture}
\subcaption{An open finite lightcone $C_f$.}
\end{subfigure}
\hfill
\begin{subfigure}[b]{0.45\textwidth}
\centering
\begin{tikzpicture}[scale=.8]
\draw[fill=gray!20,domain=1:1+2*pi,samples=200,smooth] plot (canvas polar cs:angle=\x r,radius={3*(6 + \x)});
\draw (canvas polar cs:angle=1 r,radius= {3*(6 + 1)}) -- (canvas polar cs:angle=1 r+ 2*pi r,radius={3*(6 + (1+2*pi))});
\fill (0,0) node[below] {$o$} circle (1pt);
\end{tikzpicture}
\subcaption{The closure $\overline{C_f}$ and the boundary $\partial C_f$.}
\end{subfigure}
\hfill
\caption{An example of an open finite lightcone $C_f$, its closure and boundary.}
\label{fig 9}
\end{figure}
\end{example}

\begin{example}\label{ex 2.5}
We consider two other examples which illustrate the shortcoming of the notions of the open and closed finite lightcones. See figure \ref{fig 10}.

\captionsetup[subfigure]{labelformat=simple}
\renewcommand\thesubfigure{\textit{\roman{subfigure}.}}
\begin{figure}[h]

\hfill
\begin{subfigure}[b]{0.45\textwidth}
\centering
\begin{tikzpicture}[scale=.8]
\draw[dashed,color=black,domain=1:1+2*pi,samples=200,smooth] plot (canvas polar cs:angle=\x r,radius={3*(6 + \x)});
\draw[dashed] (canvas polar cs:angle=1 r,radius= {3*(6 + 1 - 3)}) -- (canvas polar cs:angle=1 r+ 2*pi r,radius={3*(6 + (1+2*pi))})
(canvas polar cs:angle=3 r,radius= {3*(6 + 3 - 4.5)}) -- (canvas polar cs:angle=3 r,radius={3*(6 + 3)})
(canvas polar cs:angle=6 r,radius= {3*(6 + 6 - 5.5)}) -- (canvas polar cs:angle=6 r,radius={3*(6 + 6)});
\fill (0,0) node[below] {$o$} circle (1pt);
\end{tikzpicture}
\subcaption{An open finite lightcone $C_f$.}
\end{subfigure}
\hfill
\begin{subfigure}[b]{0.45\textwidth}
\centering
\begin{tikzpicture}[scale=.8]
\draw[fill=gray!20,domain=1:1+2*pi,samples=200,smooth] plot (canvas polar cs:angle=\x r,radius={3*(6 + \x)});
\draw  (canvas polar cs:angle=1 r,radius={3*(6 + 1)}) -- (canvas polar cs:angle=1 r+ 2*pi r,radius={3*(6 + (1+2*pi)+4)})
(canvas polar cs:angle=3 r,radius= {3*(6 + 3)}) -- (canvas polar cs:angle=3 r,radius={3*(6 + 3 + 8)})
(canvas polar cs:angle=6 r,radius= {3*(6 + 6)}) -- (canvas polar cs:angle=6 r,radius={3*(6 + 6 + 6)});
\fill (0,0) node[below] {$o$} circle (1pt);
\end{tikzpicture}
\subcaption{A closed finite lightcone $\calC_h$.}
\end{subfigure}
\hfill

\caption{Example \ref{ex 2.5}.}
\label{fig 10}
\end{figure}
\captionsetup[subfigure]{labelformat=parens}
\renewcommand\thesubfigure{\alph{subfigure}}

\begin{enumerate}[label=\roman*.]
\item Let $f$ be a positive bounded lower semicontinuous function which is discontinuous at a finite set of points $Z$ such that for all $\theta \in Z$,
\begin{align}
f(\vartheta) < f_{\inf,\circ} (\vartheta) \defeq \varliminf_{x\rightarrow \vartheta} f(x).
\label{eqn 2.3}
\end{align}
\item Similarly, let $h$ be a positive bounded upper semicontinuous function which is discontinuous at a finite set of points $Z$ such that for all $\theta \in Z$,
\begin{align}
h(\vartheta) > h^{\sup,\circ} (\vartheta) \defeq \varlimsup_{x\rightarrow \vartheta} h(x).
\label{eqn 2.4}
\end{align}
\end{enumerate}

Taking \textit{i} as the example, $C_f$ is an open finite lightcone and its closure $\overline{C_f}$ is the closed finite lightcone $\calC_{f_{\sup}}$. Then the interior $(\overline{C_f})^{\circ}$ of $\overline{C_f}$ is an open finite lightcone, which is larger than $C_f$, i.e. $C_f \subsetneq (\overline{C_f})^{\circ}$. See figure \ref{fig 11}.

\begin{figure}[h]

\hfill
\begin{subfigure}[b]{0.3\textwidth}
\centering
\begin{tikzpicture}[scale=.8]
\draw[dashed,color=black,domain=1:1+2*pi,samples=200,smooth] plot (canvas polar cs:angle=\x r,radius={3*(6 + \x)});
\draw[dashed] (canvas polar cs:angle=1 r,radius= {3*(6 + 1 - 3)}) -- (canvas polar cs:angle=1 r+ 2*pi r,radius={3*(6 + (1+2*pi))})
(canvas polar cs:angle=3 r,radius= {3*(6 + 3 - 4.5)}) -- (canvas polar cs:angle=3 r,radius={3*(6 + 3)})
(canvas polar cs:angle=6 r,radius= {3*(6 + 6 - 5.5)}) -- (canvas polar cs:angle=6 r,radius={3*(6 + 6)});
\fill (0,0) node[below] {$o$} circle (1pt);
\end{tikzpicture}
\subcaption{$C_f$ in \textit{i}.}
\end{subfigure}
\hfill
\begin{subfigure}[b]{0.3\textwidth}
\centering
\begin{tikzpicture}[scale=.8]
\draw[fill=gray!20,domain=1:1+2*pi,samples=200,smooth] plot (canvas polar cs:angle=\x r,radius={3*(6 + \x)});
\draw (canvas polar cs:angle=1 r,radius= {3*(6 + 1)}) -- (canvas polar cs:angle=1 r+ 2*pi r,radius={3*(6 + (1+2*pi))});
\fill (0,0) node[below] {$o$} circle (1pt);
\end{tikzpicture}
\subcaption{The closure $\overline{C_f}$.}
\end{subfigure}
\hfill
\begin{subfigure}[b]{0.3\textwidth}
\centering
\begin{tikzpicture}[scale=.8]
\draw[dashed,color=black,domain=1:1+2*pi,samples=200,smooth] plot (canvas polar cs:angle=\x r,radius={3*(6 + \x)});
\draw[dashed] (canvas polar cs:angle=1 r,radius={3*(6 + 1)}) -- (canvas polar cs:angle=1 r+ 2*pi r,radius={3*(6 + (1+2*pi))});
\fill (0,0) node[below] {$o$} circle (1pt);
\end{tikzpicture}
\subcaption{The interior $(\overline{C_f})^{\circ}$.}
\end{subfigure}
\hfill

\caption{The open finite lightcone $C_f$ in \textit{i}.}
\label{fig 11}
\end{figure}

Similarly for the closed finite lightcone $\calC_h$ in \textit{ii}., the interior $(\calC_h)^{\circ}$ is the open lightcone $C_{h_{\inf}}$, whose closure $\overline{(\calC_h)^{\circ}} = \overline{C_{h_{\inf}}}$ is a closed finite lightcone smaller than $\calC_h$, i.e. $\overline{(\calC_h)^{\circ}} \subsetneq \calC_h$. See figure \ref{fig 12}.

\begin{figure}[h]

\hfill
\begin{subfigure}[b]{0.3\textwidth}
\centering
\begin{tikzpicture}[scale=.8]
\draw[fill=gray!20,domain=1:1+2*pi,samples=200,smooth] plot (canvas polar cs:angle=\x r,radius={3*(6 + \x)});
\draw  (canvas polar cs:angle=1 r,radius={3*(6 + 1)}) -- (canvas polar cs:angle=1 r+ 2*pi r,radius={3*(6 + (1+2*pi)+4)})
(canvas polar cs:angle=3 r,radius= {3*(6 + 3)}) -- (canvas polar cs:angle=3 r,radius={3*(6 + 3 + 8)})
(canvas polar cs:angle=6 r,radius= {3*(6 + 6)}) -- (canvas polar cs:angle=6 r,radius={3*(6 + 6 + 6)});
\fill (0,0) node[below] {$o$} circle (1pt);
\end{tikzpicture}
\subcaption{$\calC_h$ in \textit{ii}.}
\end{subfigure}
\hfill
\begin{subfigure}[b]{0.3\textwidth}
\centering
\begin{tikzpicture}[scale=.8]
\draw[dashed,color=black,domain=1:1+2*pi,samples=200,smooth] plot (canvas polar cs:angle=\x r,radius={3*(6 + \x)});
\draw[dashed] (canvas polar cs:angle=1 r,radius={3*(6 + 1)}) -- (canvas polar cs:angle=1 r+ 2*pi r,radius={3*(6 + (1+2*pi))});
\fill (0,0) node[below] {$o$} circle (1pt);
\end{tikzpicture}
\subcaption{The interior $(\calC_h)^{\circ}$.}
\end{subfigure}
\hfill
\begin{subfigure}[b]{0.3\textwidth}
\centering
\begin{tikzpicture}[scale=.8]
\draw[fill=gray!20,domain=1:1+2*pi,samples=200,smooth] plot (canvas polar cs:angle=\x r,radius={3*(6 + \x)});
\draw (canvas polar cs:angle=1 r,radius= {3*(6 + 1)}) -- (canvas polar cs:angle=1 r+ 2*pi r,radius={3*(6 + (1+2*pi))});
\fill (0,0) node[below] {$o$} circle (1pt);
\end{tikzpicture}
\subcaption{The closure $\overline{(\calC_h)^{\circ}}$.}
\end{subfigure}
\hfill

\caption{The closed finite lightcone $\calC_h$ in \textit{ii}.}
\label{fig 12}
\end{figure}

A natural question raised by this example is what is the condition for $f$ and $h$ such that
\begin{align*}
\big(\overline{C_f}\big)^{\circ} = C_f,
\quad
\overline{(\calC_h)^{\circ}} = \calC_h.
\end{align*}
The answer to this question is the notion of a plump finite lightcone introduced in the next subsection.

\end{example}

\subsection{Plump finite lightcone}\label{sec 2.3}

In this subsection, we introduce a new notion of the so-called plump finite lightcone motivated by the question mentioned in the end of example \ref{ex 2.5}.
\begin{definition}\label{def 2.6}
A plump open finite lightcone in $C_0$ is an open finite lightcone $C_f$ satisfying
\begin{align}
(\overline{C_f})^{\circ} = C_f
\quad
\Leftrightarrow
\quad
(f^{\sup})_{\inf} = f.
\label{eqn 2.5}
\end{align}
A plump closed finite lightcone in $C_0$ is a closed finite lightcone $\calC_h$ satisfying
\begin{align}
\overline{(\calC_h)^{\circ}} = \calC_h
\quad
\Leftrightarrow
\quad
(h_{\inf})^{\sup} = h.
\label{eqn 2.6}
\end{align}
\end{definition}

By definition, we see that a plump open finite lightcone must be the interior of a closed finite lightcone, and a plump closed finite lightcone must be the closure of an open finite lightcone. We prove the converse of the above is also true, that a proposition saying that the closure of an open finite lightcone and the interior of a closed finite lightcone are plump.
\begin{proposition}\label{prop 2.7}
Let $C_f$ be an open finite lightcone, then its closure $\overline{C_f}$ is plump. Similarly let $\calC_h$  be a closed finite lightcone, then its interior $(\calC_h)^{\circ}$ is plump.
\end{proposition}
\begin{proof}
For the first part, by proposition \ref{prop 2.3.a}, let $f^{\sup}(\vartheta) =  \lim_{\delta \rightarrow 0^+} \sup_{B_{\vartheta}(\delta)}\{f(x)\}$, then
\begin{align*}
\overline{C_f} = \calC_{f^{\sup}}.
\end{align*}
Therefore
\begin{align*}
\overline{(\overline{C_f})^{\circ}} = \overline{C_f}
\end{align*}
follows from that $\overline{(\calC_{f^{\sup}})^{\circ}} \subset \calC_{f^{\sup}}$ and $(\overline{C_f})^{\circ} \supset C_f$.

For the second part, a similar argument shows that
$\Big(\overline{(\calC_h)^{\circ}}\Big)^{\circ}= (\calC_h)^{\circ}$.
\end{proof}

\begin{remark}\label{rem 2.8}
By propositions \ref{prop 2.3.a}, \ref{prop 2.3.b} and \ref{prop 2.7}, we have that for a lower semicontinuous function $f$ and an upper semicontinuous function $h$,
\begin{align*}
((f^{\sup})_{\inf})^{\sup} = f^{\sup},
\quad
((h_{\inf})^{\sup})_{\inf} = h_{\inf}.
\end{align*}
\end{remark}

\subsection{Perimeter of finite lightcone}\label{sec 2.4}
Let $C_f$ be an open finite lightcone. When $f$ is Lipschitz continuous, then $S_f$ has the tangent plane almost everywhere (in the sense of the usual measure of $\vartheta$ on $\mathbb{S}^{n-1}$ via the parameterisation as the graph of $f$). Let the metric restricted to $S_f$ be $\slashg_f$ and the corresponding volume form be $\dvol_f$, then
\begin{align*}
\slashg_f = f^2 \circg,
\quad
\dvol_f = f^{n-1} \dvol_{\circg},
\end{align*}
almost everywhere. Therefore the perimeter of $S_f$ is
\begin{align*}
|S_f| = \int_{S_f} 1 \cdot \dvol_f = \int_{\mathbb{S}^{n-1}} f^{n-1} \dvol_{\circg}.
\end{align*}

If $f$ is piecewise Lipschitz continuous but discontinuous, then by proposition \ref{prop 2.3.b},
\begin{align*}
\partial C_f = S_f \cup B_Z,
\end{align*}
where $Z$ is the set of discontinuous points of $f$ and $B_Z$ consists of the null segment $r\in [f(\vartheta), f^{\sup}(\vartheta)]$ at each discontinuous point $\vartheta \in Z$. Therefore if $Z$ is a Lipschitz set, then $B_Z$ is a null surface in the lightcone, thus $B_Z$ has zero perimeter. Thus for a closed finite lightcone $C_f$ where $f$ is piecewise Lipschitz continuous $h$ which is discontinuous at a Lipschitz set, the perimeter of the boundary $\partial C_f$ is also given by the integral $\int_{\mathbb{S}^{n-1}} f^{n-1} \dvol_{\circg}$.

Note that the above integral is always finite since $f$ is bounded, while the geometric meaning of the integral might be unclear since the usual notion of the perimeter of the boundary could be meaningless for a general $f$, where the tangent space of the lower envelope $S_f$ could not exist on a set with a positive $(n-1)$-dimensional measure. See example \ref{ex 2.9.a} in the following.
\addtocounter{theorem}{1}
\begin{subexample}\label{ex 2.9.a}
Let $N$ be a closed nowhere dense set with nonzero measure on the sphere and $\varphi_N$ be the characteristic function of $N$. Let $f = 2 - \varphi_N$, then $C_f$ is a closed finite lightcone. There exists no tangent space of the boundary of $C_f$ for $\vartheta\in N$, thus the usual notion of the perimeter of the boundary of $C_f$ is problematic, and the integral $\int_{\mathbb{S}^{n-1}} f^{n-1} \dvol_{\circg}$ does not measure the perimeter.
\end{subexample}

The above discussion of the perimeter of the boundary applies similarly to a closed finite lightcone. We have a similar example of a closed finite lightcone analogous to example \ref{ex 2.9.a}.
\begin{subexample}\label{ex 2.9.b}
Let $N$ be a closed nowhere dense set with nonzero measure on the sphere and $\varphi_N$ be the characteristic function of $N$. Let $h = 1 + \varphi_N$, then $\calC_h$ is a closed finite lightcone. There exists no tangent space of the boundary of $\calC_h$ for $\vartheta\in N$, thus the usual notion of the perimeter of the boundary of $\calC_h$ is problematic, and the integral $\int_{\mathbb{S}^{n-1}} h^{n-1} \dvol_{\circg}$ does not measure the perimeter.

Consider the interior of $\calC_h$, which is a plump open finite lightcone by proposition \ref{prop 2.7}. A straightforward calculation implies $h_{\inf} = 1$ then $(\calC_h)^{\circ}$ is simply $C_1$. In the next section \ref{sec 2.5}, we shall see that $\calC_h$ and $\overline{C_1}$ has the same domain of dependence, while the integral $\int_{\mathbb{S}^{n-1}} h^{n-1} \dvol_{\circg}$ is greater than $\int_{\mathbb{S}^{n-1}} 1\cdot \dvol_{\circg} = |\mathbb{S}^{n-1}|$.
\end{subexample}

Although the above examples show that the integral $\int_{\mathbb{S}^{n-1}} f^{n-1} \dvol_{\circg}$ may not measure the perimeter of a finite lightcone $C_f$, we show that this integral still has a geometric meaning.
\begin{lemma}\label{lem 2.10}
Let $C_f$ be an open finite lightcone and $\calC_h$ be a closed finite lightcone.
\begin{enumerate}[label=\alph*.]
\item
The integral $ \int_{\mathbb{S}^{n-1}} f^{n-1} \dvol_{\circg}$ is the supremum of the perimeters of all open finite lightcones contained in $C_f$ with a smooth boundary, i.e.
\begin{align*}
\int_{\mathbb{S}^{n-1}} f^{n-1} \dvol_{\circg}
=
\sup \{ |S_{\barf}| = |\partial C_{\barf}|: \barf \in C^{\infty}(\mathbb{S}^{n-1}), 0 < \barf \leq f \}.
\end{align*}

\item
The integral $ \int_{\mathbb{S}^{n-1}} h^{n-1} \dvol_{\circg}$ is the infimum of the perimeters of all open finite lightcones containing $\calC_h$ with a smooth boundary, i.e.
\begin{align*}
\int_{\mathbb{S}^{n-1}} h^{n-1} \dvol_{\circg}
=
\sup \{ |S_{\bar{h}}| = |\partial C_{\bar{h}}|: \bar{h} \in C^{\infty}(\mathbb{S}^{n-1}), h \leq \bar{h}\}.
\end{align*}
\end{enumerate}
\end{lemma}
\begin{proof}
The proof follows from Baire's theorem for semicontinuous functions (see \cite{Bai05}\cite{S81}\cite{E89}) and the dominated convergence theorem.
\end{proof}
By the above lemma, we introduce the generalised perimeter of the finite lightcone.
\begin{definition}\label{def 2.11}
Let $C_f$ be an open finite lightcone and $\calC_h$ be a closed finite lightcone, then the generalised perimeters of $C_f$ and $\calC_h$ are $\int_{\mathbb{S}^{n-1}} f^{n-1} \dvol_{\circg}$ and $\int_{\mathbb{S}^{n-1}} h^{n-1} \dvol_{\circg}$ respectively.
\end{definition}

\subsection{Domain of dependence of finite lightcone}\label{sec 2.5}

We recall the notion of the domain of dependence of a finite lightcone in $\mathbb{M}^{n,1}$. 
\begin{definition}\label{prop 2.12}
\begin{enumerate}[label=\alph*.]
\item
Let $C_f$ be an open finite lightcone. Define the domain of dependence of $C_f$ as the set of points, through which any inextensible causal curve intersects with $C_f$. Denote the domain of dependence of $C_f$ by $\calD(C_f)$.
\item
For a closed finite lightcone $\calC_h$, define its domain of dependence the same way as above and denote by $\calD(\calC_h)$.
\item
For a general achronal set $A$, i.e. $I^+(x) \cap A=\emptyset$ for any $x \in A$, its domain of dependence is also defined the same way as above.
\end{enumerate}
\end{definition}

It is easy to show the following point-set topology result for the domain of dependence of a finite lightcone.
\begin{proposition}\label{prop 2.13}
\begin{enumerate}[label=\alph*.]
\item
The domain of dependence $\calD(C_f)$ of an open finite lightcone $C_f$ is open in the causal future $J^+(o)$ of $o$.
\item
The domain of dependence $\calD(\calC_h)$ of a closed finite lightcone $\calC_h$ is closed in $J^+(o)$.
\end{enumerate}
\end{proposition}
\begin{proof}
\begin{enumerate}[label=\alph*.]
\item
If $x\in I^+(o) \cap \calD(C_f)$, then $\uC_x \cap C_0 \subset C_f$. Since $C_f$ is an open finite lightcone, there exist a future timelike vector $v$ and a sufficiently small $\delta>0$, such that $\uC_{x+\delta v} \cap C_0  \subset C_f$. Therefore $x \in J^+(o) \cap I^-(x+\delta v) \subset J^+(o) \cap \calD(C_f)$ implies that $x$ is an interior point of $\calD(C_f)$ in $J^+(o)$.

If $x \in C_f$, since $C_f$ is an open finite lightcone, there exist a future timelike vector $v$ and a sufficiently small $\delta>0$, such that $\uC_{x+\delta v} \cap C_0  \subset C_f$. Therefore $x \in J^+(o) \cap I^-(x+\delta v) \subset J^+(o) \cap \calD(C_f)$ implies that $x$ is an interior point of $\calD(C_f)$ in $J^+(o)$.

\item
We show that $J^+(o) \cap \overline{\calD(\calC_h)} = J^+(o) \cap \calD(\calC_h)$.
Let $\{x_n\}$ be a sequence in $J^+(o) \cap \calD(\calC_h)$ and $x_n \overset{n\rightarrow + \infty}{\longrightarrow} x$. Since $J^+(o)$ is closed, then $x\in J^+(o)$. We need to show $x \in J^+(o) \cap \calD(\calC_h)$.

Without loss of generality, since $J^+(o)= I^+(o) \cup C_0$, we need to simply consider two cases: either $\{x_n\} \subset I^+(o) \cap \calD(\calC_h)$ or $\{x_n\} \subset C_0 \cap \calD(\calC_h) = \calC_h$. The later case is simple, since $\calC_h$ is a closed finite lightcone.
For the former case, we have that $J^-(x_n) \cap C_0 \subset \calC_h$, then $J^-(x) \cap C_0 \subset \calC_h$ since $\calC_h$ is a closed finite lightcone, thus $x\in J^+(o) \cap \calD(\calC_h)$.

\end{enumerate}
\end{proof}

We prove a proposition comparing the domains of dependence of a closed finite lightcone $\calC_h$ and the closure of its interior $\overline{(\calC_h)^{\circ}}$, which verifies the claim in example \ref{ex 2.9.b}.
\begin{proposition}\label{prop 2.14}
Let $\calC_h$ be a closed finite lightcone and $I^+(o)$ be the timelike future of the origin $o$, then
\begin{align*}
I^+(o) \cap\calD(\calC_h) = I^+(o) \cap \calD\Big(\overline{(\calC_h)^{\circ}}\Big).
\end{align*}
\end{proposition}
\begin{proof}
Since $\overline{(\calC_h)^{\circ}} \subset \calC_h \Rightarrow \calD\big(\overline{(\calC_h)^{\circ}}\big) \subset \calD(\calC_h)$, it is sufficient to show the inclusion relation in the reversed direction.

If $x \in I^+(o) \cap \calD(\calC_h)$ and let $\uC_x$ be the past lightcone emanating from $x$, then $\uC_x\cap C_0 \subset \calC_h$, thus $\uC_x\cap C_0 \subset \overline{(\calC_h)^{\circ}}$. Hence $x \in \calD\Big(\overline{(\calC_h)^{\circ}}\Big)$. The proposition is proved.
\end{proof}

We have another similar property on the domain of dependence of a closed finite lightcone.
\begin{proposition}\label{prop 2.15}
Let $\calC_h$ be a closed finite lightcone and $(\calC_h)^{\circ}$ be the interior of $\calC_h$. We have that
\begin{align*}
I^+(o) \cap \calD(\calC_h) = I^+(o) \cap \overline{\calD((\calC_h)^{\circ})}.
\end{align*}
\end{proposition}
\begin{proof}
Since $(\calC_h)^{\circ} \subset \calC_h$, we have $\overline{\calD((\calC_h)^{\circ})} \subset \overline{\calD(\calC_h)}$, thus $I^+(o) \cap \overline{\calD((\calC_h)^{\circ})} \subset I^+(o) \cap \overline{\calD(\calC_h)} = I^+(o) \cap \calD(\calC_h)$, where the last equality follows from proposition \ref{prop 2.13} that $\calD(\calC_h)$ is closed in $J^+(o)$ thus it is also closed in $I^+(o)$.

We prove the inclusion relation in the other direction. Suppose that $x\in I^+(o) \cap \calD(\calC_h)$, then there exists a future timelike vector $v$ such that $x-\delta v \in I^+(o) \cap \calD(\calC_h)$ for all $\delta \in (0,1)$. Thus $\uC_{x-\delta v} \cap C_0 \subset (\calC_h)^{\circ}$, which implies $x-\delta v \in I^+(o) \cap \calD((\calC_h)^{\circ})$. Therefore $x\in \overline{\calD((\calC_h)^{\circ})}$.
\end{proof}

\section{Lorentz reflection and polarisation}\label{sec 3}
We introduce the Lorentz reflection and polarisation in the Minkowski spacetime, and prove several properties which will be useful for studying the isoperimetric inequality for the domain of dependence of a finite lightcone.

\subsection{Definition}\label{sec 3.1}

Let $H$ be a timelike hyperplane in $\mathbb{M}^{n,1}$ passing through the origin $o$. We define the Lorentz reflection and polarisation about $H$ as follows.
\begin{definition}[Lorentz reflection]\label{def 3.1}
Let $w$ be the unit spacelike normal vector of $H$. Denote the reflection about $H$ by $\gamma_H$, then
\begin{align*}
\gamma_H (x) = x - 2 \langle x, w \rangle w.
\end{align*}
Let $U$ be a set in $\mathbb{M}^{n,1}$, then the reflection of $U$ about $H$ is $\gamma_H (U)$. Let $f$ be a function on $\mathbb{M}^{n,1}$. Denote the reflection of $f$ about $H$ by $\gamma_H[f]$, then
\begin{align*}
\gamma_H [f](x) = f(\gamma_H (x)).
\end{align*}
\end{definition}
\begin{definition}[Lorentz polarisation]\label{def 3.2}
Let $v$ be a timelike vector at the origin and $v \notin H$. Then $H$ divides $\mathbb{M}^{n,1}$ into two open half-spacetimes $H_+$ and $H_-$, such that $v\in H_+$. 
\begin{enumerate}[label=\alph*.]
\item
Let $U$ be a set in $\mathbb{M}^{n,1}$, then we define the polarisation $U^{\gamma_H,v}$ of $U$ about the pair $(H,v)$ in the following:
\begin{align*}
\left\{
\begin{aligned}
&
U^{\gamma_H,v} \cap H_+ = (U \cup \gamma_H(U) ) \cap H_+,
\\
&
U^{\gamma_H,v} \cap H_- = (U \cap \gamma_H(U) ) \cap H_-,
\\
&
U^{\gamma_H,v} \cap H = U \cap H.
\end{aligned}
\right.
\end{align*}
Use $\rho_H^v$ to denote the map of the Lorentz polarisation about $(H,v)$, i.e. $\rho_H^v (U) = U^{\gamma_H,v}.$

\item
Let $f$ be a function on $\mathbb{M}^{n,1}$. Define the polarisation $f^{\gamma_H, v}$ of $f$ about the pair $(H,v)$ as
\begin{align*}
f^{\gamma_H,v}(x)=
\left\{
\begin{aligned}
& 
\max\{ f(x), f(\gamma_H(x)) \}, 
\quad 
x \in H_+,
\\
& 
\min\{ f(x), f(\gamma_H(x)) \}, 
\quad 
x \in H_-,
\\
& 
f(x)
\quad 
x \in H.
\end{aligned}
\right.
\end{align*}
Also use $\rho_H^v$ to denote the map of the Lorentz polarisation about $(H,v)$, i.e. $\rho_H^v [f] = f^{\gamma_H,v}$.
\end{enumerate}
\end{definition}

By the definition of the Lorentz polarisation, we have that it preserves the volume of the set.
\begin{proposition}\label{prop 3.3}
Let $U \in \mathbb{M}^{n,1}$, then the polarisation of $U$ has the same volume as $U$, i.e.
\begin{align*}
|U| = | U^{\gamma_H,v}|.
\end{align*}
\end{proposition}

Moreover we have the following proposition on the Lorentz polarisation and the symmetric difference volume.
\begin{proposition}\label{prop 3.4}
The Lorentz polarisation is a distance nonincreasing contraction map for open sets with respect to the symmetric difference volume. More concretely, let $H$ be a timelike hypersurface and $v\notin H$ a timelike vector, we have that
\begin{align*}
|E_1^{\gamma_H,v} \triangle E_2^{\gamma_H,v}| \leq |E_1 \triangle E_2|.
\end{align*}
\end{proposition}
\begin{proof}
The symmetric difference volumes have the following formulae
\begin{align*}
&
|E_1^{\gamma_H,v} \triangle E_2^{\gamma_H,v}| 
= 
|E_1^{\gamma_H,v}| + |E_2^{\gamma_H,v}| - 2 | E_1^{\gamma_H,v} \cap E_2^{\gamma_H,v}|,
\\
&
|E_1 \triangle E_2|
=
|E_1| + |E_2| - 2 |E_1 \cap E_2|.
\end{align*}
Since $E_1 \cap E_2 \subset E_1$, then $(E_1 \cap E_2 )^{\gamma_H,v} \subset E_1^{\gamma_H,v}$. Similarly $(E_1 \cap E_2 )^{\gamma_H,v} \subset E_2^{\gamma_H,v}$, thus $(E_1 \cap E_2 )^{\gamma_H,v} \subset E_1^{\gamma_H,v} \cap E_2^{\gamma_H,v}$, which implies
\begin{align*}
| E_1 \cap E_2 |
\leq
| (E_1 \cap E_2 )^{\gamma_H,v} |
\leq
| E_1^{\gamma_H,v} \cap E_2^{\gamma_H,v} |.
\end{align*}
The proposition follows.
\end{proof}

It is easy to derive the following property between the polarisations of sets and functions.
\begin{proposition}\label{prop 3.5}
For $t\in \mathbb{R}$, let $L_f(t)$ be the upper level set of $f$ that
\begin{align*}
L_f(t) = \{x: f(x) >t \}.
\end{align*}
Then the polarisation of the level set $L_f(t)$ is equal to the level set of the polarisation of $f$,
\begin{align*}
(L_f(t))^{\gamma_H,v} = L_{f^{\gamma_H,v}}(t).
\end{align*}
\end{proposition}

\subsection{Lorentz reflection and polarisation of finite lightcone and causal diamond}
We apply the Lorentz reflection and polarisation to the finite lightcone as an example. 
\begin{example}\label{ex 3.6}
Let $H$ be a timelike hyperplane in $\mathbb{M}^{n,1}$ and $v\notin H$ be a timelike vector. Let $\gamma_H$ be the reflection about $H$. Suppose that $C_f$ is an open finite lightcone.
\begin{enumerate}[label=\alph*.]
\item
The reflection of $C_f$ about $H$ is also an open finite lightcone $\gamma_H(C_f)$. Suppose that $\gamma_H(C_f)$ is also parameterised by some function on the sphere. The Lorentz reflection $\gamma_H$ in fact defines a transformation between the functions on the sphere, which transforms $f$ to the parameterisation function of the finite lightcone $\gamma_H(C_f)$. This transformation is actually related to the conformal reflection on the sphere induced by the Lorentz reflection $\gamma_H$.

\item
Similarly as the reflection of $C_f$ about $H$, the polarisation of $C_f$ about the pair $(H,v)$ is also an open finite lightcone $(C_f)^{\gamma_H,v}$. Again $(C_f)^{\gamma_H,v}$ is parametrised by some function on the sphere, thus it introduces the transformation from a function $f$ to the parameterisation function of $(C_f)^{\gamma_H,v}$.
\end{enumerate}
\end{example}

We summarise the above transformations between the parameterisation functions of finite open lightcones in the following definition.
\begin{definition}\label{def 3.7}
Let $H$ be a timelike hyperplane in $\mathbb{M}^{n,1}$ and $v\notin H$ be a timelike vector. Let $\gamma_H$ be the reflection about $H$.
\begin{enumerate}[label=\alph*.]
\item 
Define the conformal reflection $\gamma_{c,H}$ on the sphere induced by the Lorentz reflection $\gamma_H$ as follows. Suppose $f$ is the parameterisation function of an open finite lightcone $C_f$. Then define $\gamma_{c,H}[f]$ as the parameterisation function of the Lorentz reflection $\gamma_H (C_f)$ of $C_f$.

\item 
Define the conformal polarisation $\rho_{c,H}^v$ on the sphere induced by the pair $(H,v)$ as follows. Suppose $f$ is the parameterisation function of an open finite lightcone $C_f$. Then define $\gamma_{c,H}[f]$ as the parameterisation function of the Lorentz polarisation $(C_f)^{\gamma_H,v}$ of $C_f$.
\end{enumerate}
In the above notations $\gamma_{c,H}$ and $\rho_{c,H}^v$, the subscript $c$ represents the word ``conformal''.
\end{definition}

As the Lorentz polarisation preserving the volume, the Lorentz polarisation similarly preserves the generalised perimeter of the finite lightcone.
\begin{proposition}\label{prop 3.8}
Let $C_f$ be an open finite lightcone and $\rho_{c,H}^v[f]$ be the parameterisation function of the polarisation $(C_f)^{\gamma_H,v}$ of $C_f$. We have that
\begin{align*}
\int_{\mathbb{S}^{n-1}} f^{n-1} \dvol_{\circg} = \int_{\mathbb{S}^{n-1}} (\rho_{c,H}^v[f])^{n-1} \dvol_{\circg}.
\end{align*}
\end{proposition}

\begin{example}\label{ex 3.9}
Let $q$ be a point in the future of the origin $o$. Let $J^+(o)$ be the causal future of $o$ and $J^-(q)$ be the causal past of $q$, then the causal diamond $J(o,q)$ with vertices $o$ and $q$ is defined as the intersection of $J^+(o)$ and $J^-(q)$, i.e.
\begin{align*}
J(o,q) = J^+(o) \cap J^-(q).
\end{align*}
Let $H$ and $v$ be as in example \ref{ex 3.6}, then the reflection and polarisation of the causal diamond are given by
\begin{align*}
\gamma_H(J(o,q)) = J(o,\gamma_H(q)),
\quad
(J(o,q))^{\gamma_H,v} = J(o, q^{\gamma_H,v})
=
\left\{
\begin{aligned}
&
J(o,q),
\quad
q\in H_+ \cup H
\\
&
J(o,\gamma_H(q))
\quad
q\in H_-.
\end{aligned}
\right.
\end{align*}
See figure \ref{fig 13}.
\begin{figure}[h]
\centering
\begin{tikzpicture}[scale=.8]
\draw (-8 *.7,8 *.7) -- (0,0) node[below] {$o$} -- (8 *.7,8 *.7);
\draw (-2 *.7,2 *.7) -- (4 *.7,8 *.7) node[above] {$q$} -- (6 *.7,6 *.7);
\path[fill=gray,opacity=0.3] 
(0,0) -- (-2 *.7,2 *.7) -- (4 *.7,8 *.7) -- (6 *.7,6 *.7) circle;
\node at (4 *.7,6 *.7) {$J(o,q)$};
\draw (-4 *.7,4 *.7) -- (-1 *.7,7 *.7) node[above] {$\gamma_H(q)$} -- (3 *.7,3 *.7);
\path[fill=blue,opacity=0.2] 
(0,0) -- (-4 *.7,4 *.7) -- (-1 *.7,7 *.7) -- (3 *.7,3 *.7) circle;
\node at (-1 *.7,5 *.7) {$(J(o,q))^{\gamma_{H,v}}$};
\node[above] at (4 *.7+1,8 *.7) {$H_{-}$};
\node[above] at (-1 *.7-2,7 *.7) {$H_{+}$};
\draw (0,0) -- (1 *.7,5 *.7) -- (1.6 *.7,8 *.7) node[right] {$H$};
\end{tikzpicture}
\caption{Example \ref{ex 3.9}: $q \in H_-$.}
\label{fig 13}
\end{figure}
\end{example}

\subsection{Monotonicity for domain of dependence under Lorentz polarisation}\label{sec 3.2}
We prove that the domain of dependence of the finite lightcone satisfies a useful monotone property under the Lorentz polarisation.
\begin{proposition}\label{prop 3.10}
Let $H$ be a timelike hyperplane in $\mathbb{M}^{n,1}$ and $v\notin H$ be a timelike vector. Let $\gamma$ be the reflection about $H$.\footnote{We drop the subscript $H$ when $H$ is understood.} Suppose that $C_f$ is an open finite lightcone, and $\calD(C_f)$ is the domain of dependence of $C_f$. We have that the Lorentz polarisation of $\calD(C_f)$ is contained in the domain of dependence of the Lorentz polarisation of $C_f$, i.e.
\begin{align*}
(\calD(C_f))^{\gamma,v} \subset \calD((C_f)^{\gamma,v}).
\end{align*}
\end{proposition}
\begin{proof}
Suppose $q\in(\calD(C_f))^{\gamma,v}$, there are three cases $q \in (\calD(C_f))^{\gamma,v} \cap H_-$, $q \in (\calD(C_f))^{\gamma,v} \cap H_+$ and $q \in (\calD(C_f))^{\gamma,v} \cap H$
\begin{enumerate}[label=\textit{\roman*.}]
\item 
If $q \in (\calD(C_f))^{\gamma,v} \cap H_-$, then $q \in \calD(C_f) \cap H_-$ and $\gamma(q) \in \calD(C_f) \cap H_+$, hence $J(o,q) \subset \calD(C_f)$ and $J(o,\gamma(q)) \subset \calD(C_f)$. Thus $J(o,q) \cap C_0 \subset C_f$ and $J(o,\gamma(q)) \cap C_0 \subset C_f$, which implies $J(o,q) \cap C_0 \subset (C_f)^{\gamma,v}$. Therefore $q \in \calD((C_f)^{\gamma,v})$.

\item 
If $q \in (\calD(C_f))^{\gamma,v} \cap H_+$, then $q \in \calD(C_f) \cap H_+$ or $\gamma(q) \in \calD(C_f) \cap H_-$, hence $J(o,q) \subset \calD(C_f)$ or $J(o,\gamma(q)) \subset \calD(C_f)$. Thus $J(o,q) \cap C_0 \subset C_f$ or $J(o,\gamma(q)) \cap C_0 \subset C_f$, which implies $J(o,q) \cap C_0 \subset (C_f)^{\gamma,v}$. Therefore $q \in \calD((C_f)^{\gamma,v})$.

\item
If $q \in (\calD(C_f))^{\gamma,v} \cap H$, then $q \in \calD(C_f) \cap H$, hence $J(o,q) \subset \calD(C_f)$. Thus $J(o,q) \cap C_0 \subset C_f$ which implies $J(o,q) \cap C_0 = (J(o,q) \cap C_0)^{\gamma,v} \subset (C_f)^{\gamma,v}$. Therefore $q \in \calD((C_f)^{\gamma,v})$.
\end{enumerate}
Summary of the conclusions in two cases proves the proposition.
\end{proof}

As a corollary of proposition \ref{prop 3.10}, we have the monotonicity of the volume of the domain of dependence of the finite lightcone under Lorentz polarisation.
\begin{corollary}\label{coro 3.11}
The volume of the domain of dependence $\calD(C_f)$ of an open finite lightcone $C_f$ is no more than the volume of the domain of dependence $\calD((C_f)^{\gamma,v})$ of the polarisation of $C_f$, i.e.
\begin{align*}
|\calD(C_f)| \leq |\calD((C_f)^{\gamma,v})|.
\end{align*}
\end{corollary}

\section{Variational problem of the isoperimetric inequality}\label{sec 4}
We formulate again the variational problem of the isoperimetric inequality.
\begin{problem}\label{prob 4.1}
Define the following set $F_1$ of the open finite lightcone with the perimeter $n\omega_n$, the area of $(\mathbb{S}^{n-1},\circg)$, or equivalently with the parameterisation function of the $L^{n-1}$ norm being $n\omega_n$:
\begin{align*}
F_1 = \{C_f: \int_{\mathbb{S}^{n-1}} f^{n-1} \dvol_{\circg} = n\omega_n \}.
\end{align*}
Introduce the volume functional $\calV$ of the domain of dependence $\calD(C_f)$ of the open finite lightcone $C_f \in F_1$,
\begin{align*}
\calV(C_f) = |\calD(C_f)|.
\end{align*}
Consider the variational problem of the volume functional $\calV$ in the set $F_1$. 
\begin{enumerate}[label=\alph*.]
\item
Find the supremum $\sup_{F_1} \calV$ and determine whether the supremum can be achieved. 
\item
If the supremum can be achieved, find the extremal point of $\calV$ in $F_1$.
\end{enumerate}
\end{problem}
In this section, we solve part \textit{a.} of the above variational problem.

\subsection{Precompactness of set of finite lightcones constructed by Lorentz polarisation}\label{sec 4.1}
Let $C_f$ be an open finite lightcone. Proposition \ref{prop 3.8} and corollary \ref{coro 3.11} implies that the Lorentz polarisation could increase the volume of the domain of dependence while preserving the $L^{n-1}$ norm of the parameterisation function. Thus we introduce the following set of the open finite lightcone by the Lorentz polarisation.
\begin{definition}
Let $v$ be a timelike vector at the origin $o$. Let $C_f$ be an open finite lighcone. Define the following set $P_f^v$ of the open finite lightcone obtained by applying the Lorentz polarisation for finite times:
\begin{align*}
P_f^v
=
\{ (C_f)^{\gamma_1, \gamma_2,\cdots, \gamma_k, v}: k\in \mathbb{N}, v\notin H_1, \cdots, v \notin H_k\},
\end{align*}
where $(C_f)^{\gamma_1, \gamma_2,\cdots, \gamma_k, v}$ is $(\cdots((C_f)^{\gamma_1,v})^{\gamma_2,v}\cdots)^{\gamma_k,v}$.
\end{definition}

We show that the finite lightcones in $P_f^v$ are uniformly bounded.
\begin{lemma}\label{lem 4.3}
There exists a positive constant $R$ depending on $f$ and $v$ such that for any open finite lightcone $C_{\barf} \in P_f^v$, we have
\begin{align*}
\barf \leq R.
\end{align*}
\end{lemma}
\begin{proof}
Without loss of generality, we can assume that $v= e_0=(1,0,\cdots,0)$. Assume that $f\leq R$, then for any $C_{\barf} \in P_f^v$, we have
\begin{align*}
\barf \leq R,
\end{align*}
since $(C_R)^{\gamma,e_0} = C_R$ implies that $C_{\barf} \subset C_R \Rightarrow (C_{\barf})^{\gamma,e_0} \subset C_R$ for any reflection $\gamma$ with $e_0 \neq \gamma(e_0)$. See figure \ref{fig 14}.
\begin{figure}[h]
\centering
\begin{tikzpicture}[scale=.7]
\draw (-8 *.7,8 *.7) -- (0,0) -- (8 *.7,8 *.7);
\node[below] at (0,-0.075) {$o$};
\draw[->] (0,0) -- (0,1.5) node[above] {$e_0$};
\node[above] at (3 *.7,6 *.7) {$H_{-}$};
\node[above] at (-1 *.7,6 *.7) {$H_{+}$};
\draw (0,0) -- (1 *.7,5 *.7) -- (1.6 *.7,8 *.7) node[right] {$H$};
\draw[red,thick] (-2 *.7 +.05*1.5 ,2 *.7 +.05*1.5 ) -- (0,.1*1.5 ) -- (6 *.7 -.05*1.5 ,6 *.7 +.05*1.5 ) node[left] {$C_{\barf}$};
\draw[violet,thick] (-4 *.7 +.1*1.5 , 4 *.7 +.1*1.5 ) node[right] {$(C_{\barf})^{\gamma,e_0}$} -- (0,.2*1.5 ) -- (3 *.7 -.1*1.5 ,3 *.7 +.1*1.5 ) ;
\draw[blue,thick] (-7 *.7 -.05*1.5 , 7 *.7 -.05*1.5 ) -- (0,-.1*1.5 ) -- (7 *.7 +.05*1.5 , 7 *.7 -.05*1.5 ) node[below right] {$C_R$};
\end{tikzpicture}
\caption{$C_{\barf} \subset C_R \Rightarrow (C_{\barf})^{\gamma,e_0} \subset C_R$.}
\label{fig 14}
\end{figure}
\end{proof}

Before proving the precompactness of the set $P_f^v$, we introduce the strict future boundary of the domain of dependence of an open finite lightcone.
\begin{definition}\label{def 4.4}
\begin{enumerate}[label=\alph*]
\item
Let $C_f$ be an open finite lightcone and $\calD(C_f)$ be the domain of dependence of $C_f$. Define the strict future boundary of $\calD(C_f)$, denoted by $\partial_+ \calD(C_f)$, as the following set
\begin{align*}
\partial_+ \calD(C_f) 
= 
\{x \in \partial \calD(C_f): I^+(x) \cap \calD(C_f) = \emptyset , I^-(x) \cap \calD(C_f) \neq \emptyset \},
\end{align*}
where $I^+(x)$ is the chronological future of $x$, the set of points which can be reached by a future-directed timelike curve emanating from $x$.

\item
More generally, for a causally convex set $\calE$, define its strict future boundary $\partial_+ \calE$ and strict past boundary $\partial_- \calE$ by
\begin{align*}
&
\partial_+ \calE
= 
\{x \in \partial \calE: I^+(x) \cap \calE = \emptyset , I^-(x) \cap \calE \neq \emptyset \},
\\
&
\partial_- \calE
= 
\{x \in \partial \calE: I^-(x) \cap \calE = \emptyset , I^+(x) \cap \calE \neq \emptyset \}.
\end{align*}

\item
For a causally convex set $\calE$, define its future boundary $\bar{\partial}_+ \calE$ and past boundary $\bar{\partial}_- \calE$ by
\begin{align*}
&
\bar{\partial}_+ \calE
= 
\{x \in \partial \calE: I^+(x) \cap \calE = \emptyset \},
\\
&
\bar{\partial}_- \calE
= 
\{x \in \partial \calE: I^-(x) \cap \calE = \emptyset \}.
\end{align*}

\end{enumerate}

\end{definition}

We can prove the precompactness of the set $P_f^v$ now.
\begin{lemma}\label{lem 4.5}
Let $\{C_{\barf_k}\}_{k\in \mathbb{N}}$ be a sequence of open finite lightcones in $P_f^v$, then there exist a subsequence $\{C_{\barf_{k_i}}\}_{i\in \mathbb{N}}$ and an open finite lightcone $C_{\barf}$ such that
\begin{align*}
\lim_{i\rightarrow +\infty} | \calD(C_{\barf_{k_i}}) \triangle \calD(C_{\barf}) | = 0,
\quad
|\calD(C_{\barf})| = \lim_{i\rightarrow +\infty} | \calD(C_{\barf_{k_i}})|.
\end{align*}
\end{lemma}
\begin{proof}
We parameterise the future and strict past boundary of $\calD(C_{\barf_k})$ and prove the convergence of the parameterisation functions. Let $\{t, x_1, \cdots, x_n\}$ be the rectangular coordinate system of $\mathbb{M}^{n,1}$. By lemma \ref{lem 4.3}, suppose that $C_{\barf_k} \subset C_R$. Then for each $\calD(C_{\barf_k})$, we construct the function $u_k$ on the ball $B_R \subset \mathbb{R}^n$, such that
\begin{align*}
\calD(C_{\barf_k}) = \{ (t,x) \in \mathbb{M}^{n,1}:  r\leq t< u_k(x) \},
\quad
r=\sqrt{(x_1)^2 + \cdots + (x_n)^2}.
\end{align*}
We observe that $u_k$ is simply the parameterisation function of the strict future boundary $\partial_+ \calD(C_{\barf_k})$ of $\calD(C_{\barf_k})$ with the continuous extension by the function $r$ at the boundary $\partial(C_{\barf_k})$. Since $\{u_k\}_{k\in \mathbb{N}}$ is uniformly bounded and uniformly Lipschitz with the Lipschitz constant $1$, there exists a subsequence $\{u_{k_i}\}_{i\in \mathbb{N}}$ converging uniformly to a Lipschitz function $u$ with the Lipschitz constant $1$ and $u|_{\partial B_R} = R$. Define the set $\calD$ by
\begin{align*}
\calD = \{ (t,x) \in \mathbb{M}^{n,1}:  r\leq t< u(x) \}.
\end{align*}
We show that $\calD$ is the domain of dependence of some open finite lightcone which will be determined in the proof.
\begin{enumerate}[label=\textit{\roman*}.]
\item
\underline{\textsc{Claim:}} if $q \in \calD$, then $J(o,q) \subset \calD$.

Suppose that $q=(t_q,x_q)$. Since $q \in \calD$, then $ r(q) \leq t_q < u(x_q)$. Thus there exists a small positive number $\delta$ such that for $k_i$ sufficiently large, $r(q) \leq t_q + \delta < u_{k_i}(x_q)$. Define $q_{+\delta} = (t_q + \delta,x_q)$, then
\begin{align*}
J(o,q_{+\delta}) \subset \calD(C_{\barf_{k_i}}) = \{ (t,x) \in \mathbb{M}^{n,1}:  r\leq t < u_{k_i}(x) \}.
\end{align*}
Taking the limit $i\rightarrow +\infty$, we have
\begin{align*}
J(o,q_{+\delta}) \subset \overline{\calD} = \{ (t,x) \in \mathbb{M}^{n,1}:  r\leq t \leq u(x) \},
\end{align*}
therefore
\begin{align*}
J(o,q) \subset \calD = \{ (t,x) \in \mathbb{M}^{n,1}:  r\leq t < u(x) \},
\end{align*}

\item
Define $C_{\barf}$ as the strict past boundary of $\calD$, which is $\calD \cap C_0$. \textit{i.} implies that $\calD$ is open in $J^+(o)$, thus $\calD \cap C_0$ is an open finite lightcone, hence $C_{\barf}$ is well defined.

\item
\underline{\textsc{Claim:}} if $J(o,q) \cap C_0 \subset C_{\barf}$, then $q \in \calD$.

Note $J(o,q) \cap C_0$ is a closed set, thus $J(o,q) \cap C_0 \subset C_{\barf} \subset \calD$ implies that there exists a small positive $\delta$ such that  for sufficiently large $k_i$
\begin{align*}
J(o,q_{+\delta}) \cap C_0 \subset \calD(C_{\barf_{k_i}}) = \{ (t,x) \in \mathbb{M}^{n,1}:  r\leq t < u_{k_i}(x) \}.
\end{align*}
Hence $q_{+\delta} \in \overline{\calD}$. Therefore $q \in \calD$.
\end{enumerate}
\textit{i. ii.} \& \textit{iii.} imply that $\calD = \calD(C_{\barf})$. The lemma is proved.
\end{proof}

\subsection{Existence of extremal finite lightcone maximising volume of domain of dependence}\label{sec 4.2}
We study the volume functional $\calV$ of the domain of dependence in the set $P_f^v$. Introduce the closure of the set $P_f^v$ under the distance of the volume of the symmetric difference between the domain of dependences.
\begin{definition}\label{def 4.6}
Introduce the distance $d_{\calD}$ between two open finite lightcones by
\begin{align*}
d_{\calD} ( C_{f_1}, C_{f_2} )= | \calD(C_{f_1}) \triangle \calD(C_{f_2}) |
\end{align*}
Define the set $\overline{P_f^v}$ as the closure of $P_f^v$ under the distance $d_{\calD}$, i.e. the set of the open finite lightcones $C_{\barf}$ that there exists a sequence $\{C_{\barf_k}\}_{k \in \mathbb{N}}$ in $P_f^v$ such that
\begin{align*}
\lim_{k\rightarrow +\infty} d_{\calD}( C_{\barf_k}, C_{\barf} )
= 
\lim_{k\rightarrow +\infty} | \calD(C_{\barf_k}) \triangle \calD(C_{\barf}) | 
=0.
\end{align*}
\end{definition}

By lemma \ref{lem 4.5}, we know that $\overline{P_f^v}$ is complete under the distance $d_{\calD}$.
\begin{proposition}\label{prop 4.7}
$(\overline{P_f^v}, d_{\calD})$ is a complete compact metric space.
\end{proposition}
\begin{proof}
The completeness and compactness follow from lemma \ref{lem 4.5}.
\end{proof}

Moreover we can show that $\overline{P_f^v}$ is invariant under the Lorentz polarisation. 
\begin{proposition}\label{prop 4.8}
$\overline{P_f^v}$ is invariant under the Lorentz polarisation about any pair $(H,v)$ where $H$ is a timelike hyperplane through the origin $o$ not containing $v$.
\end{proposition}
\begin{proof}
In order to prove the proposition, we need to show that $C_{\barf} \in \overline{P^v_f} \Rightarrow (C_{\barf})^{\gamma,v} \in  \overline{P^v_f}$ where $\gamma$ is the reflection about $H$. We prove this by three steps.

\vspace{5pt}
\noindent{\underline{\textsc{Step 1.}} }
Let $C_{\barf} \in \overline{P^v_f}$ be the limit of the sequence $\{C_{f_k}\}_{k\in \mathbb{N}}$ under the distance $d_{\calD}$. We obtain a description of $(C_{\barf})^{\gamma, v}$ as a certain limit of $(C_{f_k})^{\gamma, v}$.

Introduce $(\calD(C_{f_k}))^{\gamma,v}$, the polarisation of $\calD(C_{f_k})$, where $\gamma_H$ is the reflection about $H$. We have that $(\calD(C_{f_k}))^{\gamma,v}$ is causally convex. Let $\partial_+ (\calD(C_{f_k}))^{\gamma,v}$ be the strict future boundary of $(\calD(C_{f_k}))^{\gamma,v}$, thus $\partial_+ (\calD(C_{f_k}))^{\gamma,v}$ is achronal. Define the function $\nu_k$ by
\begin{align*}
(\calD(C_{f_k}))^{\gamma,v} = \{ (t,x) \in \mathbb{M}^{n,1}:  r\leq t< \nu_k(x) \},
\end{align*}
similarly as the definition of the function $u_k$. The sequence $\{\nu_k\}_{k\in \mathbb{N}}$ is uniformly Lipschitz and $\nu_k(x) = R$ for a sufficiently large radius $R$. Then any subsequence of $\{\nu_k\}_{k\in \mathbb{N}}$ has a further converging subsequence. By proposition \ref{prop 3.4}, $\{ (\calD(C_{f_k}))^{\gamma,v} \}_{k\in\mathbb{N}}$ converges to $(\calD(C_{\barf}))^{\gamma,v}$ with respect to the symmetric difference volume. Define the function $\nu$ by
\begin{align*}
(\calD(C_{\barf}))^{\gamma,v} = \{ (t,x) \in \mathbb{M}^{n,1}:  r\leq t< \nu(x) \}.
\end{align*}
Thus any converging subsequence of $\{\nu_k\}_{k\in \mathbb{N}}$ converges to the same limit $\nu$. Therefore $\{ \nu \}_{k\in \mathbb{N}}$ converges to $\nu$. 

Introduce the notations $\calE_k$ and $\calE$ by
\begin{align*}
\calE_k = (\calD(C_{f_k}))^{\gamma,v},
\quad
\calE = (\calD(C_{\barf}))^{\gamma,v}.
\end{align*}
The strict future boundary $\partial_+ \calE$ of $\calE$ is the part of the graph of $\nu$ which is not in the lightcone $C_0$ that
\begin{align*}
\partial_+ \calE = \{ (t,x) \in \mathbb{M}^{n,1}: t= \nu(x) >r \}
\end{align*}
We define the following set from the domain of dependence  of $\partial_+ \calE$
\begin{align*}
\calD = \calD(\calE) = \mathrm{Int}(\calD(\partial_+ \calE)) \cup \partial_- (\calD(\partial_+ \calE)).
\end{align*}
$\partial_- (\calD(\partial_+ \calE))$ is an open finite lightcone, denoted by $C_{\barf_{\nu}}$. Then we have that
\begin{align*}
\calD = \calD(C_{\barf_{\nu}}).
\end{align*}
The goal is to show that the above constructed $C_{\barf_{\nu}}$ belongs to $\overline{P^v_f}$ and is actually $(C_{\barf})^{\gamma, v}$.

\vspace{5pt}
\noindent{\underline{\textsc{Step 2.}} }
We show that $C_{\barf_{\nu}} \in \overline{P^v_f}$. 

Consider the sequence of open finite lightcones $\{ (C_{f_k})^{\gamma,v}\}$ and the corresponding domains of dependence $\{ \calD((C_{f_k})^{\gamma,v})\}$. Following the definition of $\calD$, we introduce $\partial_+ \calE_k$ and $\calD_k$ by
\begin{align*}
\partial_+ \calE_k = \{ (t,x) \in \mathbb{M}^{n,1}: t= \nu_k (x) >r \},
\quad
\calD_k = \calD(\calE_k) = \mathrm{Int}(\calD(\partial_+ \calE_k)) \cup \partial_- \calD(\partial_+ \calE_k).
\end{align*}
Since $\calD_k = \calD((C_{f_k})^{\gamma,v})$, it is sufficient to prove that $\lim_{k\rightarrow+\infty}|\calD_k \triangle \calD| =0$, then we have that $\lim_{k\rightarrow +\infty}d_{\calD} ((C_{f_k})^{\gamma,v}, C_{\barf_{\nu}} ) =0$ which implies that $C_{\barf_{\nu}} \in \overline{P^v_f}$. 

Introduce the functions $\baru_k$ and $\baru$ by
\begin{align*}
\calD_k = \{ (t,x) \in \mathbb{M}^{n,1}:  r\leq t< \baru_k(x) \},
\quad
\calD = \{ (t,x) \in \mathbb{M}^{n,1}:  r\leq t< \baru (x) \}.
\end{align*}
We prove that  $\{ \baru_k \}_{k\in\mathbb{N}}$ converges to $\baru$.
\begin{enumerate}[label=\textit{\roman*}.]
\item
Show 
$
\baru(x) \geq \lim_{k\rightarrow +\infty} \baru_k(x).
$

Let $q \in J^+(o)$ with the coordinate $(t,x) = (\baru(x), x) \in \partial_{+} \calD \cup C_0$. Define $q_{+\delta} \in I^+(o)$ with the coordinate $(t,x) = (\baru(x) + \delta, x), \delta>0$. Then $\uC_{q_{+\delta}} \cap C_0 \not\subset C_{\barf} \cup S_{\barf}$. 

Let $p \in \uC_{q_{+\delta}} \cap C_0$ and $p\notin C_{\barf} \cup S_{\barf}$. Let $P_{opq_{+\delta}}$ be the $2$-dimensional plane determined by three points $o$, $p$, $q_{+\delta}$. Thus the null segment $\overline{pq_{+\delta}} \cap \overline{\calD \cap P_{opq_{+\delta}}} = \emptyset$, which implies that $\overline{pq_{+\delta}} \cap \overline{\calE \cap P_{opq_{+\delta}}} = \emptyset$. Therefore for $k$ sufficiently large, $\overline{pq_{+\delta}} \cap \overline{\calE_k \cap P_{opq_{+\delta}}} = \emptyset$, hence $q_{+\delta} \notin \overline{\calD_k}$ and $\baru(x) + \delta > \baru_k(x)$. By taking the limit as $k\rightarrow +\infty$, we obtain that
\begin{align*}
\baru(x) + \delta \geq \lim_{k\rightarrow +\infty} \baru_k(x).
\end{align*}
Since $\delta$ is arbitrary, we have that
\begin{align*}
\baru(x) \geq \lim_{k\rightarrow +\infty} \baru_k(x).
\end{align*}

\item
Show 
$
\baru(x) \leq \lim_{k\rightarrow +\infty} \baru_k(x).
$

If $\baru(x) = r$, then $\baru(x) = r \leq \baru_k(x) \Rightarrow \baru(x) \leq \lim_{k\rightarrow +\infty} \baru_k(x)$.

Otherwise if $\baru(x) > r$, let $q \in I^+(o)$ with the coordinate $(t,x) = (\baru(x), x) \in \partial_{+} \calD$. Define  $q_{-\delta} \in I^+(o)$ with the coordinate $(t,x) = (\baru(x) - \delta, x), \delta>0$. Then $\uC_{q_{-\delta}} \cap C_0 \subset C_{\barf}$. Then for sufficiently large $k$, $\uC_{q_{-\delta}} \cap C_0 \subset \calE_k$ which implies that $\uC_{q_{-\delta}} \cap C_0 \subset \calD_k$. Thus for sufficiently large $k$,
\begin{align*}
\baru(x) - \delta < \baru_k(x).
\end{align*}
Taking the limit $k\rightarrow +\infty$, we obtain that
\begin{align*}
\baru(x) - \delta \leq \lim_{k\rightarrow +\infty} \baru_k(x).
\end{align*}
Since $\delta$ is arbitrary, we have that
\begin{align*}
\baru(x) \leq \lim_{k\rightarrow +\infty} \baru_k(x).
\end{align*}
\end{enumerate}
Therefore
\begin{align*}
\baru(x) = \lim_{k\rightarrow +\infty} \baru_k(x)
\quad
\Rightarrow
\quad
\lim_{k\rightarrow+\infty}|\calD_k \triangle \calD| =0
\quad
\Rightarrow
\quad
\lim_{k\rightarrow +\infty}d_{\calD} ((C_{f_k})^{\gamma,v}, C_{\barf_{\nu}} ) =0.
\end{align*}
Hence $C_{\barf_{\nu}} \in \overline{P^v_f}$.

\vspace{5pt}
\noindent{\underline{\textsc{Step 3.}} }
We show that $C_{\barf_{\nu}}= (C_{\barf})^{\gamma,v}$. It is sufficient to show that $\calD(C_{\barf_{\nu}})= \calD((C_{\barf})^{\gamma,v})$.
\begin{enumerate}[label=\textit{\roman*}.]
\item
$\calD(C_{\barf_{\nu}}) \subset \calD((C_{\barf})^{\gamma,v})$. Note $\calD(C_{\barf_{\nu}}) = \calD(\calE) = \calD((\calD(C_{\barf}))^{\gamma,v})$. Since $(\calD(C_{\barf}))^{\gamma,v} \subset \calD((C_{\barf})^{\gamma,v})$, we have $\calD((\calD(C_{\barf}))^{\gamma,v}) \subset \calD(\calD((C_{\barf})^{\gamma,v})) = \calD((C_{\barf})^{\gamma,v})$, hence
\begin{align*}
	\calD(C_{\barf_{\nu}}) \subset \calD((C_{\barf})^{\gamma,v}).
\end{align*}

\item
$\calD(C_{\barf_{\nu}}) \supset \calD((C_{\barf})^{\gamma,v})$. Since $(C_{\barf})^{\gamma,v} \subset (\calD(C_{\barf}))^{\gamma,v} = \calE$, then $\calD((C_{\barf})^{\gamma,v}) \subset \calD(\calE) = \calD(C_{\barf_{\nu}})$.
\end{enumerate}
Combining the above three steps, we prove that $(C_{\barf})^{\gamma,v} \in \overline{P^v_f}$. 
\end{proof}

We consider the variational problem of finding the supremum of the functional $\calV$ in the metric space $\overline{P_f^v}$. Formally $\calV(C_{\barf}) = d_{\calD}(C_{\barf}, \emptyset)$. The existence of the extremal point of $\calV$ in $\overline{P_f^v}$ follows easily from the compactness of $(\overline{P_f^v}, d_{\calD})$.
\begin{proposition}\label{prop 4.9}
Any maximising sequence of the functional $\calV$ in $\overline{P_f^v}$ has a converging subsequence in $(\overline{P_f^v}, d_{\calD})$ whose limit achieves the maximum of $\calV$ in $\overline{P_f^v}$.
\end{proposition}

\subsection{Extremal set of volume functional}

We introduce the subset of $\overline{P^{v}_f}$ which maximises the volume functional $\calV$ and study the property of this subset.
\begin{definition}\label{def 4.16}
Define the set $M\overline{P^{v}_f}$ as the set of the open finite lightcone $C_{\barf}$ which maximises $\calV$ in $\overline{P^{v}_f}$, i.e.
\begin{align*}
M\overline{P^{v}_f} = \{C_{\barf} \in \overline{P^{v}_f} : \calV(C_{\barf}) = \max_{\overline{P^{v}_f}} \calV \}.
\end{align*}
\end{definition}
\begin{lemma}\label{lem 4.11}
The set $M\overline{P^{v}_f}$ has the following properties.
\begin{enumerate}[label=\alph*.]
\item
$M\overline{P^{v}_f}$ is a closed set under the distance $d_{\calD}$.

\item
$M\overline{P^{v}_f}$ is invariant under the Lorentz polarisation about any pair $(H,v)$ where $H$ is a timelike hyperplane through the origin $o$ not containing $v$.

\item
For any $C_{\barf} \in M\overline{P^{v}_f}$, we have that
\begin{align*}
\calD( (C_{\barf})^{\gamma,v}) = (\calD(C_{\barf}))^{\gamma,v}.
\end{align*}
\end{enumerate}
\end{lemma}
\begin{proof}
\begin{enumerate}[label=\textit{\alph*}.]
\item
This simply follows from that
\begin{align*}
|\calD(C_{f_1})| \geq |\calD(C_{f_2})| - |\calD(C_{f_1}) \triangle \calD(C_{f_2})|
\quad
\Leftrightarrow
\quad
\calV(C_{f_1}) \geq \calV(C_{f_2}) - d_{\calD}(C_{f_1}, C_{f_2}).
\end{align*}

\item This follows from corollary \ref{coro 3.11} the monotonicity of the volume of the domain of dependence of the open finite lightcone under the Lorentz polarisation, and proposition \ref{prop 4.8} the invariance of $\overline{P^{v}_f}$ under the Lorentz polarisation.

\item
This follows from corollary \ref{coro 3.11} and the maximising property of the volume functional $\calV$ in definition \ref{def 4.16} of $M\overline{P^{v}_f}$. 
\end{enumerate}
\end{proof}

\subsection{Lorentz polarisation in extremal set and polarisation in hyperbolic space}\label{subsec 4.4}
We present the following connection between the Lorentz polarisation in $\mathbb{M}^{n,1}$ and the polarisation in the hyperbolic space.
\begin{lemma}\label{lem 4.12}
Let $v$ be a future-directed vector and $H$ be a timelike hyperplane through $o$ not containing $v$. Let $S_{-\delta}$ be the hyperboloid
\begin{align*}
-t^2 + r^2 = -\delta, 
\quad 
\delta>0.
\end{align*}
$(S_{-\delta}, \eta|_{S_{-\delta}})$ is a hyperbolic space of constant negative sectional curvature $-\delta^{-1}$. Introduce the intersection of $S_{-\delta}$ with $H$, denoted by $L^H_{-\delta}$. $L^H_{-\delta}$ is a totally geodesic submanifold of $(S_{-\delta}, \eta|_{S_{-\delta}})$. Also introduce the intersection of $S_{-\delta}$ with the line through $o$ in $v$ direction, denoted by $o^{v}_{-\delta}$.
\begin{enumerate}[label=\alph*.]
\item
Let $\gamma_H$ be the Lorentz reflection about $H$. Denote the restriction of $\gamma_H$ on $S_{-\delta}$ by $\gamma^{H}_{-\delta}$. Then $\gamma^{H}_{-\delta}$ is the reflection of the hyperbolic space $S_{-\delta}$ about $L^H_{-\delta}$, which is an isometry of $(S_{-\delta}, \eta|_{S_{-\delta}})$ with $L^H_{-\delta}$ as the fixed point set.

\item
The polarisation of sets in $S_{-\delta}$ about the pair $(L^H_{-\delta}, o^{v}_{-\delta})$ is the restriction of the polarisation of sets in $\mathbb{M}^{n,1}$ about the pair $(H,v)$. Denote the polarisation of a set $U\subset S_{-\delta}$ about $(L^H_{-\delta}, o^{v}_{-\delta})$ by $U^{\gamma^{H}_{-\delta}, o^v_{-\delta}}$.
\end{enumerate}
\end{lemma}
The proof is straightforward, so we omit here. We now establish the connection between the Lorentz polarisation in the extremal set $M\overline{P^{v}_f}$ and the polarisation in the hyperbolic space.
\begin{lemma}\label{lem 4.13}
Let $C_{\barf} \in M\overline{P^{v}_f}$. Let $H$ be a timelike hyperplane through $o$ not containing $v$. Then we have that
\begin{align*}
(\calD(C_{\barf}) \cap S_{-\delta})^{\gamma^{H}_{-\delta},o^v_{-\delta}}
=
\calD(C_{\barf}^{\gamma,v}) \cap S_{-\delta}.
\end{align*}
\end{lemma}
\begin{proof}
By lemma \ref{lem 4.12},
\begin{align*}
(\calD(C_{\barf}) \cap S_{-\delta})^{\gamma^{H}_{-\delta},o^v_{-\delta}} 
= 
(\calD(C_{\barf}))^{\gamma,v} \cap S_{-\delta}.
\end{align*}
By lemma \ref{lem 4.11}, we have that $\calD((C_{\barf})^{\gamma,v}) =\calD(C_{\barf}))^{\gamma,v} $, thus
\begin{align*}
(\calD(C_{\barf}) \cap S_{-\delta})^{\gamma^{H}_{-\delta},o^v_{-\delta}} 
= 
\calD((C_{\barf})^{\gamma,v}) \cap S_{-\delta}.
\end{align*}
\end{proof}

Consider the function $\phi$ from $I^+(o)$ to $\mathbb{R}_-$ that
\begin{align*}
\phi:
\quad
I^+(o) \rightarrow \mathbb{R}_-
\quad
(t,x) \mapsto -t^2 + r^2.
\end{align*}
$\phi$ is Lipschitz on any bounded set of $I^+(o)$. Let $C_{\barf}$ be an open finite lightcone. Then we have that for almost all $y \in \mathbb{R}_-$, $\phi^{-1}(y) \cap \partial_+ \calD(C_{\barf})$ is of finite $(n-1)$-dimensional Hausdorff measure and $\phi^{-1}(y) \cap \calD(C_{\barf})$ is a set of finite perimeter or empty, where the Hausdorff measure and the perimeter are taken with respect to the Euclidean metric $g= \ed t^2 + \ed x_1^2 + \cdots + \ed x_n^2$. Then we have the following lemma for $\calD(C_{\barf}) \cap S_{-\delta}$.
\begin{lemma}\label{lem 4.14}
Let $C_{\barf}$ be an open finite lightcone. For almost all $\delta>0$, $\calD(C_{\barf}) \cap S_{-\delta}$ is either an open set of finite perimeter in $(S_{-\delta}, \eta|_{S_{-\delta}})$ or empty. 
\end{lemma}
\begin{proof}
We assume that $C_{\barf} \subset \calC_R$. Then consider the Euclidean metric $g= \ed t^2 + \ed x_1^2 + \cdots + \ed x_n^2$ in the compact set $\calD(\calC_R)$. Recall Eilenberg's inequality (see \cite{BZ88} theorem 13.3.1) for the Lipschitz map $\phi(t,x) = -t^2 + r^2$ in $\calD(\calC_R)$
\begin{align*}
	\overline{\int}_{\mathbb{R}_-} H_{n-1}((\partial_+ \calD(C_{\barf})) \cap S_{-\delta}) \ed \delta
	\leq
	\frac{v_1 v_{n-1}}{v_n} [ \mathrm{Lip}(\phi|_{\calD(\calC_R)}) ] H_n (\partial_+ \calD(C_{\barf})),
\end{align*}
where $\overline{\int}_{\mathbb{R}_-}$ denotes the upper Lebesgue integral. Therefore we have that
\begin{align*}
	H_{n-1}((\partial_+ \calD(C_{\barf})) \cap S_{-\delta})
	<
	+\infty,
\end{align*}
for almost all $\delta >0$. Since $\partial(\calD(C_{\barf}) \cap S_{-\delta}) \subset (\partial_+ \calD(C_{\barf})) \cap S_{-\delta}$ in $S_{-\delta}$, by a result of Federer (see \cite{F69} theorem 4.5.11, \cite{AFP00} proposition 3.62), $(\partial_+ \calD(C_{\barf})) \cap S_{-\delta}$ is either a set of finite perimeter in $(S_{-\delta}, \eta|_{S_{-\delta}})$ or empty for almost all $\delta>0$.
\end{proof}

\begin{lemma}\label{lem 4.15}
Let $C_{\barf} \in M\overline{P^v_f}$. Consider the set $\overline{P^v_{\barf}}$, which is a closed subset of $M\overline{P^v_f}$. We have that for almost all $\delta \in \mathbb{R}_+$, there exists an open finite lightcone in $\overline{P^v_{\barf}}$, denoted by $C_{\barf_\delta}$, that $\calD(C_{\barf_{\delta}}) \cap S_{-\delta}$ is either an open geodesic ball at $o^v_{-\delta}$ of the hyperbolic space $(S_{-\delta}, \eta|_{S_{-\delta}})$ or empty.
\end{lemma}
\begin{proof}
By lemma \ref{lem 4.14}, for almost all $\delta>0$, $\calD(C_{\barf}) \cap S_{-\delta}$ is either a set of finite perimeter in $(S_{-\delta}, \eta|_{S_{-\delta}})$ or empty. 

If $\calD(C_{\barf}) \cap S_{-\delta} = \emptyset$, then by lemma \ref{lem 4.13}, for any $C_{\barf'} \in \overline{P^v_{\barf}}$, $\calD(C_{\barf'}) \cap S_{-\delta}$ is empty.

If $\calD(C_{\barf}) \cap S_{-\delta}$ is an open set of finite perimeter, then one can approximate the geodesic ball of the same volume with respect to the volume of symmetric difference by applying polarisations in the hyperbolic space $(S_{-\delta}, \eta|_{S_{-\delta}})$ to the set $\calD(C_{\barf}) \cap S_{-\delta}$. More precisely, let $B^v_{-\delta}(l)$ be the geodesic ball of radius $l$ at $o^v_{-\delta}$ in $S_{-\delta}$ with the same volume as $C_{\barf} \cap S_{-\delta}$. Let $C_{\barf'_i} \in \overline{P^v_{\barf}}$ converges to $C_{\barf'}$ such that
\begin{align*}
	\lim_{i\rightarrow +\infty} 
	( \vert (\calD(C_{\barf'_i})\cap S_{-\delta}) \setminus B^v_{-\delta}(l) \vert 
	+ \vert B^v_{-\delta}(l) \setminus (\calD(C_{\barf'_i})\cap S_{-\delta}) \vert )
	=0.
\end{align*}
Let $\baru'_i$ and $\baru'$ be the parameterisation functions of $\partial_+ \calD(C_{\barf'_i})$ and $\partial_+ \calD(C_{\barf'})$, i.e.
\begin{align*}
	\calD(C_{\barf'_i})
	=
	\{ (t,x) \in \mathbb{M}^{n,1}:  r\leq t< \baru'_i(x) \},
	\quad
	\calD(C_{\barf'})
	=
	\{ (t,x) \in \mathbb{M}^{n,1}:  r\leq t< \baru'(x) \},
\end{align*}
therefore
\begin{align*}
	\calD(C_{\barf'_i}) \cap S_{-\delta}
	=
	\{ (t,x)\in S_{-\delta}: \baru'_i(x) > \sqrt{\delta + |x|^2} \}.
\end{align*}
Define $X^v_{-\delta}(l)$ be the set of $x$ that
\begin{align*}
	X^v_{-\delta}(l)
	=
	\{ x \in \mathbb{R}^n: (t,x) = (\sqrt{\delta + |x|^2}, x) \in B^v_{-\delta}(l)\},
\end{align*}
Since $\baru'(x) = \lim_{i\rightarrow +\infty} \baru'_i$, then for almost all $x \in X^v_{-\delta}(l)$, we have
\begin{align*}
	\baru'(x) = \lim_{i \rightarrow +\infty} \baru'_i(x) \geq \sqrt{\delta + |x|^2}.
\end{align*}
Therefore
\begin{align*}
	\vert B^v_{-\delta}(l) \setminus (\overline{\calD(C_{\barf'})}\cap S_{-\delta}) \vert 
	=0
	\quad
	\Rightarrow
	\quad
	B^v_{-\delta}(l) \subset (\overline{\calD(C_{\barf'})}\cap S_{-\delta}),
\end{align*}
then
\begin{align*}
	B^v_{-\delta}(l)
	\subset
	\calD(C_{\barf'})\cap S_{-\delta}.
\end{align*}
Since $B^v_{-\delta}(l)$ and $\calD(C_{\barf'})\cap S_{-\delta}$ have the same volume, we obtain that 
\begin{align*}
	B^v_{-\delta}(l) = \calD(C_{\barf'})\cap S_{-\delta}.
\end{align*}
The lemma is proved.
\end{proof}

Applying lemma \ref{lem 4.15} to a sequence $\{ \delta_k \}_{k \rightarrow +\infty}$ converging to $0$, we can show that there exists an open finite lightcone truncated by a spacelike hyperplane in $M\overline{P^v_f}$.
\begin{proposition}\label{prop 4.16}
Let $l$ be a positive number. Introduce the spacelike hyperplane through the point $l\frac{v}{|v|}$ and orthogonal to $v$, denoted by $H_l^v$. Also introduce the open finite lightcone in $C_0$ truncated by the hyperplane $H_l^v$, denoted by $C_l^v$, i.e.
\begin{align*}
C_l^v = \{ x \in C_0:  \langle x, v \rangle > \langle l\frac{v}{|v|}, v \rangle \}.
\end{align*}
There exists a positive number $l>0$ such that $C_l^v \in M \overline{P^v_f}$.
\end{proposition}
\begin{proof}
Choose a decreasing sequence $\{ \delta_k \}_{k \rightarrow +\infty}$ converging to $0$, such that there exists an open finite lightcone $C_{\barf_k} \in M\overline{P^v_f}$ that $\calD(C_{\barf_k}) \cap S_{-\delta_k}$ is a geodesic ball of $(S_{-\delta_k},\eta|_{S_{-\delta_k}})$ centering at $o^v_{-\delta_k}$. The future domain of dependence of such a geodesic ball is simple, which is equal to the part of $\calD(C_{l_k}^v)$ in the future of $S_{-\delta_k}$ for some $l_k>0$. Thus for any $K\in \mathbb{N}$, $\calD(C_{\barf_k}) \cap I^+(S_{-\delta_K}) = \calD(C_{l_k}^v) \cap I^+(S_{-\delta_K}) $ for all $k\geq K$.

By the completeness and compactness of $M\overline{P^v_f}$, the sequence of open finite lightcones $\{ C_{\barf_k} \}_{k\in \mathbb{N}}$ converges to an open finite lightcone $C_{\barf}$ in $M\overline{P^v_f}$. By the previous description of $\calD(C_{\barf_k}) \cap I^+(S_{-\delta_k})$, we see that the limit open finite lightcone $C_{\barf}$ must be $C_l^v$ for some $l>0$.
\end{proof}

\subsection{Semicontinuity of $L^{n-1}$ norm for parameterisation function}

In order to apply the previously obtained results on the Lorentz polarisation to the variational problem of the isoperimetric inequality, we need to study the limit of the $L^{n-1}$ norm of the parameterisation function under the convergence of open finite lightcones with respect to the distance $d_{\calD}$.

Concerning the regularity of the parameterisation function of an open finite lightcone, examples \ref{ex 2.5}, \ref{ex 2.9.a} and \ref{ex 2.9.b} show that the parameterisation function could be rather wild in the lightcone. Thus it is not easy to study the convergence of the parameterisation functions of a convergent sequence of open finite lightcones. In order to overcome this difficulty, we make use of the structure of the domain of dependence of an open finite lightcone to prove the following helpful lemmas.

\begin{lemma}\label{lem 4.17}
Let $C_f$ be an open finite lightcone and $\calD(C_f)$ be the domain of dependence of $C_f$. Let $e_0=\partial_t$ and $C_{\delta e_0}$ be the future lightcone emanating from the vertex $\delta e_0$. Introduce the open finite lightcone
\begin{align*}
C_{\delta e_0} \cap \calD(C_f).
\end{align*}
Define the parameterisation function $f_{\delta}$ of $C_{\delta e_0} \cap \calD(C_f)$ similar as in proposition \ref{def 2.1.a}: let $\{ r, \vartheta \}$ be the coordinate system on $C_{\delta e_0}$, then the parameterisation function $f_{\delta}$ is defined as
\begin{align*}
C_{\delta e_0} \cap \calD(C_f) = \{(r, \vartheta): r < f_{\delta}(\vartheta)  \} \cup \{ \delta e_0 \}.
\end{align*}
We have that
\begin{align*}
f_{\delta} \leq f - \frac{\delta}{2},
\quad
\lim_{\delta \rightarrow 0^+} f_{\delta} = f.
\end{align*}
Moreover for any $0\leq \eta \leq \delta$, we have that
\begin{align*}
f_{\delta} \leq f_{\eta}  - \frac{\delta - \eta}{2},
\quad
\lim_{\delta \rightarrow \eta^+} f_{\delta} = f_{\eta}.
\end{align*}
\end{lemma}
\begin{proof}
Let $q$ be the point of the coordinate 
\begin{align*}
(t,r,\vartheta) = (\delta + f(\vartheta) - \frac{\delta}{2}+\epsilon, f(\vartheta) - \frac{\delta}{2} +\epsilon, \vartheta)
\end{align*} 
where $\epsilon >0$. 
Let $\vec{l}$ be the past-directed null vector from $o$ to the point $(t,r,\vartheta) = (-\frac{\delta}{2}, \frac{\delta}{2}, \vartheta)$. See figure \ref{fig 15}. Then the point $q+\vec{l} \in J^-(q)$. Note the coordinate of $q+\vec{l}$ is
\begin{align*}
(t,r,\vartheta) = (f(\vartheta) +\epsilon, f(\vartheta)+\epsilon, \vartheta) \notin C_f.
\end{align*}
Thus $q \notin \calD(C_f)$. Hence $f_{\delta} < f - \frac{\delta}{2} +\epsilon$ for all $\epsilon >0$, which implies that 
\begin{align*}
f_{\delta} \leq f - \frac{\delta}{2}.
\end{align*}
\begin{figure}[h]
\centering
\begin{tikzpicture}[scale=0.9]
\draw (-4.5, 4.5) 
	-- (0,0) 
	node[below] {$o$} 
	-- (4.5, 4.5);
\draw (-4.5, 4.5+1.6) 
	-- (0,1.6) 
	node[right] {$\delta e_0$} 
	-- (4.5, 4.5+1.6); 
\draw[->] (0,0) 
	-- (0,3);
\draw[fill=black] (4,4) 
	circle [radius=0.04];
\node[below right] at (4,4) 
	{\footnotesize $(f(\vartheta), f(\vartheta), \vartheta)$};
\draw[->,thick] (4,4) -- (4-0.7,4+0.7) -- (4-0.8,4+0.8);
\draw[fill=black] (4-0.8,4+0.8) 
	circle [radius=0.04];
\node[above left] at (4-0.8,4+0.8) 
	{\footnotesize $(f(\vartheta)+\frac{\delta}{2}, f(\vartheta)-\frac{\delta}{2}, \vartheta)$};
\draw[->,thick] (4-0.8,4+0.8) -- (4-0.8+0.6,4+0.8+0.6);
\draw[fill=black] (4-0.8+0.6,4+0.8+0.6) 
	circle [radius=0.04];
\node[above left] at (4-0.8+0.6,4+0.8+0.6) 
	{\footnotesize $(f(\vartheta)+\frac{\delta}{2}+\epsilon, f(\vartheta)-\frac{\delta}{2}+\epsilon, \vartheta)=q$};
\draw[->,thick] (4-0.8+0.6,4+0.8+0.6) -- (4+0.6,4+0.6);
\draw[fill=black] (4+0.6,4+0.6) 
	circle [radius=0.04];
\node[below right] at (4+0.6,4+0.6) 
	{\footnotesize $q+\vec{l}=(f(\vartheta)+\epsilon, f(\vartheta)+\epsilon, \vartheta)$};
\end{tikzpicture}
\caption{$f_{\delta} \leq f-\frac{\delta}{2}$.}
\label{fig 15}
\end{figure}

Let $p$ be the point of the coordinate 
\begin{align*}
(t,r,\vartheta) = (f(\vartheta_p) - \epsilon, f(\vartheta_p) -\epsilon, \vartheta_p),
\end{align*}
then $p \in C_f$. There exist a neighbourhood $U \subset \mathbb{S}^{n-1}$ of $\vartheta_p$ and a positive number $r_0$, then $f(\vartheta) > r_0$, for all $\vartheta \in \mathbb{S}^{n-1}$ and
\begin{align*}
f(\vartheta) > f(\vartheta_p) - \epsilon,
\quad
\vartheta \in U.
\end{align*}
Define $f_{\epsilon}$ by
\begin{align*}
f_{\epsilon}(\vartheta)
=
\left\{
\begin{aligned}
&
f(\vartheta_q) - \epsilon,
\quad
\vartheta \in U,
\quad
\\
&
r_0,
\quad
\vartheta \in \mathbb{S}^{n-1} \setminus U.
\end{aligned}
\right.
\end{align*}
Then $f_{\epsilon} < f$ and $C_{f_{\epsilon}} \subset C_f$.  Let $q_{\delta}$ be the point of the coordinate 
\begin{align*}
q_{\delta}: \quad (t,r,\vartheta) = (\delta + f(\vartheta_q) - \frac{\delta}{2}-\epsilon, f(\vartheta_q) - \frac{\delta}{2} -\epsilon, \vartheta_q) \in C_{\delta e_0}. 
\end{align*}
See figure \ref{fig 16}.
\begin{figure}[h]
\centering
\begin{tikzpicture}
\draw (-5,5) 
	-- (0,0) 
	node[below] {$o$} 
	-- (5,5);
\draw (-5,5+1.6) 
	-- (0,1.6) 
	node[right] {$\delta e_0$} 
	-- (5,5+1.6); 
\draw[->] (0,0) 
	-- (0,3);
\draw[fill=black] (4,4) 
	circle [radius=0.04];
\node[below right] at (4,4) 
	{\footnotesize $(f(\vartheta_p), f(\vartheta_p), \vartheta_p)$};
\draw[->,thick] (4,4) -- (4-0.7,4+0.7) -- (4-0.8,4+0.8);
\draw[fill=black] (4-0.8,4+0.8) 
	circle [radius=0.04];
\node[above left] at (4-0.8,4+0.8) 
	{\footnotesize $(f(\vartheta_p)+\frac{\delta}{2}, f(\vartheta_p)-\frac{\delta}{2}, \vartheta_p)$};
\draw[fill=black] (4-0.8-0.6,4+0.8-0.6) 
	circle [radius=0.04];
\node[above left] at (4-0.8-0.6,4+0.8-0.6) 
	{\footnotesize $(f(\vartheta_p)+\frac{\delta}{2}-\epsilon, f(\vartheta_p)-\frac{\delta}{2}-\epsilon, \vartheta_p)=q_{\delta}$};
\draw[<-,thick] (4-0.8-0.6,4+0.8-0.6) -- (4-0.6,4-0.6);
\draw[->,thick] (4,4) -- (4-0.6,4-0.6);
\draw[fill=black] (4-0.6,4-0.6) 
	circle [radius=0.04];
\node[below right] at (4-0.6,4-0.6) 
	{\footnotesize $p=(f(\vartheta_p)-\epsilon, f(\vartheta_p)-\epsilon, \vartheta_p)$};
\end{tikzpicture}
\caption{$p$ and $q_{\delta}$.}
\label{fig 16}
\end{figure}
For $\delta$ sufficiently small, we have that
\begin{align*}
J^-(q_{\delta}) \cap C_0  \subset C_{f_{\epsilon}} \subset C_f.
\end{align*}
See figure \ref{fig 17}.
\begin{figure}[h]
\centering
\begin{tikzpicture}
\draw[fill] (0,0) node[right] {$o$} circle[radius=0.01];
\draw plot [smooth cycle] coordinates {(4.3,0) (3,2) (1,2) (-1,2) (-2,-1) (-2,-3) (1,-1.5) (3,-1.8)};
\draw[fill] (4.3,0) node[right] {\footnotesize $(f(\vartheta_p), \vartheta_p)$} circle[radius=0.04];
\draw[->] (5,1) node[right] {\footnotesize $p:(f(\vartheta_p)-\epsilon, \vartheta_p)$} to [out=180,in=60] (4.05,0.05);
\draw[fill] (4,0) circle[radius=0.04];
\draw (0.984*4,-0.174*4) arc(-10:10:4);
\draw (0.984*1.6,-0.174*1.6) arc(-10:-350:1.6);
\draw[dashed] (0.984*1.6,-0.174*1.6) -- (0.984*4,-0.174*4);
\draw[dashed] (0.984*1.6,0.174*1.6) -- (0.984*4,0.174*4);
\draw[<->] (0,0.1) -- (0,0.75) node[right] {$r_0$} -- (0,1.5);
\draw[domain=0:360,smooth,variable=\t]
  plot ({4*((1-0.95)/(1-0.95*cos(\t)))*cos(\t)},{1.4*((1-0.95)/(1-0.95*cos(\t)))*sin(\t)});
\node[below] at ({4*((1-0.95)/(1-0.95*cos(270)))*cos(270)},{1.4*((1-0.95)/(1-0.95*cos(270)))*sin(270)}) {\footnotesize $J^{-}(q_{\delta}) \cap C_0$};
\node[above right] at ({1.6*cos(45)},{1.6*sin(45)}) {\footnotesize $f_{\epsilon}$};
\end{tikzpicture}
\caption{$J^{-}(q_{\delta}) \cap C_0 \subset C_{f_{\epsilon}}$}
\label{fig 17}
\end{figure}
Hence for $\delta$ sufficiently small, $q_{\delta} \in \calD(C_f)$, thus
\begin{align*}
f_{\delta}(\vartheta_q) > f(\vartheta_q) - \frac{\delta}{2}-\epsilon
\quad
\Rightarrow
\lim_{\delta \rightarrow 0^+} f_{\delta}(\vartheta_q)  \geq  f(\vartheta_q) -\epsilon.
\end{align*}
Since $\epsilon$ is an arbitrary small positive number and $\vartheta_q$ is arbitrary in $\mathbb{S}^{n-1}$, we obtain that
\begin{align*}
\lim_{\delta \rightarrow 0^+} f_{\delta}(\vartheta)  \geq  f(\vartheta),
\quad
\forall \vartheta \in \mathbb{S}^{n-1}.
\end{align*}
Combining with $f_{\delta} \leq f- \frac{\delta}{2}$, we obtain that 
\begin{align*}
\lim_{\delta \rightarrow 0^+} f_{\delta}= f.
\end{align*}

The case of $0\leq \eta \leq \delta$ is proved in the same way by replacing $C_f$ with $C_{\eta e_0} \cap \calD(C_f)$.
\end{proof}

Lemma \ref{lem 4.17} implies that the parameterisation function $f$ of any open finite lightcone $C_f$ can be approximated by the parameterisation functions $f_{\delta}$ of other open finite lightcones $C_{\delta e_0} \cap \calD(C_f)$ in the future of $C_f$. We shall use this approximation to obtain the semicontinuity of the $L^{n-1}$ norm of the parametrisation function.
\begin{lemma}\label{lem 4.18}
Let $\{ C_{f_k} \}$ be a convergent sequence of open finite lightcones under the distance $d_{\calD}$, which converges to $C_f$. Suppose that $\{ C_{f_k} \}$ is uniformly bounded by $C_R$, and $C^{e_0}_l \subset C_f$. For any $0\leq \delta < l$, consider the parameterisation function $f_{\delta}$ of the open finite lightcone $C_{\delta e_0} \cap \calD(C_f)$, and similarly the parameterisation function $f_{k,\delta}$ of $C_{\delta e_0} \cap \calD(C_{f_k})$ (which always exists for $k$ sufficiently large). We have that
\begin{align*}
f_{\delta} \leq \liminf_{k\rightarrow +\infty} f_{k,\delta-\eta} - \frac{\eta}{2},
\end{align*}
for any $0< \eta \leq \delta < l$. In particular in the case $\eta = \delta$, we have that
\begin{align*}
f_{\delta} \leq \liminf_{k\rightarrow +\infty} f_k - \frac{\delta}{2}.
\end{align*}
\end{lemma}
\begin{proof}
Let $q$ be the point in the boundary $C_{\delta e_0} \cap \partial_+ \calD(C_f)$ of the open finite lightcone $C_{\delta e_0} \cap \calD(C_f)$ with the coordinate 
\begin{align*}
q: 
\quad 
(t, r, \vartheta) = (\delta + f_{\delta}(\vartheta), f_{\delta}(\vartheta), \vartheta).
\end{align*} 
Then the point $q-\epsilon e_0$ with the coordinate 
\begin{align*}
q-\epsilon e_0:
\quad
(t, r, \vartheta) = (\delta - \epsilon + f_{\delta}(\vartheta), f_{\delta}(\vartheta), \vartheta)
\end{align*} 
lies in the domain of dependence $\calD(C_f)$ for any sufficiently small positive number $\epsilon$. Thus for any sufficiently large $k$, $q-\epsilon e_0$ also lies in the domain of dependence $\calD(C_{f_k})$. Hence for any sufficiently large $k$, we have that
\begin{align*}
f_{\delta}(\vartheta) 
< 
f_{k,\delta-\epsilon} (\vartheta).
\end{align*}
By lemma \ref{lem 4.17}, we have that for any $k$ sufficient large and $\epsilon \leq \eta \leq \delta$,
\begin{align*}
f_{\delta}(\vartheta) 
< 
f_{k,\delta-\epsilon} (\vartheta) 
\leq 
f_{k,\delta-\eta} (\vartheta) - \frac{\eta - \epsilon}{2}.
\end{align*}
Thus taking the limit that $k\rightarrow +\infty$,
\begin{align*}
f_{\delta}(\vartheta) 
\leq 
\liminf_{k\rightarrow +\infty} f_{k,\delta-\eta} (\vartheta) - \frac{\eta - \epsilon}{2}.
\end{align*}
Since $\epsilon$ can be arbitrarily small, then we have that for $0 < \eta \leq \delta <l$,
\begin{align*}
f_{\delta}(\vartheta) 
\leq 
\liminf_{k\rightarrow +\infty} f_{k,\delta-\eta} (\vartheta) - \frac{\eta}{2}.
\end{align*}
The lemma is proved.
\end{proof}

Now we are ready to prove the semicontinuity of the $L^{n-1}$ norm for the parameterisation function of the open finite ligthcone under the distance $d_{\calD}$.
\begin{proposition}\label{prop 4.19}
Let $\{ C_{f_k} \}$ be a convergent sequence of open finite lightcones under the distance $d_{\calD}$, which converges to $C_f$. Suppose that $\{ C_{f_k} \}$ is uniformly bounded by $C_R$. We have that
\begin{align*}
\int_{\mathbb{S}^{n-1}} f^{n-1} \dvol_{\circg} 
\leq 
\liminf_{k\rightarrow +\infty} \int_{\mathbb{S}^{n-1}} f_k^{n-1} \dvol_{\circg}.
\end{align*}
\end{proposition}
\begin{proof}
Suppose that $C_l \subset C_f$, then by lemma \ref{lem 4.18} and Fatou's lemma
\begin{align*}
\int_{\mathbb{S}^{n-1}} f_{\delta}^{n-1} \dvol_{\circg} 
\leq 
\liminf_{k\rightarrow +\infty} \int_{\mathbb{S}^{n-1}} (f_k-\frac{\delta}{2})^{n-1} \dvol_{\circg}.
\end{align*}
Taking the limit $\delta \rightarrow 0^+$, by lemma \ref{lem 4.17} and the bounded convergence theorem, we obtain that
\begin{align*}
\int_{\mathbb{S}^{n-1}} f^{n-1} \dvol_{\circg} 
\leq 
\liminf_{k\rightarrow +\infty} \int_{\mathbb{S}^{n-1}} f_k^{n-1} \dvol_{\circg}.
\end{align*}
The proposition is proved.
\end{proof}

By propositions \ref{prop 3.8}, \ref{prop 4.16}, \ref{prop 4.19}, we have the following result for the set $\overline{P^v_f}$.
\begin{proposition}\label{prop 4.20}
Let $C_{\barf}$ be an open finite lightcone in the set $\overline{P^v_f}$. We have that
\begin{align*}
\int_{\mathbb{S}^{n-1}} \barf^{n-1} \dvol_{\circg} 
\leq 
\int_{\mathbb{S}^{n-1}} f^{n-1} \dvol_{\circg}.
\end{align*}
Let $l_f$ be the positive number that $n \omega_n l_f^{n-1} = \int_{\mathbb{S}^{n-1}} f^{n-1} \dvol_{\circg}$. Then there exists a positive number $l\in (0, l_f]$ such that $C^v_l \in M\overline{P^v_f}$.
\end{proposition}

\subsection{Identify extremal volume of domain of dependence of finite lightcone}\label{sec 4.3}
Finally in this section we can solve part \textit{a.} of the variational problem \ref{prob 4.1}. We state the answer in the following proposition.
\begin{proposition}\label{prop 4.21}
Recall the set $F_1$
\begin{align*}
F_1 = \{C_f: \int_{\mathbb{S}^{n-1}} f^{n-1} \dvol_{\circg} = n\omega_n \}.
\end{align*}
and the volume functional $\calV$
\begin{align*}
\calV(C_f) = |\calD(C_f)|.
\end{align*}
We have that
\begin{align*}
\max_{F_1} \calV = \frac{2}{n+1} \omega_n,
\end{align*}
which is equivalent to that
\begin{align*}
\calV(C_f) \leq \calV(C_1),
\quad
\forall C_f \in F_1.
\end{align*}
\end{proposition}

By proposition \ref{prop 4.20}, it is sufficient to prove that for any $C_f \in F_1$, $\max_{\overline{P^{e_0}_f}} \calV \leq  \calV(C_1)$. We have the following lemma, which is an immediate consequence of proposition \ref{prop 4.20}.
\begin{lemma}\label{lem 4.22}
For any $C_f \in F_1$, there exists a positive number $0< l_f \leq 1$ such that $C_{l_f} \in \overline{P^{e_0}_f}$ and
\begin{align*}
\max_{\overline{P^{e_0}_f}} \calV \leq  \calV(C_{l_f}) \leq \calV(C_1) = \frac{2}{n+1} \omega_n.
\end{align*}
\end{lemma}
Therefore proposition \ref{prop 4.21} follows from lemma \ref{lem 4.22} and we answer the part \textit{a.} of the variational problem \ref{prob 4.1}.

\section{Case of equality}\label{sec 5}
In this section, we solve part \textit{b.} of the variational problem \ref{prob 4.1}. We show that if an open finite lightcone in $F_1$ achieves the maximal volume of the domain of dependence, then the open finite lightcone must be $C^v_1$ for some future-directed timelike vector $v$.

\subsection{Equal perimeter hyperplane separation}\label{sec 5.1}
We introduce the following construction of reflection symmetric open finite lightcones from an arbitrary open finite lightcone by Lorentz reflection.
\begin{definition}\label{def 5.1}
Let $w$ be a spacelike vector and $H_w$ be the timelike hyperplane through the origin $o$ orthogonal to $w$. $H_w$ separates the spacetime into two parts. Let $H_w^+$ be the open half spacetime containing $w$ and $H_w^-$ be the other open half part. Let $U$ be a set of $\mathbb{M}^{n,1}$. Define the following two sets invariant under the reflection $\gamma$ about $H_w$:
\begin{align*}
	&
	U_{\gamma}^{w,+} 
	= 
	(U\cap H_w^+)
	\ \cup\ 
	(U\cap H)
	\ \cup\ 
	(\gamma(U) \cap H_w^-),
\\
	&
	U_{\gamma}^{w,-} 
	= 
	(U\cap H_w^-)
	\ \cup\ 
	(U\cap H)
	\ \cup\ 
	(\gamma(U) \cap H_w^+).
\end{align*}
We call $U_{\gamma}^{w,\pm} $ the positive or negative reflection symmetrisation of $U$ about $(w,H_w)$. 

Let $C_f$ be an open finite lightcone. Then introduce the parameterisation functions of the open finite lightcones $(C_f)_{\gamma}^{w,\pm}$. Denote the parameterisation functions of $(C_f)_{\gamma}^{w,\pm}$ by $f_{\gamma}^{w,\pm}$ respectively, i.e.
\begin{align*}
(C_f)_{\gamma}^{w,\pm} = C_{f_{\gamma}^{w,\pm}}.
\end{align*}
We call $f_{\gamma}^{w,\pm}$ the positive or negative reflection symmetrisation of the parameterisation function $f$ about $(w,H_w)$ respectively.
\end{definition}
We have the following lemma on the above constructed reflection symmetric open finite lightcone.
\begin{lemma}\label{lem 5.2}
Let $C_f$ be an open finite lightcone. Let $w$ be a spacelike vector and $H_w$ be its timelike orthogonal hyperplane through the origin $o$. We have that the reflection symmetrisation of $\calD(C_f)$ is contained in the domain of dependence of the reflection symmetrisation of $C_f$, i.e.
\begin{align*}
(\calD(C_f))_{\gamma}^{w,\pm} \subset \calD((C_f)_{\gamma}^{w,\pm}) = \calD(C_{f_{\gamma}^{w,\pm}}).
\end{align*}
\end{lemma}
\begin{proof}
It is easy to check the lemma for any spacelike hyperplane truncated open finite lightcone $C_l^v$. The general case follows easily from this special case.
\end{proof}

For an open finite lightcone, we introduce a particular class of timelike hyperplanes which separate the perimeter of its boundary equally.
\begin{definition}\label{def 5.3}
Let $C_f$ be an open finite lightcone. Let $w$ be a spacelike vector and $H_w$ be its timelike orthogonal hyperplane through the origin $o$. $H_w$ separates the sphere into two spherical caps $(\mathbb{S}^{n-1})_w^{\pm}$ corresponding to the domain of variables of the parameterisation of $C_f \cap H_n^{\pm}$. We call that $H_w$ an equal perimeter separation hyperplane if
\begin{align*}
\int_{(\mathbb{S}^{n-1})_w^+} f^{n-1} \dvol_{\circg} 
= 
\int_{(\mathbb{S}^{n-1})_w^-} f^{n-1} \dvol_{\circg} 
=
\frac{1}{2} \int_{\mathbb{S}^{n-1}} f^{n-1} \dvol_{\circg}.
\end{align*}
\end{definition}
We can easily construct equal perimeter separation hyperplanes by the following lemma.
\begin{lemma}\label{lem 5.4}
Let $C_f$ be an open finite lightcone. Choose an arbitrary $2$-dimensional timelike plane $P$. Consider the class of timelike hyperplanes orthogonal to $P$, denoted by $\mathcal{H}_P$. Equivalently 
\begin{align*}
\mathcal{H}_P = \{ H_w: w\in P \text{ and $w$ is spacelike} \}.
\end{align*}
There exists a unique equal perimeter separation hyperplane $H_w$ in $\mathcal{H}_P$.
\end{lemma}
\begin{proof}
Let $u_+(H_w) = \int_{(\mathbb{S}^{n-1})_w^+} f^{n-1} \dvol_{\circg}$ be an injective continuous function of $H_w$. The range of $u_+=$ is $(0, \int_{\mathbb{S}^{n-1}} f^{n-1} \dvol_{\circg})$. Then the lemma follows.
\end{proof}

We have a nice property for an extremal open finite lightcone of the volume functional $\calV$ in $F_1$ and its equal perimeter separation hyperplane.
\begin{lemma}\label{lem 5.5}
Let $C_f$ be an extremal open finite lightcone of the volume functional $\calV$ in $F_1$, i.e.
\begin{align*}
\calV(C_f) = |\calD(C_f)| = \max_{F_1} \calV = \calV(C_1) = \frac{2}{n+1} \omega_n.
\end{align*}
Introduce the notation $V_1 = \calV(C_1) = \frac{2}{n+1} \omega_n$. Suppose that $H_w$ is an equal perimeter separation hyperplane of $C_f$, then we have that
\begin{align*}
|\calD(C_f) \cap H_w^+ | = |\calD(C_f) \cap H_w^- | = \frac{1}{2} V_1.
\end{align*}
\end{lemma}
\begin{proof}
Assume the lemma is false. Without loss of generality, suppose that $|\calD(C_f) \cap H_w^+ | > \frac{1}{2} V_1$. Then consider the reflection symmetrisation $(C_f)_{\gamma}^{w,+}$ of $C_f$. We have that $(C_f)_{\gamma}^{w,+} \in F_1$ since $H_w$ is an equal perimeter separation hyperplane. However we arrive at a contradiction that
\begin{align*}
	&
	\phantom{\Rightarrow}
	(\calD(C_f))_{\gamma}^{w,+} \subset \calD((C_f)_{\gamma}^{w,+}) 
\\
	&
	\Rightarrow
	\calV((C_f)_{\gamma}^{w,+}) 
	= |\calD((C_f)_{\gamma}^{w,+})| 
	\geq |(\calD(C_f))_{\gamma}^{w,+}| 
	= 2 | \calD(C_f) \cap H_w^+ | 
	> V_1.
\end{align*}
Thus the lemma is true.
\end{proof}

By lemmas \ref{lem 5.2} and \ref{lem 5.5}, we have the following result for an extremal open finite lightcone and its equal perimeter separation hyperplane.
\begin{lemma}\label{lem 5.6}
Let $C_f$ be an extremal open finite lightcone of the volume functional $\calV$ in $F_1$. Suppose that $H_w$ is an equal perimeter separation hyperplane of $C_f$, then we have that
\begin{align*}
\calD(C_{f_{\gamma}^{w,\pm}})
=
\calD((C_f)_{\gamma}^{w,\pm})
=
(\calD(C_f))_{\gamma}^{w,\pm}.
\end{align*}
\end{lemma}

\subsection{Equal perimeter separation hyperplane noncrossing property}\label{sec 5.2}
We show that the null generator in the strict future boundary of an extremal open finite lightcone will not cross any equal perimeter separation hyperplane. The following lemma illustrates the key reason of this noncrossing property.
\begin{lemma}\label{lem 5.7}
Let $q_1, q_2$ be two points in $I^+$ and $|\overrightarrow{oq_1}| = |\overrightarrow{oq_2}|=2l$. Let $w$ be the vector $\overrightarrow{q_2 q_1}$. Suppose that $w$ is spacelike, which is equivalent to that $q_1, q_2$ are spacelike to each other. Let $H_w$ be the timelike hyperplane through the origin $o$ orthogonal to $w$. 
\begin{figure}[h]
\centering
\begin{tikzpicture}
\draw (-8 *.7,8 *.7) -- (0,0) -- (8 *.7,8 *.7);
\draw (-2 *.7,2 *.7) -- (4 *.7,8 *.7) node[above] {$q_1$} -- (6 *.7,6 *.7);
\draw (-4 *.7,4 *.7) -- (-1 *.7,7 *.7) node[above] {$q_2$} -- (3 *.7,3 *.7);
\node[above] at (4 *.7+1,8 *.7) {$H_{w}^{+}$};
\node[above] at (-1 *.7-2,7 *.7) {$H_{w}^{-}$};
\node[right] at (1 *.7,5 *.7) {$q$};
\draw (0,0) -- (1 *.7,5 *.7) -- (1.6 *.7,8 *.7) node[right] {$H_w$};
\draw[red,thick] (-2 *.7 +.05*1.5 ,2 *.7 +.05*1.5 ) -- (0,.1*1.5 ) -- (6 *.7 -.05*1.5 ,6 *.7 +.05*1.5 ) node[left] {$C_l^{v_1}$};
\draw[violet,thick] (-4 *.7 +.1*1.5 , 4 *.7 +.1*1.5 ) node[right] {$C_l^{v_2}$} -- (0,.2*1.5 ) -- (3 *.7 -.1*1.5 ,3 *.7 +.1*1.5 ) ;
\draw[olive,thick] (-2 *.7 -.05*1.5 , 2 *.7 -.05*1.5 ) -- (0,-.1*1.5 ) -- (3 *.7 +.05*1.5 , 3 *.7 -.05*1.5 ) node[below right] {$(C_l^{v_1})_{\gamma}^{w,-}$};
\draw[blue,thick] (-4 *.7 -.1*1.5 , 4 *.7 -.1*1.5 ) -- (0,-.2*1.5 ) -- (6 *.7 +.1*1.5 ,6 *.7 -.1*1.5 ) node[right] {$(C_l^{v_1})_{\gamma}^{w,+}$};
\end{tikzpicture}
\caption{Illustration of lemma \ref{lem 5.7}}
\label{fig 18}
\end{figure}

\begin{enumerate}[label=\alph*.]
\item 
$q_1$ and $q_2$ are symmetric about $H_w$, i.e. the reflection $\gamma$ about $H_w$ maps $q_1$, $q_2$ to each other. We have that $q_1 \in H_w^+$ and $q_2 \in H_w^-$. 

\item 
Introduce two vectors $v_1 = \frac{1}{2}\overrightarrow{oq_1}$ and $v_2 = \frac{1}{2} \overrightarrow{oq_2}$. We have $|v_1| = |v_2| = l$. Then the open finite lightcone $I^-(q_a) \cap C_0$ is the spacelike hyperplane truncated open finite lightcone $C_l^{v_a}$, $a=1,2$. The domain of dependence of $C_l^{v_a}$ and its strict future boundary are
\begin{align*}
\calD(C_l^{v_a}) = J^+(o) \cap I^-(q_a),
\quad
\partial_+ \calD(C_l^{v_a}) = \uC_{q_a} \cap I^+(o).
\end{align*}
The intersection of two strict future boundaries $\partial_+ \calD(C_l^{v_a})$ lies in $H_w$. We have that
\begin{align*}
\partial_+ \calD(C_l^{v_1})  \cap \partial_+ \calD(C_l^{v_2})
&=
\uC_{q_1} \cap \uC_{q_2} \cap I^+(o)
\\
&=
H_w \cap \uC_{q_1} \cap I^+(o)
= 
H_w \cap \uC_{q_2} \cap I^+(o).
\end{align*}

\item
Consider the reflection symmetrisations of $C_l^{v_1}$. We have that
\begin{align*}
(C_l^{v_1})_{\gamma}^{w,+} = C_l^{v_1} \cup C_l^{v_2},
\quad
(C_l^{v_1})_{\gamma}^{w,-} = C_l^{v_1} \cap C_l^{v_2}.
\end{align*}
For the domain of dependence of the reflection symmetrisations $(C_l^{v_1})_{\gamma}^{w,\pm}$, we have that
\begin{align*}
	\calD((C_l^{v_1})_{\gamma}^{w,-}) 
	= 
	(\calD(C_l^{v_1}))_{\gamma}^{w,-} 
	&
	= 
	(\calD(C_l^{v_1})) \cap (\calD(C_l^{v_2}))
\\
	&
	=
	J^+(o) \cap I^-(q_1) \cap I^-(q_2),
\\
	\calD((C_l^{v_1})_{\gamma}^{w,+}) 
	\supsetneq
	(\calD(C_l^{v_1}))_{\gamma}^{w,+} 
	&
	= 
	(\calD(C_l^{v_1})) \cup (\calD(C_l^{v_2}))
\\
	&
	=
	J^+(o) \cap (I^-(q_1) \cup I^-(q_2)),
\end{align*}
Therefore for the volume of the domain of dependence, we have that
\begin{align*}
|\calD(C_l^{v_1})_{\gamma}^{w,-}) | = 2 |\calD((C_l^{v_1}) \cap H_w^- |,
\quad
|\calD((C_l^{v_1})_{\gamma}^{w,+}) | > 2 |\calD(C_l^{v_1}) \cap H_w^+ |,
\end{align*}

\item
The intersection of the strict future boundaries $\uC_{q_1} \cap \uC_{q_2} \cap I^+(o) \subset H_w$ is a subset of the domain of dependence $\calD((C_l^{v_1})_{\gamma}^{w,+}) $ of $(C_l^{v_1})_{\gamma}^{w,+}$. Let $q$ be a point in the intersection of the strict future boundaries $\uC_{q_1} \cap \uC_{q_2} \cap I^+(o)$. Then for $\delta>0$ sufficiently small, we have that $q+ \delta \overrightarrow{oq} \in \calD((C_l^{v_1})_{\gamma}^{w,+})$.
\end{enumerate}
\end{lemma}
\begin{proof}
The proof is straightforward by figure \ref{fig 18}. We just mention the proof of part \textit{d.}. Let $p_a, a=1,2$ be the intersection of the line $\overline{qq_a}$ with $C_0$. Note that 
\begin{align*}
(\uC_q \cap C_0)\setminus \{p_1, p_2\} \subset C_l^{v_1} \cap C_l^{v_2} \subset C_l^{v_1} \cup C_l^{v_2}.
\end{align*}
Moreover since $p_1 \in C_l^{v_2}$ and $p_2 \in C_l^{v_1}$, we have $\{p_1, p_2\} \subset C_l^{v_1} \cup C_l^{v_2}$. Therefore
\begin{align*}
\uC_q \cap C_0 \subset C_l^{v_1} \cup C_l^{v_2}
\quad
\Rightarrow
\quad
q\in \calD ( C_l^{v_1} \cup C_l^{v_2} ) = \calD((C_l^{v_1})_{\gamma}^{w,+}).
\end{align*}
Since $\calD((C_l^{v_1})_{\gamma}^{w,+})$ is open in $I^+(o)$, $q+ \delta \overrightarrow{oq} \in \calD((C_l^{v_1})_{\gamma}^{w,+})$ for sufficiently small $\delta$.
\end{proof}

We state the equal perimeter separation hyperplane noncrossing property of the null generator and prove it by the above lemma.
\begin{proposition}\label{prop 5.8}
Let $C_f$ be an extremal open finite lightcone of the volume functional $\calV$ in $F_1$. Let $q$ be a point in the strict future boundary $\partial_+ \calD(C_f)$ of the domain of dependence $\calD(C_f)$. Then $\uC_q\cap S_f$, the intersection of $\uC_q$ with the lower envelope $S_f$ of $C_f$, is nonempty. Let $p \in \uC_q\cap S_f$. We have that the null segment $\overline{pq}$ lies in the strict future boundary $\partial_+ \calD(C_f)$, thus $\overline{pq}$ is a null generator of $\partial_+ \calD(C_f)$. Every null generator of $\partial_+ \calD(C_f)$ can be obtained by the above construction.

Let $H_w$ be an equal perimeter separation hyperplane. Then the null generator $\overline{pq}$ will not cross $H_w$, i.e. either $\overline{pq} \cap H_w^+ = \emptyset$ or $\overline{pq} \cap H_w^- = \emptyset$.
\end{proposition}
\begin{proof}
We prove the proposition by the method of contradiction. Assume that $\overline{pq}$ crosses $H_w$.

Without loss of generality, suppose $q \in H_w^+$ and $p \in H_w^-$. Let $q_w$ be the intersection of $\overline{pq}$ with $H_w$. Let $\gamma$ be the reflection about $H_w$. Apply lemma \ref{lem 5.7} to the points $q, \gamma(q)$. Let $v=\overrightarrow{op}$ and $|v| =2l$. Then considering the reflection symmetrisation $(C_l^{v})_{\gamma}^{w,+}$ of the open finite ligthcone $C_l^{v}$, we have that $q_w \in \calD((C_l^{v})_{\gamma}^{w,+})$. Since
\begin{align*}
C_l^v \subset C_f
\quad
\Rightarrow
\quad
(C_l^{v})_{\gamma}^{w,+} \subset (C_f)_{\gamma}^{w,+},
\end{align*}
then
\begin{align*}
q_w \in \calD((C_l^{v})_{\gamma}^{w,+})
\quad
\Rightarrow
\quad
q_w \in \calD((C_f)_{\gamma}^{w,+}).
\end{align*}
This implies that $\calD((C_f)_{\gamma}^{w,+})$ contains the reflection symmetrisation $(\calD(C_f))_{\gamma}^{w,+}$ as a proper subset. Note that
\begin{align*}
	&
	(\calD(C_f))_{\gamma}^{w,+}
	=
	(\calD(C_f) \cap H_w^+) 
	\ \cup\ 
	(\calD(C_f) \cap H) 
	\ \cup\
	(\gamma(\calD(C_f)) \cap H_w^-)
\\
	\Rightarrow
	\quad
	&
	|(\calD(C_f))_{\gamma}^{w,+}| 
	= 
	2 |(\calD(C_f) \cap H_w^+) | 
	= 
	V_1,
\end{align*}
by lemma \ref{lem 5.5}. Thus we arrive at the contradiction that $(C_f)_{\gamma}^{w,+} \in F_1$ while
\begin{align*}
\calV((C_f)_{\gamma}^{w,+})
= 
|\calD((C_f)_{\gamma}^{w,+})|
>
|(\calD(C_f))_{\gamma}^{w,+}|
=
V_1.
\end{align*}
Then the assumption is false, thus the proposition is true.
\end{proof}

\subsection{$(\mathbb{Z}_2)^{\times n}$-reflection symmetric extremal finite lightcone}\label{sec 5.3}
We introduce a kind of open finite lightcones which has $n$ orthogonal equal perimeter separation hyperplanes.
\begin{definition}\label{def 5.9}

\begin{enumerate}[label=\alph*.]
\item
Let $C_f$ be an open finite lightcone. Suppose that there is an orthogonal set of $n$ spacelike vectors $\{w_1, \cdots, w_n \}$ such that $C_f$ is symmetric with respect to the reflection about each hyperplane $H_{w_i}$, $i=1,\cdots,n$. Let $\gamma_i$ be the reflection about $H_{w_i}$, then $C_f$ is invariant under the group action of $(\mathbb{Z}_2)^{\times n} = \{ \rmI, \gamma_1\} \times \cdots \times \{ \rmI, \gamma_n\}$. We call $C_f$ a $(\mathbb{Z}_2)^{\times n}$-reflection symmetric open finite lightcone.

\item
We call a set $\calE$ being $(\mathbb{Z}_2)^{\times n}$-reflection symmetric if there exists an orthogonal set of $n$ spacelike vectors $\{w_1, \cdots, w_n \}$, such that $\calE$ is invariant under the reflections $\{ \gamma_1,\cdots, \gamma_n\}$ as in a..
\end{enumerate}
\end{definition}
We show that all equal perimeter separation hyperplanes of a $(\mathbb{Z}_2)^{\times n}$-reflection symmetric open finite lightcone pass through a timelike line.
\begin{lemma}\label{lem 5.10}
Let $C_f$ be $(\mathbb{Z}_2)^{\times n}$-reflection symmetric and $\{w_1, \cdots, w_n \}$ be the corresponding orthogonal set of spacelike vectors in definition \ref{def 5.9}. 
\begin{enumerate}[label=\alph*.]
\item
Let $v$ be the future-directed timelike vector orthogonal to $w_1, \cdots, w_n$. Introduce the line $l_v$ through the origin $o$ in the direction of $v$. The line inverse $\psi_v$ about the line $l_v$ is defined by
\begin{align*}
\psi_v(w)
=
2\langle w, \frac{v}{|v|}\rangle \frac{v}{|v|} -w.
\end{align*}
We have that $\psi_v = \gamma_1 \circ \cdots \circ \gamma_n$. $C_f$ is symmetric about the line $l_v$. We call $l_v$ the axis of reflection symmetry of $C_f$.

\item
$H$ is an equal perimeter separation hyperplane of $C_f$ if and only if $H$ passes through $l_v$.

\end{enumerate}
\end{lemma}
\begin{proof}
A vector $w$ can be expressed by the formula
\begin{align*}
w 
= 
\langle w, \frac{v}{|v|}\rangle \frac{v}{|v|} 
+ 
\langle w, \frac{w_1}{|w_1|}\rangle \frac{w_1}{|w_1|} 
+ 
\cdots 
+ 
\langle w, \frac{w_n}{|w_n|}\rangle \frac{w_n}{|w_n|}.
\end{align*}
Then part \textit{a.} follows from the above formula and direct calculations. Part \textit{b.} follows from \textit{a.}.
\end{proof}

Applying the equal perimeter separation hyperplane noncrossing property to an $(\mathbb{Z}_2)^{\times n}$-reflection symmetric extremal open finite lightcone, we shall show that the extremal open finite lightcone must be the spacelike hyperplane truncated open finite lightcone $C_1^v$ where $v$ is the future direction of the axis of reflection symmetry.
\begin{proposition}\label{prop 5.11}
Let $C_f$ be $(\mathbb{Z}_2)^{\times n}$-reflection symmetric. Let $v$ be the future-directed timelike vector that $l_v$ is the axis of reflection symmetry of $C_f$. Suppose that $C_f$ is an extremal open finite lightcone of the volume functional $\calV$ in the set $F_1$, then $C_f$ is the spacelike hyperplane truncated open finite lightcone $C_1^v$.
\end{proposition}
\begin{proof}
We prove the proposition by three steps.

\vspace{5pt}
\noindent{\underline{\textsc{Step 1.}} }
Let $P$ be a $2$-dimensional plane passing through $l_v$. Assume $P$ intersects with the lower envelope $S_f$ of the open finite lightcone $C_f$ at two points $p_1$, $p_2$. Let $q \notin l_v$ be a point of the strict future boundary $\partial_+ \calD(C_f)$ inside $P$. We claim that $\uC_q \cap S_f \subset \{ q_1, q_2\}$. As a corollary of the claim, either the segment $\overline{p_1 q}$ or $\overline{p_2 q}$ lies in the closure of the strict future boundary $\overline{\partial_+ \calD(C_f) \cap P}$ in $P$.

This claim is proved as follows. Since $q \in \partial_+ \calD(C_f)$, then $\uC_q \cap S_f \neq \emptyset$. If $p \in \uC_q \cap S_f$, then the segment $\overline{pq}$ is a null generator of the strict future boundary $\partial_+ \calD(C_f)$. The points $p$, $q$ and the axis of reflection symmetry $l_v$ must lie in one $2$-dimensional plane, otherwise we can find a timelike hyperplane $H_w$ containing $l_v$ but separating $p$ and $q$ in two sides, for example $p\in H_w^+$ while $q\in H_w^-$. Such a timelike hyperplane $H_w$ does not exist since it contradicts with proposition \ref{prop 5.8}, the equal perimeter separation hyperplane noncrossing property of the null generator $\overline{pq}$. Hence $p$ lies in the $2$-dimensional plane spanned by $q$ and $l_v$, which is $P$ since $q \notin l_v$. Thus $p$ is either $p_1$ or $p_2$.

\vspace{5pt}
\noindent{\underline{\textsc{Step 2.}} }
Let $P$, $p_1$, $p_2$, $q \notin l_v$ be as in \textit{Step 1.}. Without loss of generality, assume that the null segment $\overline{p_1q}$ is a null generator of the strict future boundary $\partial_+ \calD(C_f)$. We claim that $\overline{p_1q} \cap l_v$ must be empty. Otherwise one can easily choose a timelike hyperplane $H_w$ passing through $l_v$ such that $\overline{p_1q}$ crosses $H_w$, which again arrives at a contradiction with proposition \ref{prop 5.8}.

\vspace{5pt}
\noindent{\underline{\textsc{Step 3.}} }
Let $q_v = \partial_+ \calD(C_f) \cap l_v$, the intersection of the strict future boundary and the axis of reflection symmetry. \textit{Step 1.} and \textit{Step 2.} implies that for any $2$-dimensional plane $P$ passing through $l_v$, $\partial_+ \calD(C_f) \cap P = \uC_{q_v} \cap I^+(o) \cap P$. Therefore $\partial_+ \calD(C_f) = \uC_{q_v} \cap I^+(o)$ and $C_f = C_l^v$ for some $l$. Since $C_f \in F_1$,  we obtain that $C_f = C_1^v$.
\end{proof}

\subsection{Identification of case of equality}\label{sec 5.4}
Now we can identity the extremal open finite lightcone cone of the volume function $\calV$ in the set $F_1$. It follows from the following construction of $(\mathbb{Z}^2)^{\times n}$-reflection symmetric open finite lightcones from a given open finite lightcone without any symmetry assumption.
\begin{construction}\label{con 5.12}
Let $C_f$ be an open finite lightcone. We construct a $(\mathbb{Z}^2)^{\times n}$-reflection symmetric open finite lightcone by the following inductive steps.
\begin{enumerate}[label=\arabic*.]
\item
Choose a $2$-dimensional timelike plane $P_1$. Choose $w_1 \in P_1$ such that $H_{w_1}$ is the equal perimeter separation hyperplane of $C_f$. Define $f_1$ as the positive reflection symmetrisation of $f$ about $(w_1, H_{w_1})$, i.e. $f_1 = f_{\gamma_1}^{w_1,+}$ where $\gamma_1$ is the reflection about $H_{w_1}$. $C_{f_1}$ is the positive reflection symmetrisation of $C_f$.

\item
Assume that we obtain the orthogonal set of spacelike vectors $\{ w_1, \cdots, w_k\}$ and the open finite lightcone $C_{f_k}$.

Choose a $2$-dimensional timelike plane $P_{k+1}$ orthogonal to the set $\{ w_1, \cdots, w_k\}$. Then we find a spacelike vector $w_{k+1} \in P_{k+1}$, such that $H_{w_{k+1}}$ is the equal perimeter separation hyperplane of $C_{f_k}$. Adding $w_{k+1}$ to the set $\{ w_1, \cdots, w_k\}$, we obtain a larger orthogonal set of spacelike vectors $\{ w_1, \cdots, w_{k+1}\}$. Define $f_{k+1}$ as the positive reflection symmetrisation of $f_k$ about $(w_{k+1}, H_{w_{k+1}})$. The open finite lightcone $C_{f_{k+1}}$ is the positive reflection symmetrisation of $C_{f_k}$ about $(w_{k+1}, H_{w_{k+1}})$. 
Repeat the above construction until it stops at $k=n$.

\end{enumerate}
We have that
\begin{align*}
f_k = ( \cdots (f_k)_{\gamma_1}^{w_1,+} \cdots)_{\gamma_2}^{w_k,+},
\quad
C_{f_k} = ( \cdots (C_f)_{\gamma_1}^{w_1,+} \cdots)_{\gamma_2}^{w_k,+}.
\end{align*}
An inductive argument shows that $C_{f_k}$ is invariant under the reflection $\gamma_i$ about $H_{w_i}$, $i = 1, \cdots, k$. Thus $C_{f_n}$ is a $(\mathbb{Z}^2)^{\times n}$-reflection symmetric open finite lightcone, invariant under the reflection $\gamma_i$ about $H_{w_i}$, $i = 1, \cdots, n$. Moreover we have that
\begin{align*}
C_f \cap H_{w_1}^+ \cap \cdots \cap H_{w_k}^+ = C_{f_k} \cap H_{w_1}^+ \cap \cdots \cap H_{w_k}^+.
\end{align*}

If $C_f$ is an extremal open finite lightcone of the volume functional $\calV$ in the set $F_1$, then so is $C_{f_i}$, $i=1, \cdots, n$.
\end{construction}

Finally, we answer part \textit{b.} of problem \ref{prob 4.1} by the following result identifying the extremal open finite lightcone.
\begin{proposition}\label{prop 5.13}
Let $C_f$ be an extremal open finite lightcone of the volume functional $\calV$ in the set $F_1$. Then $C_f$ must be a spacelike hyperplane truncated open finite lightcone, which is $C_1^v$ for some future-directed timelike vector $v$.
\end{proposition}
\begin{proof}
For any $p$ in the lower envelope $S_f$ of $C_f$, we can construct a $(\mathbb{Z}_2)^{\times n}$-reflection symmetric open finite lightcone $C_{\barf}$ such that $C_{\barf}$ coincides with $C_f$ in a neighbourhood of $p$. Then by proposition \ref{prop 5.11}, $S_f$ is a spacelike hyperplane section of $C_0$ in a neighbourhood of $p$. By a continuity argument, we obtain that $S_f$ is a spacelike hyperplane section of $C_0$, thus $C_f$ is a spacelike hyperplane truncated open finite lightcone.
\end{proof}

\section{Summary on isoperimetric inequality for domain of dependence of finite lightcone}

We summarise the isoperimetric inequality in the Minkowski spacetime proved previously, and use it to obtain a new geometric inequality in the Euclidean space.

\subsection{Isoperimetric inequality for domain of dependence of finite lightcone}

We state the isoperimetric inequality for the domain of dependence of an open finite lightcone.
\begin{theorem}\label{thm 6.1}
Let $C_f$ be an open finite lightcone in Minkowski spacetime $\mathbb{M}^{n,1}$ and $\calD(C_f)$ be the domain of dependence of $C_f$. Then we have the following inequality
\begin{align*}
\frac{|\calD(C_f)|}{|\calD(C_1)|} 
\leq 
\Big(\frac{\int_{\mathbb{S}^{n-1}} f^{n-1} \dvol_{\circg}}{|\mathbb{S}^{n-1}|} \Big)^{\frac{n+1}{n-1}}.
\end{align*}
If $f$ is piecewise Lipschitz, then $\int_{\mathbb{S}^{n-1}} f^{n-1} \dvol_{\circg}$ is the perimeter of the boundary of $C_f$. The equality is achieved if and only if the open finite lightcone $C_f$ is a spacelike hyperplane truncated open finite lightcone.
\end{theorem}
The theorem follows from propositions \ref{prop 4.21} and \ref{prop 5.13}. We can also prove an analogous isoperimetric inequality for closed finite lightcones, stated as the following theorem.
\begin{theorem}\label{thm 6.2}
Let $\calC_h$ be an open finite lightcone in Minkowski spacetime $\mathbb{M}^{n,1}$ and $\calD(\calC_h)$ be the domain of dependence of $\calC_h$. Then we have the following inequality
\begin{align*}
\frac{|\calD(\calC_f)|}{|\calD(\calC_1)|} \leq \Big(\frac{\int_{\mathbb{S}^{n-1}} h^{n-1} \dvol_{\circg}}{|\mathbb{S}^{n-1}|} \Big)^{\frac{n+1}{n-1}}.
\end{align*}
If $h$ is piecewise Lipschitz, then $\int_{\mathbb{S}^{n-1}} h^{n-1} \dvol_{\circg}$ is the perimeter of the boundary of $\calC_h$. The equality is achieved if and only if the closed finite lightcone $\calC_h$ is a spacelike hyperplane truncated closed finite lightcone.
\end{theorem}
\begin{proof}
The theorem follows from theorem \ref{thm 6.2} and proposition \ref{prop 2.15}. Note that the parameterisation function of the interior $(\calC_h)^{\circ}$ is $h_{\inf}$, which is no more than $h$, hence
\begin{align*}
\int_{\mathbb{S}^{n-1}} (h_{\inf})^{n-1} \dvol_{\circg}
\leq
\int_{\mathbb{S}^{n-1}} h^{n-1} \dvol_{\circg}.
\end{align*}
By proposition \ref{prop 2.15},
\begin{align*}
|\calD(\calC_h)| = |\overline{\calD((\calC_h)^{\circ})}| = |\calD((\calC_h)^{\circ})|,
\end{align*}
where the last identity follows from that $\partial \calD((\calC_h)^{\circ})$ is Lipschitz. Thus the theorem follows from theorem \ref{thm 6.1} for open finite lightcones.
\end{proof}

\subsection{Implication in Euclidean space}
We show that the isoperimetric inequality for the domain of dependence of a finite lightcone has an implication in Euclidean space. Choose a rectangular coordinate system $\{t,x_1, \cdots, x_n\}$ of $(\mathbb{M}^{n,1}, \eta)$ where $\eta = -\ed t^2 + \ed x_1^2 + \cdots + \ed x_n^2$. Then consider the Euclidean metric $g = \ed t^2 + \ed x_1^2 + \cdots + \ed x_n^2$. We have the following correspondence between the geometric concepts in $(\mathbb{M}^{n,1}, \eta)$ and $(\mathbb{E}^{n+1}, g)$.
\begin{enumerate}[label=\alph*.]
\item
Let $p$ be the point with the coordinate $(t(p),x_1(p), \cdots, x_n(p))$. 
\begin{align*}
&
C_p 
= 
\{q = (t,x_1, \cdots, x_n): t\geq t(p), |\overrightarrow{pq}|_{g} = \sqrt{2} (t-t(p)) \},
\\
&
\uC_p 
= 
\{q = (t,x_1, \cdots, x_n): t\leq t(p), |\overrightarrow{pq}|_{g} = \sqrt{2} (t(p)-t) \},
\\
&
I^+(p)
=
\{q = (t,x_1, \cdots, x_n): t> t(p), |\overrightarrow{pq}|_{g} < \sqrt{2} (t-t(p)) \},
\\
&
I^-(p) 
= 
\{q = (t,x_1, \cdots, x_n): t< t(p),  |\overrightarrow{pq}|_{g} < \sqrt{2} (t(p)-t) \},
\\
&
J^+(p) 
= 
\{q = (t,x_1, \cdots, x_n): t \geq t(p), |\overrightarrow{pq}|_{g} \geq \sqrt{2}(t-t(p)) \},
\\
&
J^-(p) 
= 
\{q = (t,x_1, \cdots, x_n): t \leq t(p), |\overrightarrow{pq}|_{g} \leq \sqrt{2}(t(p)-t) \}.
\end{align*}

\item
The definitions \ref{def 2.1.a}, \ref{def 2.1.b} of an open finite lightcone and a closed finite lightcone can be applied to define the Euclidean finite cone of angle $45^{\circ}$ with the axis parallel to $\partial_t$ without any difficulty.

\item
By the correspondence in a., the domain of dependence of $C_f$ can defined by
\begin{align*}
\calD(C_f) = \{ p \in \mathbb{M}^{n,1} = \mathbb{E}^{n+1}: \uC_p \cap C_o \subset C_f \}
\end{align*}
in the Euclidean space.

\item
The volume forms of $\eta$ and $g$ are the same, i.e. $\dvol_{\eta} = \dvol_g$.

\item
Let $f$ be a piecewise Lipschitz continuous function. Then $\eta|_{\partial C_f} \leq g|_{\partial C_f}$ on any tangent plane of $\partial C_f$. Thus $\dvol_{\eta|_{\partial C_f}} \leq \dvol_{g|_{\partial C_f}}$, and $|\partial C_f|_{\eta} \leq |\partial C_f|_g$. The equality is achieved if and only if $f$ is constant.
\end{enumerate}
Therefore we have the following theorem in the Euclidean space.
\begin{theorem}\label{thm 6.3}
Let $f>0$ be a bounded piecewise Lipschitz continuous function and $C_f$ be the open Euclidean finite cone at $o$ of angle $45^{\circ}$ with the axis parallel to $\partial_t$. Then
\begin{align*}
\frac{|\calD(C_f)|}{|\calD(C_1)|} 
\leq 
\Big(\frac{|\partial C_f|_g}{|\mathbb{S}^{n-1}|} \Big)^{\frac{n+1}{n-1}}.
\end{align*}
The equality is achieved if and only if $f$ is constant.
\end{theorem}

\section{Isoperimetric inequality for domain of dependence of set in hyperboloid}\label{sec 7}
The Lorentz polarisation can be applied similarly to the domain of dependence of sets in the hyperboloid $S_{-1} = \{ (t,x) \in \mathbb{M}^{n,1}: -t^2 + |x|^2 = -1 \}$ to prove an analogous isoperimetric inequality, see theorem \ref{thm 7.7}.

\subsection{Lorentz polarisation for domain of dependence of set in hyperboloid}
For the domain of dependence of a set in the hyperboloid $S_{-1}$, we have the following analogous monotonicity proposition as proposition \ref{prop 3.10}.
\begin{proposition}\label{prop 7.1}
Let $H$ be a timelike hyperplane in $\mathbb{M}^{n,1}$ and $v\notin H$ be a timelike vector. Let $\gamma$ be the reflection about $H$. Suppose that $E$ is a set in $(S_{-1}, \eta|_{S_{-1}})$, and $\calD(E)$ is the domain of dependence of $E$. We have that the Lorentz polarisation of $\calD(E)$ is contained in the domain of dependence of the Lorentz polarisation of $E$, i.e.
\begin{align*}
(\calD(E))^{\gamma,v} \subset \calD(E^{\gamma,v}).
\end{align*}
As a corollary, we have that
\begin{align*}
	|\calD(E)| \subset |\calD(E^{\gamma,v})|.
\end{align*}
\end{proposition}

\begin{definition}\label{def 7.2}
Let $v$ be a future-directed timelike vector and $o^v_{-1}$ be the intersection of $S_{-1}$ with the line through $o$ in $v$ direction. Let $E$ be a set in $(S_{-1}, \eta|_{S_{-1}})$. Define $P_E^v$ as the set of sets in $S_{-1}$ obtained by applying the Lorentz polarisation for finite times:
\begin{align*}
	P_E^v
	=
	\{ E^{\gamma_1, \gamma_2,\cdots, \gamma_k, v}: k\in \mathbb{N}, v\notin H_1, \cdots, v \notin H_k\},
\end{align*}
where $E^{\gamma_1, \gamma_2,\cdots, \gamma_k, v}$ is $(\cdots(E^{\gamma_1,v})^{\gamma_2,v}\cdots)^{\gamma_k,v}$.
\end{definition}

Similarly, we can introduce the distance $d_{\calD}$ for sets in $S_{-1}$ by the volume of the symmetric difference of the domains of dependence.
\begin{definition}\label{def 7.2}
Introduce the distance $d_{\calD}$ between two sets in $S_{-1}$ by
\begin{align*}
d_{\calD} ( E_1, E_2 )= | \calD(E_1) \triangle \calD(E_2) |
\end{align*}
Define the set $\overline{P_E^v}$ as the closure of $P_E^v$ under the distance $d_{\calD}$.
\end{definition}
The set $(\overline{P_E^v}, d_{\calD})$ for bounded open $E$ has the similar metric property as $(\overline{P_f^v}, d_{\calD})$.
\begin{proposition}\label{prop 7.4}
$(\overline{P_E^v}, d_{\calD})$ with bounded $E$ is a complete compact metric space. Moreover $\overline{P_E^v}$ is invariant under the Lorentz polarisation about any pair $(H,v)$ where $v\not\in H$.
\end{proposition}
To prove proposition \ref{prop 7.4}, we introduce the parameterisation function for the boundaries $\partial_{\pm} \calD(E) = \partial_{\pm} \calD(E^{\circ}) $,
\begin{align*}
	\calD(E^{\circ})
	=
	\{ (t,x) \in \mathbb{M}^{n,1}: u_-(x) < t < u_+(x), u_-(x) \leq \sqrt{1+|x|^2} \leq u_+(x) \}.
\end{align*}
We sketch the proof of proposition \ref{prop 7.4}.
\begin{proof}
Let $\{ E_k \}$ be a sequence in $(P_E^v, d_{\calD})$. Let $u_{k,\pm}$ be the parameterisation functions of $\partial_{\pm} \calD(E_k)$. Since $E$ is bounded, there exist convergent subsequences $\{ u_{k_i,\pm} \}$ of $\{ u_{k_i,\pm} \}$. Suppose $\baru_{\pm} = \lim_{i \rightarrow +\infty} u_{k_i,\pm}$, then define
\begin{align*}
	\barE
	=
	\{ (t,x) \in S_{-1}: \baru_-(x) < \baru_+(x) \},
	\quad
	\calD
	=
	\{ (t,x) \in \mathbb{M}^{n,1} : \baru_-(x) < t < \baru_+(x) \}.
\end{align*}
By similar arguments in the proof of lemma \ref{lem 4.5}, we can show that
\begin{align*}
	\calD = \calD(\barE).
\end{align*}
Then $\lim_{i \rightarrow +\infty} d_{\calD} (E_{k_i}, \barE) = 0$ and the compactness of $\overline{P_E^v}$ follows.

The proof of the invariance of $\overline{P_E^v}$ under the Lorentz polarisation follows the similar route as in the proof of proposition \ref{prop 4.8}. Let $\barE \in \overline{P_E^v}$ and $\{E_k\}\subset P_E^v$ converges to $\barE$ under the distance $d_{\calD}$. Without loss of generality, we can assume that $\barE$ is open, since $d_{\calD}(\barE, \barE^{\circ}) = 0$. Consider the sequence $\{(E_k)^{\gamma,v}\}$ and $\calE_k = (\calD(E_k))^{\gamma,v}$. Introduce the parameterisation functions $\nu_{k,\pm}$ of $\partial_{\pm} \calE_k$,
\begin{align*}
	\calE_k
	=
	\{ (t,x) \in \mathbb{M}^{n,1}: \nu_{k,-}(x) < t < \nu_{k,+}(x), \nu_{k,-}(x) \leq \sqrt{1+|x|^2} \leq \nu_{k,+}(x)  \}
\end{align*}
$\{ \calE_k \}$ converges with respect to the symmetric volume difference by proposition \ref{prop 3.4}. Then $\calE = (\calD(\barE))^{\gamma,v}$ is the limit of $\{ \calE_k \}$. Introduce the parameterisation function $\nu_{\pm}$ of $\partial_{\pm} \calE$,
\begin{align*}
	\calE 
	= 
	\{ (t,x) \in \mathbb{M}^{n,1}: \nu_-(x) < t < \nu_+(x), \nu_-(x) \leq \sqrt{1+|x|^2} \leq \nu_+ (x)  \},
\end{align*}
then we have that
\begin{align*}
	\nu_{\pm} = \lim_{k\rightarrow +\infty} \nu_{k,\pm}.
\end{align*}

Define $\barE_{\nu} = \calE \cap S_{-1}$, we have
\begin{align*}
	\barE_{\nu}
	=
	\{ (t,x) \in S_{-1}: \nu_- (x) < \nu_+ (x) \}.
\end{align*}
We can show that
\begin{align*}
	\lim_{k\rightarrow +\infty} d_{\calD}((\barE_k)^{\gamma,v}, \barE_{\nu}) =0,
	\quad
	\barE_{\nu} = \barE^{\gamma,v}.
\end{align*}
The proof is sketched in the following.

\vspace{5pt}
\noindent{\underline{\textsc{Proof of $\lim_{k\rightarrow +\infty} d_{\calD}((\barE_k)^{\gamma,v}, \barE_{\nu}) =0$}.} }
Introduce $\calD_k = \calD((\barE_k)^{\gamma,v}) = \calD(\calE_k)$, $\calD = \calD(\barE_{\nu}) = \calD(\calE)$ and their parameterisation functions $\baru_{k,\pm}$, $\baru_{\pm}$,
\begin{align*}
	&
	\calD_k
	= 
	\{ (t,x) \in \mathbb{M}^{n,1}: 
	\baru_{k,-}(x) < t < \baru_{k,+}(x), 
	\baru_{k,-}(x) \leq \sqrt{1+|x|^2} \leq \baru_{k,+} (x)  \},
	\\
	&
	\calD
	= 
	\{ (t,x) \in \mathbb{M}^{n,1}: 
	\baru_-(x) < t < \baru_+(x), 
	\baru_-(x) \leq \sqrt{1+|x|^2} \leq \baru_+ (x)  \}.
\end{align*}
By the similar argument as in \textit{Step 2.} of the proof of proposition \ref{prop 4.8}, we can show that
\begin{align*}
	\baru_{\pm} = \lim_{k\rightarrow +\infty} \baru_{k,\pm},
\end{align*}
which implies that $\lim_{k\rightarrow +\infty} d_{\calD}((\barE_k)^{\gamma,v}, \barE_{\nu}) =\lim_{k\rightarrow +\infty} |\calD_k \triangle \calD|= 0$.

\vspace{5pt}
\noindent{\underline{\textsc{Proof of $\barE_{\nu} = \barE^{\gamma,v}$}.} }
It is sufficient to show that $\calD(\barE_{\nu}) = \calD(\barE^{\gamma,v})$.
\begin{enumerate}[label=\textit{\roman*.}]
\item
$\calD(\barE_{\nu}) \subset \calD(\barE^{\gamma,v})$. Note $\barE_{\nu} = \calE \cap S_{-1}$, then
\begin{align*}
	\calD(\barE_{\nu}) = \calD(\calE) = \calD((\calD(\barE))^{\gamma,v}).
\end{align*}
Since $(\calD(\barE))^{\gamma,v} \subset \calD(\barE^{\gamma,v})$, we have $\calD((\calD(\barE))^{\gamma,v}) \subset \calD(\calD(\barE^{\gamma,v})) = \calD(\barE^{\gamma,v})$, hence
\begin{align*}
	\calD(\barE_{\nu})
	\subset
	\calD(\barE^{\gamma,v}).
\end{align*}

\item
$\calD(\barE_{\nu}) \supset \calD(\barE^{\gamma,v})$. Since $\barE^{\gamma,v} \subset (\calD(\barE))^{\gamma,v} = \calE$, then $\calD(\barE^{\gamma,v}) \subset \calD(\calE) = \calD(\barE_{\nu})$.
\end{enumerate}
\end{proof}

Then we can define the extremal set $M\overline{P_E^v} \subset P_E^v$ which maximises the volume of the domain of dependence.
\begin{definition}\label{def 7.5}
Let $E \subset S_{-1}$ be bounded. Define the set $M\overline{P_E^v}$ as the set of $\barE \in \overline{P_E^v}$ which maximises $|\calD(\barE)|$ in $\overline{P_E^v}$, i.e.
\begin{align*}
	M\overline{P_E^v} 
	= 
	\{ \barE' \in \overline{P_E^v} : 
	|\calD(\barE)| = \max_{\barE \in \overline{P_E^v}} |\calD(\barE)| \}.
\end{align*}
\end{definition}
We have the following metric property for $M\overline{P_E^v}$ similar to lemma \ref{lem 4.11}.
\begin{lemma}\label{lem 7.6}
The set $M\overline{P_E^v}$ is a closed set under the distance $d_{\calD}$ and is invariant under the Lorentz polarisation about any pair $(H,v)$. For any $\barE \in M\overline{P_E^v}$, we have that
\begin{align*}
\calD( \barE^{\gamma,v}) = (\calD(\barE))^{\gamma,v}.
\end{align*}
\end{lemma}

\subsection{Statement of isoperimetric inequality and proof}
By the relation between the Lorentz polarisation in Minkowski spacetime and the polarisation on the hyperboloid $S_{-1}$, we can obtain the following isoperemetric inequality for the domain of dependence of sets in $S_{-1}$.
\begin{theorem}
\label{thm 7.7}
We state the results in two cases.
\begin{enumerate}[label=\alph*.]
\item
Let $E \subset S_{-1}$ be a set of finite perimeter. The following isoperimetric inequality holds that
\begin{align}
	\frac{|\calD(E)|}{2\omega_n/(n+1)} 
	\leq
	\Big(\frac{P(E)}{n \omega_n} \Big)^{\frac{n+1}{n-1}}.
	\label{eqn 7.1}
\end{align}
where $P(E)$ is the perimeter of $E$. The equality is achieved when $E^{\circ}$ is a geodesic ball $B$ in $(S_{-1}, \eta|_{S_{-1}})$ and $|E \setminus B | =0$.

\item
Let $E \subset S_{-1}$ be a set with finite measure. The following inequality holds that
\begin{align}
	|\calD(E)|
	\leq
	|\calD(B_{r_{|E|}})|,
	\label{eqn 7.2}
\end{align}
where $B_{r_{|E|}}$ is a geodesic ball in $(S_{-1}, \eta|_{S_{-1}})$ with the same measure as $E$. The equality is achieved when $E^{\circ}$ is a geodesic ball $B$ in $(S_{-1}, \eta|_{S_{-1}})$ and $|E\setminus B| =0$.
\end{enumerate}
\end{theorem}
\begin{proof}[{\bf Proof of isoperimetric inequality \eqref{eqn 7.1}}]
Inequality \eqref{eqn 7.1} follows from inequality \eqref{eqn 7.2} and the isoperimetric inequality in the hyperbolic space.
\end{proof}

\begin{proof}[{\bf Proof of inequality \eqref{eqn 7.2} for bounded set}]
We consider the case of bounded $E$ first. Choose $v=e_0=(1,0,\cdots,0)$. Let $B$ be the geodesic ball of $S_{-1}$ centred at $e_0$ with $|B| = |E|$. Then there exists a convergent sequence of sets $E_k$ in $P_E^v$ converging to a set $\barE\in \overline{P_E^v}$  with respect to the distance $d_{\calD}$,
\begin{align*}
	\lim_{k \rightarrow + \infty} |E_k \triangle B| = 0.
\end{align*}
Then we have that
\begin{align*}
	\barE^{\circ} \subset B.
\end{align*}
This is because that if $x \in \barE^{\circ}$, then there exists a small number $\delta>0$, such that the geodesic ball of radius $\delta$ entered at $x$ $B_{\delta}(x)\subset \barE$, which implies that
\begin{align*}
	\calD(B_{\delta}(x)) \subset \calD(\barE).
\end{align*}
Then for sufficiently large $k$, we have $\calD(B_{\delta/2}(x)) \subset \calD(E_k)$, which implies $B_{\delta/2}(x) \subset E_k$. Hence 
\begin{align*}
	|B_{\delta/2}(x) \setminus B| \leq \lim_{k\rightarrow} |E_k \setminus B| = 0,
\end{align*}
which implies that $x\in B_{\delta/2}(x)  \subset B$. Thus $\barE^{\circ} \subset B$.

Applying the monotonicity property of the volume of the domain of dependence under the Lorentz polarisation, proposition \ref{prop 7.1}, to the sequence $\{E_k\}$, we have that
\begin{align*}
	|\calD(E)| \leq \lim_{k\rightarrow + \infty} |\calD(E_k)| = |\calD(\barE)| = |\calD(\barE^{\circ})| \leq |\calD(B)|.
\end{align*}
Then inequality \eqref{eqn 7.2} follows.
\end{proof}

\begin{proof}[{\bf Proof of inequality \eqref{eqn 7.2} for general case}]
Now we consider the general case. The following estimate is useful.

\vspace{5pt}
\noindent{\underline{\textsc{Claim 1.}} } 
\textit{Let $E$ be a set of finite perimeter in $S_{-1}$. Let $\delta_{|E|}$ be the radius of geodesic ball in $(S_{-1},\eta|_{S_{-1}})$ having the volume $|E|$. Then
\begin{align*}
	|\calD(E)|
	\leq
	\frac{|E|}{n+1} [(\cosh \delta_{|E|} + \sinh \delta_{|E|})^{n+1} - (\cosh \delta_{|E|} - \sinh \delta_{|E|})^{n+1}]
\end{align*}}

We prove the inequality assuming Claim 1. and verifying it later. For any $\epsilon >0$, there exists a sufficiently large $R$ such that 
\begin{align*}
	| E\setminus B_R(e_0)| < \epsilon, 
	\quad
	\int_{\partial B_{2R}(e_0)} (\mathbf{1}_E)^- \ed \sigma < \epsilon,
\end{align*}
where $(\mathbf{1}_E)_-$ is the trace of $\mathbf{1}_E$ on $\partial B_{2R}(e_0)$ from inward\footnote{See remark 2.14 in \cite{G84}.},
and $4r_{\epsilon} <R$ where
\begin{align*}
	|B_{r_{\epsilon}}(e_0)| = \epsilon.
\end{align*}

\vspace{5pt}
\noindent{\underline{\textsc{Claim 2.}} }
\textit{Define 
$E_{2R} = E \cap B_{2R}(e_0)$, $E_R^c = E \setminus \overline{B_{R}(e_0)}$,
then we have
\begin{align*}
	\calD(E) \setminus \calD(E_R^c)
	\subset
	\calD(E_{2R}).
\end{align*}}

Assuming Claim 2., we have that
\begin{align*}
	P(E_{2R})
	\leq
	P(E) 
	- P(E_{2R}^c) 
	+ \int_{\partial B_{2R}(e_0)} (\mathbf{1}_E)^- \ed \sigma
	\leq
	P(E) + \epsilon,
\end{align*}
and
\begin{align*}
	|\calD(E)|
	&
	\leq
	|\calD(E_{2R})| + |\calD(E_R^c)|
\\
	&
	\leq
	|\calD(E_{2R})| 
	+ \frac{|\epsilon|}{n+1} [(\cosh \delta_{\epsilon} + \sinh \delta_{\epsilon})^{n+1} - (\cosh \delta_{\epsilon} - \sinh \delta_{\epsilon})^{n+1}].
\end{align*}
Then applying the inequality to $E_{2R}$ and $\calD(E_{2R})$ and taking $\epsilon \rightarrow 0^+$, we obtain the inequality for $E$.
\end{proof}

\begin{proof}[Proof of \textsc{Claim 1}]
First, we have
\begin{align*}
	\calD(E) 
	\subset
	\{ l v: v\in E , l >0  \}.
\end{align*}
Let $v \in S_{-1}$. Then for $l \in (0,\cosh \delta_{|E|} - \sinh \delta_{|E|}]$,
\begin{align*}
	B_{\delta_{|E|}}(v)
	\subset
	I^+(lv) \cap S_{-1},
\end{align*}
and for $l \in [\cosh \delta_{|E|} + \sinh \delta_{|E|}, +\infty)$,
\begin{align*}
	B_{\delta_{|E|}}(v)
	\subset
	I^-(lv) \cap S_{-1}.
\end{align*}
\begin{figure}[h]
\centering
\begin{tikzpicture}
\draw (-5,5) -- (0,0) -- (5,5);
\draw[->,line width=0.3mm,blue] (0,0) -- (0,3);
\node[below right] at (0,3) {$v$};
\draw[domain=-1:1,smooth,variable=\t]
	plot ({3*sinh(\t)},{3*cosh(\t)});
\draw[domain=-.8:.8,smooth,line width=0.75mm,red,variable=\t]
	plot ({3*sinh(\t)},{3*cosh(\t)});
\node[above] at ({3*sinh(0.4)},{3*cosh(0.4)}) {\tiny $B_{\delta_{|E|}}(v)$};
\draw ({3*sinh(0.8)},{3*cosh(0.8)}) -- (0,{3*cosh(0.8) - 3*sinh(0.8)}) -- ({-3*sinh(0.8)},{3*cosh(0.8)});
\node[right] at (0,{3*cosh(0.8) - 3*sinh(0.8)}) {\tiny$[\cosh(\delta_{|E|}) - \sinh(\delta_{|E|})] v $};
\draw ({3*sinh(0.8)},{3*cosh(0.8)}) -- (0,{3*cosh(0.8) + 3*sinh(0.8)}) -- ({-3*sinh(0.8)},{3*cosh(0.8)});
\node[right] at (0,{3*cosh(0.8) + 3*sinh(0.8)}) {\tiny$[\cosh(\delta_{|E|}) + \sinh(\delta_{|E|})] v $};
\draw (0,0) -- (0,{3*cosh(0.8) + 3*sinh(0.8)});
\end{tikzpicture}
\caption{$B_{\delta_{|E|}}(v) \subset I_{\pm}(lv) \cap S_{-1}$.}
\label{fig 19}
\end{figure}
See figure \ref{fig 19}. Therefore for $l \in (0,\cosh \delta_{|E|} - \sinh \delta_{|E|}] \cup [\cosh \delta_{|E|} + \sinh \delta_{|E|}, +\infty)$,
\begin{align*}
	lv \notin \calD.
\end{align*}
This is because $|E| = |B_{\delta_{|E|}}(v)| \leq |I_{\pm}(lv) \cap S_{-1}|$, which implies $I_{\pm}(lv) \cap S_{-1} \not\subset E$. Hence we obtain that
\begin{align*}
	\calD(E)
	\subset
	\{ l v: v\in E , l \in (\cosh \delta_{|E|} - \sinh \delta_{|E|}, \cosh \delta_{|E|} + \sinh \delta_{|E|})  \}.
\end{align*}
Then the inequality of $|\calD(E)|$ follows.
\end{proof}

\begin{proof}[Proof of \textsc{Claim 2}]
It is equivalent to show $\calD(E) \setminus [\calD(E_R^c) \cup \calD(E_{2R})] = \emptyset$. Suppose $x \in \calD(E) \setminus [\calD(E_R^c) \cup \calD(E_{2R})]$, then
\begin{align*}
	[ I_{\pm}(x) \cap S_{-1} ] \cap \overline{B_R(e_0)} \neq \emptyset,
	\quad
	[ I_{\pm}(x) \cap S_{-1} ] \setminus B_{2R}(e_0) \neq \emptyset.
\end{align*}
\begin{figure}[h]
\centering
\begin{tikzpicture}
\draw (0,0) circle [radius=2];
\draw (0,0) -- ({2*cos(120)},{2*sin(120)});
\node[right] at ({1*cos(120)},{1*sin(120)}) {$R$};
\draw (0,0) circle [radius=4];
\draw (0,0) -- ({4*cos(150)},{4*sin(150)});
\node[above right] at ({3*cos(150)},{3*sin(150)}) {$2R$};
\draw ({3.2*cos(30)},{3.2*sin(30)}) circle [radius=1.3];
\node[right] at ({3.2*cos(30)+1.3*cos(30)},{3.2*sin(30)+1.3*sin(30)}) {\small $I_{\pm}(x) \cap S_{-1} $};
\draw ({3*cos(30)},{3*sin(30)}) circle [radius=0.8];
\draw ({3*cos(30)},{3*sin(30)}) node[right] {\small $y$} -- ({3*cos(30)+.8*cos(120)},{3*sin(30)+.8*sin(120)});
\node[right] at ({3*cos(30)+.4*cos(120)},{3*sin(30)+.4*sin(120)}) {\small $r_{\epsilon}$};
\end{tikzpicture}
\caption{$B_{r_{\epsilon}}(y) \subset [ I_{\pm}(x) \cap S_{-1} ] \setminus \overline{B_R (e_0)}$.}
\label{fig 20}
\end{figure}
See figure \ref{fig 20}. Then there exists a point $y \in S_{-1}$ such that
\begin{align*}
	&
	B_{r_{\epsilon}}(y)
	\subset
	[ I_{\pm}(x) \cap S_{-1} ] \setminus \overline{B_R (e_0)}
	\subset
	E \setminus B_R (e_0)
\\
	\Rightarrow
	\quad
	&
	\epsilon
	=
	|B_{r_{\epsilon}}(y)|
	\leq
	|E \setminus \overline{B_R (e_0)}|
	<
	\epsilon.
\end{align*}
Thus the above contradiction implies that such $x$ does not exist, hence $\calD(E) \setminus [\calD(E_R^c) \cup \calD(E_{2R})] = \emptyset$.
\end{proof}

\begin{proof}[{\bf Proof of the case of equality}]
The case of equality in inequality \eqref{eqn 7.1} follows from the one in inequality \eqref{eqn 7.2} and the isoperimetric inequality in the hyperbolic space.

For $E \in S_{-1}$, we define the equal measure separation hyperplane by the hyperplane separating the set to two parts with the same measure. Then the case of equality follows the same route as in section \ref{sec 5} by the equal measure separation hyperplane taking the role of the equal perimeter separation hyperplane.
\end{proof}

\section{Isoperimetric inequality for achronal hypersurface with boundary in lightcone}\label{sec 8}
In this section, we apply the method to Lorentz polarisation to show that the spacelike hyperplane has the maximal area among all the achronal hypersurfaces with the boundary of the same perimeter in the lightcone.
We enlarge the set of achronal hypersurfaces for the comparison of the area such that the suitable convergence limit stays inside the set.
\begin{theorem}\label{thm 8.1}
Let $\Sigma \subset I^+(o)$ be a closed achronal hypersurface in $I^+(o)$. Assume that $\Sigma$ is contained in the domain of dependence of an open finite lightcone $C_f$. Then we have the following isoperimetric inequality
\begin{align}
	\frac{|\Sigma|}{\omega_n}
	\leq
	\Big( \frac{\int_{\mathbb{S}^{n-1}} f^{n-1} \dvol_{\circg}}{n \omega_n} \Big)^{\frac{n}{n-1}}.
	\label{eqn 8.1}
\end{align}
The equality is achieved if and only if $C_f$ is a spacelike hyperplane $H$ truncated open finite lightcone and $\Sigma= H\cap I^+(o)$.
\end{theorem}

\subsection{Lorentz polarisation for strict future boundary of causal convex set}
We show that the Lorentz polarisation preserves the area of a closed achronal hypersurface.
\begin{lemma}\label{lem 8.2}
Let $\calE$ be an open causal convex set, i.e. $\calE$ is open and the causal diamond $J(p,q) \subset \calE$ for all $p,q\in \calE$. Assume that $\calE \subset I^+(o)$ and the strict past boundary $\partial_- \calE \subset C_0$. Let $v$ be future timelike vector and $H$ be a timelike hyperplane through the origin $o$ not containing $v$. Then the Lorentz polarisation $\calE^{\gamma,v}$ of $\calE$ about the pair $(H,v)$ is also an open causal convex set. Moreover the strict past boundary is transformed to the strict past boundary of the Lorentz polarisation, i.e.
\begin{align*}
	(\partial_- \calE)^{\gamma,v}
	=
	\partial_- (\calE^{\gamma,v}),
\end{align*}
and the strict future boundaries of $\calE$ and $\calE^{\gamma,v}$ have the same area,
\begin{align*}
	|\partial_+ (\calE^{\gamma,v})|
	=
	|\partial_+ \calE|.
\end{align*}
\end{lemma}
\begin{proof}
\noindent{\underline{\textsc{$\calE^{\gamma,v}$ is open causal convex.}} }
We show that if $p,q \in \calE^{\gamma,v}$, then $J(p,q) \subset \calE^{\gamma,v}$. Note that $J(p,q) \subset J(o,q)$, thus it is sufficient to show $J(0,q) \subset \calE^{\gamma,v}$. Clearly that $\calE^{\gamma,v}$ is closed. The argument for the causal convexity of $\calE^{\gamma,v}$ is similar to the proof of proposition \ref{prop 3.10}.
\begin{enumerate}[label=\textit{\roman*}.]
\item
If $q\in \calE^{\gamma,v} \cap H_-$, then $q\in \calE \cap H_-$ and $\gamma(q) \in \calE\cap H_+$, hence $J(o,q) \subset \calE$ and $J(o,\gamma(q)) \subset \calE$. Thus $J(o,q) \subset \calE^{\gamma,v}$.

\item
If $q\in \calE^{\gamma,v} \cap H_+$, then $q\in \calE \cap H_-$ or $\gamma(q) \in \calE\cap H_+$, hence $J(o,q) \subset \calE$ or $J(o,\gamma(q)) \subset \calE$. Thus $J(o,q) = J(o,q)^{\gamma,v}  \subset \calE^{\gamma,v}$.

\item
If $q\in \calE^{\gamma,v} \cap H$, then $q\in \calE \cap H$, hence $J(o,q) \subset \calE$. Thus $J(o,q) = J(o,q)^{\gamma,v}  \subset \calE^{\gamma,v}$.
\end{enumerate}

\vspace{5pt}
\noindent{\underline{\textsc{$(\partial_- \calE)^{\gamma,v} = \partial_- (\calE^{\gamma,v})$.}} }

\begin{enumerate}[label=\textit{\alph*}.]
\item
$(\partial_- \calE)^{\gamma,v} \subset \partial_- ( \calE^{\gamma,v} )$: if $q\in (\partial_- \calE)^{\gamma,v}$, then three cases occur.
\begin{enumerate}[label=\textit{a.\roman*}.]
\item
If $q\in (\partial_- \calE)^{\gamma,v} \cap H_-$, then $q \in \partial_- \calE$ and $\gamma(q) \in \partial_- \calE$. Therefore there exists a point $q' \in I^+(q)$ such that $q' \in \calE$ and $\gamma(q') \in \calE$. Thus $J(o,q') \subset \calE$ and $J(o,\gamma(q')) \subset \calE$, which implies that $J(o,q') \subset \calE^{\gamma,v}$. Hence $q \in \partial_- ( \calE^{\gamma,v} )$.

\item
If $q\in (\partial_- \calE)^{\gamma,v} \cap H_+$, then $q \in \partial_- \calE$ or $\gamma(q) \in \partial_- \calE$. Therefore there exists a point $q' \in I^+(q) \cap H_+$ such that $q' \in \calE$ or $\gamma(q') \in \calE$. Thus $J(o,q') \subset \calE$ or $J(o,\gamma(q')) \subset \calE$, which implies that $J(o,q') = (J(o,q'))^{\gamma,v} \subset \calE^{\gamma,v}$. Hence $q \in \partial_- ( \calE^{\gamma,v} )$.

\item
If $q\in (\partial_- \calE)^{\gamma,v} \cap H = \partial_- \calE \cap H$. Therefore there exists a point $q' \in I^+(q) \cap H$ such that $q' \in \calE \cap H$. Thus $J(o,q') \subset \calE$, which implies that $J(o,q') = (J(o,q'))^{\gamma,v} \subset \calE^{\gamma,v}$. Hence $q \in \partial_- ( \calE^{\gamma,v} )$.
\end{enumerate}

\item
$\partial_- ( \calE^{\gamma,v} ) \subset (\partial_- \calE)^{\gamma,v}$: if $q \in \partial_- ( \calE^{\gamma,v} )$, then three cases occur.
\begin{enumerate}[label=\textit{b.\roman*}.]
\item
If $q \in \partial_- ( \calE^{\gamma,v} ) \cap H_-$, then there exists a point $q' \in I^+(q) \cap H_-$ such that $q' \in \calE^{\gamma,v} \cap H_-$. Thus $q' \in \calE$ and $\gamma(q') \in \calE$, which implies that $J(o,q') \subset \calE$ and $J(o,\gamma(q')) \subset \calE$. Hence $q \in \partial_- \calE$ and $\gamma(q) \in \partial_- \calE$. Therefore $q \in (\partial_- \calE)^{\gamma,v}$.

\item
If $q \in \partial_- ( \calE^{\gamma,v} ) \cap H_+$, then there exists a point $q' \in I^+(q) \cap H_+$ such that $q' \in \calE^{\gamma,v} \cap H_+$. Thus $q' \in \calE$ or $\gamma(q') \in \calE$, which implies that $J(o,q') \subset \calE$ or $J(o,\gamma(q')) \subset \calE$. Hence $q \in \partial_- \calE$ or $\gamma(q) \in \partial_- \calE$. Therefore $q = {q}^{\gamma,v} \in (\partial_- \calE)^{\gamma,v}$.

\item
If $q \in \partial_- ( \calE^{\gamma,v} ) \cap H$, then there exists a point $q' \in I^+(q) \cap H$ such that $q' \in \calE^{\gamma,v} \cap H$. Thus $q' \in \calE \cap H$, which implies that $J(o,q') \subset \calE$. Hence $q \in \partial_- \calE \cap H$. Therefore $q \in (\partial_- \calE)^{\gamma,v}$.
\end{enumerate}
\end{enumerate}

\vspace{5pt}
\noindent{\underline{\textsc{$\partial_+ ( \calE^{\gamma,v} ) = \partial_+ \calE$.}} }
Introduce the parameterisation function $\nu$ of $\partial_+ \calE$,
\begin{align*}
	\calE
	=
	\{ (t,x) \in \mathbb{M}^{n,1}: r < t < \nu(x) \}.
\end{align*}
Without loss of generality, we assume that $H=\{x^1=0\}$, $H_+ = \{ x^1 >0 \}$, $H_- = \{ x^1 <0 \}$. Define $\bar{\nu}$ as
\begin{align*}
	\bar{\nu}(x)
	=
	\left\{
	\begin{aligned}
		&
		\max \{ \nu(x^1,\cdots, x^n), \nu(-x^1, x^2, \cdots, x^n) \},
		\quad
		x^1 \geq 0,
	\\
		&
		\min \{ \nu(x^1,\cdots, x^n), \nu(-x^1, x^2, \cdots, x^n) \},
		\quad
		x^1 < 0.
	\end{aligned}
	\right.
\end{align*}
Then $\bar{\nu}$ is the parameterisation function of $(\calE^{\gamma,v})^{\circ}$,
\begin{align*}
	(\calE^{\gamma,v})^{\circ}
	=
	\{ (t,x) \in \mathbb{M}^{n,1}: r < t < \bar{\nu}(x) \}.
\end{align*}

Let $e_1=(1,0,\cdots,0) \in \mathbb{R}^n$ and $K=\{x_1=0\} \subset \mathbb{R}^n$. Let $\sigma$ be the reflection in $\mathbb{R}^n$ about $K$. Then $\bar{\nu} = \nu^{\sigma,e_1}$, the polarisation of $\nu$ about the pair $(K,e_1)$ in $\mathbb{R}^n$. Define $E \subset \mathbb{R}^n$ to be the set $\{x \in E: \nu(x) >r \}$ and $\bar{E} \subset \mathbb{R}^n$ to be the set $\{ x \in \bar{\Omega}: \bar{\nu}(x) >r \}$. Then $\bar{E} = E^{\sigma,e_1}$, the polarisation of $E$ about the pair $(K,e_1)$ in $\mathbb{R}^n$.

Since $\nu$ and $\bar{\nu}=\nu^{\sigma,e_1}$ are both Lipschitz functions with Lipschitz constant $1$, we have that the areas of $\partial_+ \calE$ and $\partial_+ ( \calE^{\gamma,v} )$ are given by
\begin{align*}
	| \partial_+ \calE |
	=
	\int_{\Omega} \sqrt{1 - |\nabla \nu|^2} \ed x,
	\quad
	| \partial_+ (\calE^{\gamma,v}) |
	=
	\int_{\bar{\Omega}} \sqrt{1 - |\nabla \nu^{\sigma,e_1}|^2} \ed x.
\end{align*}
Let $\Omega = \{ x: \nu(x) = \nu( \sigma(x) ) \}$. To show that $| \partial_+ \calE | = | \partial_+  (\calE^{\gamma,v}) |$, it is sufficient to prove that
\begin{align*}
	\int_{\Omega} \sqrt{1 - |\nabla \nu|^2} \ed x
	=
	\int_{\Omega} \sqrt{1 - |\nabla (\nu)^{\sigma,e_1}|^2} \ed x.
\end{align*}
Note that $\nabla \nu(x) = \sigma (\nabla \nu(\sigma(x))) \Rightarrow \nabla \nu (x) = \nabla \nu^{\sigma,e_1}(x)$, therefore the set $\{ x \in \Omega: \nabla \nu (x) \neq \nabla \nu^{\sigma,e_1}(x) \}$ is of measure zero in $\mathbb{R}^n$,
\begin{align*}
	| \{ x \in \Omega: \nabla \nu (x) \neq \nabla \nu^{\sigma,e_1}(x) \} |
	=
	0.
\end{align*}
Hence $| \partial_+ \calE | = | \partial_+ (\calE^{\gamma,v}) |$ follows.
\end{proof}

\subsection{Simple case of spacelike hyperplane truncated lightcone}
Next we show that the spacelike hyperplane maximising the area among all the spacelike hypersurfaces with the same boundary in the lightcone.
\begin{lemma}\label{lem 8.3}
Let $v$ be a future-directed unit timelike vector and $l>0$. Let $H_l^v$ be the spacelike hyperplane orthogonal to $v$ and passing through the point $lv$. Let $C_l^v$ be the open finite lightcone truncated by $H_l^v$ (same notations as in proposition \ref{prop 4.16}) and $D_l^v = H_l^v \cap I^+(o)$. For any archronal hypersurface $\Sigma \subset \calD(C_l^v)$, we have
\begin{align*}
	|\Sigma| \leq |D_l^v|,
\end{align*}
and the equality is achieved if and only if $\Sigma = D_l^v$. In other words,
\begin{align*}
	|D_l^v| = \max\{ |\Sigma|: \Sigma \subset \calD(C_l^v) \textit{ and is achronal} \},
\end{align*}
and the maximal area hypersurface is unique.
\end{lemma}
\begin{proof}
Without loss of generality, we assume that $v=e_0 = (1,0,\cdots,0)$. Any achronal hypersurface $\Sigma \subset \calD(C_l^v)$ can be parameterised as follows
\begin{align*}
	\Sigma
	=
	\{ (t,x) \in \mathbb{M}^{n,1}: t=\nu(x), x\in E \subset D_l \},
\end{align*}
where $D_l$ is the disk of radius $l$. Then the area of $\Sigma$ is given by
\begin{align*}
	|\Sigma|
	=
	\int_E \sqrt{1-|\nabla \nu|^2} \ed x
	\leq
	|D_l|.
\end{align*}
Note $|D_l^v| = |D_l|$ and the last equality holds if and only if $E=D$ and $\nabla \nu=0$, then the uniqueness follows.
\end{proof}

\subsection{Proof of the isoperimetric inequality for achronal hypersurface}
With the help of lemmas \ref{lem 8.2} and \ref{lem 8.3}, we can prove the isoperimetric inequality \eqref{eqn 8.1} in theorem \ref{thm 8.1}.
\begin{proof}[{\bf Proof of the isoperimetric inequality \eqref{eqn 8.1}}]
Assume that $C_f \subset C_R$. Consider the set $P_h^{e_0}$ of open finite lightcones obtained from $C_f$ by finite times of polarisations. There exists a sequence $\{ C_{f_k} \}_{k\in \mathbb{N}} \subset P_h^{e_0}$ converging to $C_l^{e_0}$ with 
\begin{align}
	n \omega_n l^{n-1}
	\leq
	\int_{\mathbb{S}^{n-1}} f^{n-1} \ed x,
	\label{eqn 8.2}
\end{align}
by proposition \ref{prop 4.20}. Suppose $\{\gamma_{1,k}, \cdots, \gamma_{l_k,k}\}$ is the sequence of reflections such that
\begin{align*}
	C_{f_k}
	=
	(C_f)^{\gamma_{1,k},\cdots,\gamma_{l_k,k},e_0}.
\end{align*}
Denote $\{\gamma_{1,k}, \cdots, \gamma_{l_k,k}\}$ by $\Gamma_k$ and the polarisation $\Omega^{\gamma_{1,k},\cdots,\gamma_{l_k,k},e_0}$ of a set $\Omega$ by $\Omega^{\Gamma_k,e_0}$. 

We want to construct a sequence of achronal hypersurfaces $\{ \Sigma^{\Gamma_k,e_0} \}$ by applying the polarisations corresponding to $\Gamma_k$ and show that $\Sigma^{\Gamma_k,e_0}$ converges to a closed achronal hypersurface in $C_l^v$ as the convergence of the graphs of Lipschitz functions over $x$-variables in the $\{t, x\}$ rectangular coordinate system. We want to adopt the construction in lemma \ref{lem 8.2}, therefore we introduce the following open causal convex set $I(o,\Sigma)$ and the corresponding null extension of $\Sigma$:
\begin{enumerate}[label=\roman*.]
\item
Let $I^-(\Sigma)$ be the causal past of $\Sigma$. Define $I(o,\Sigma)$ by
\begin{align*}
	I(o,\Sigma) = I^+(o) \cap I^-(\Sigma).
\end{align*}
Then $I(o,\Sigma)$ is open causal convex and $\Sigma \cap I_+(o) \subset \partial_+ I(o,\Sigma)$.

\item
Introduce the parameterisation function $\nu$ for $\partial_+ I(o,\Sigma)$,
\begin{align*}
	\partial_+ I(o,\Sigma)
	=
	\{ (t,x) \in \mathbb{M}^{n,1}: |x| < t = \nu(x) \}.
\end{align*}
Then we extend $\nu$ to $\mathbb{R}^n$ by setting $\nu(x) = |x|$ beyond the domain covered by $\partial_+ I(o,\Sigma)$. Then the graph of the extended $\nu$ is the extension of $\partial_+ I(o,\Sigma)$ by the lightcone $C_0$. Define the graph of the extended $\nu$ as the achronal hypersurface $\Sigma_{n.e.}$, which is an extension of $\Sigma$ by null hypersurfaces.
\end{enumerate}
From the above construction, we have that
\begin{align*}
	|\Sigma| = |\partial_+ I(o,\Sigma)| = |\Sigma_{n.e.} \cap C_R|
\end{align*}

Now we apply the corresponding polarisations of $\Gamma_k$ to $I(o,\Sigma)$ to obtain a sequence of open causal convex sets $\{ I(o,\Sigma)^{\Gamma_k,e_0} \}$ that
\begin{align*}
	I(o,\Sigma)^{\Gamma_k,e_0} 
	\subset 
	(\calD(C_f))^{\Gamma_k,e_0} 
	\subset
	\calD((C_f)^{\Gamma_k,e_0})
	= 
	\calD(C_{f_k}).
\end{align*}
By lemma \ref{lem 8.2}, we have
\begin{align*}
	|\partial_+ I(o,\Sigma)^{\Gamma_k,e_0}|
	=
	|\partial_+ I(o,\Sigma)|
	=
	|\Sigma|
\end{align*}
Let $\nu_k$ be the parameterisation function of $\partial_+ I(o,\Sigma)^{\Gamma_k,e_0}$ and extend $\nu_k$ by $\nu_k = |x|$ as in ii. above. Then there exists a subsequence $\{\nu_{k'}\}$ converges to a limit function $\bar{\nu}$. Then the sequence $\{ I(o,\Sigma)^{\Gamma_{k'},e_0} \}$ converges to a limit open causal convex set $\bar{I}$ w.r.t. the volume of symmetric difference. We have the following conclusions for $\bar{I}$ and $\bar{\nu}$:
\begin{enumerate}[label=\arabic*.]
\item
$\partial_+ \bar{I} = \{ (t,x)\in \mathbb{M}^{n,1}: t= \nu(x) > |x| \}$.

\item
$\partial_+ \bar{I} \subset \calD(C_l^{e_0})$.

\item
$|\partial_+ \bar{I}| \geq |\Sigma|$.
\end{enumerate}
Conclusions 1. and 2. simply follow from the convergence of $\nu_{k'}$. Conclusion 3. follows from the formulae of the areas $|\partial_+ \bar{I}|$, $\partial_+ I(o,\Sigma)^{\Gamma_k,e_0}$ that
\begin{align*}
	|\partial_+ I|
	=
	\int_{D_R} \sqrt{1- | \nabla \bar{\nu} |^2} \ed x,
	\quad
	|\partial_+ I(0,\Sigma)^{\Gamma_{k'},e_0}|
	=
	\int_{D_R} \sqrt{1- | \nabla \nu_{k'} |^2} \ed x.
\end{align*}
Since $-\sqrt{1-|p|^2} > -1$ is convex for $p \in D_1$, then the integral $I(u) = \int_{D_R} -\sqrt{1-|\nabla u|^2}$ is lower semi-continuous (see theorem 1.6 in \cite{St00}), therefore
\begin{align*}
	|\partial_+ \bar{I}|
	\geq
	\limsup_{k'\rightarrow +\infty} |\partial_+ I(0,\Sigma)^{\Gamma_{k'},e_0}|
	=
	|\Sigma|.
\end{align*}
Since $\partial_+ \bar{I} \subset \calD(C_l^{e_0})$, by lemma \ref{lem 8.3} we have
\begin{align*}
	|\Sigma|
	\leq
	|\partial_+ \bar{I}|
	\leq
	\max \{ |\Sigma'|: \Sigma' \subset \calD ( C_l^{e_0} ) \textit{ and is achronal}  \}
	=
	\omega_n l^n
\end{align*}
Substituting inequality \eqref{eqn 8.2}, we prove inequality \eqref{eqn 8.1}.
\end{proof}

\subsection{$(\mathbb{Z}_2)^{\times n}$-reflection symmetric achronal hypersurface}
To prove the case of equality for \eqref{eqn 8.1} in theorem \ref{thm 8.1}, we first show that the case can be reduced to the rotational symmetric case. Here we need a construction similar to the one of $(\mathbb{Z}_2)^{\times n}$-reflection symmetric lightcone in section \ref{sec 5.3}.

\begin{definition}\label{def 8.4}
Let $\calE$ be an open causal convex set in $I^+(o)$ with $\partial_- \calE \subset C_0$. Let $H$ be a timelike hyperplane passing through the origin $o$. We call $H$ an equal perimeter separation hyperplane of $\calE$ if $H$ separates $\partial_- \calE$, which is an open finite lightcone, to two parts with the same perimeter.
\end{definition}

Similar to lemma \ref{lem 5.5}, we have that the equal perimeter hyperplane also separates the area of strict future boundary equally for the set in the case of equality.
\begin{lemma}\label{lem 8.5}
Let $\calE$ be an open causal convex set in $J^+(o)$ with $\partial_- \calE \subset C_0$. Suppose that $\partial_+ \calE$ and $\partial_- \calE$ achieves the equality in \eqref{eqn 8.1}. Let $w$ be a spacelike vector and $H_w$ be the orthogonal timelike hyperplane of $w$ passing through $o$. If $H_w$ is an equal perimeter separation hyperplane of $\calE$, then $H_w$ separates $\partial^+ \calE$ to two parts with the same area, i.e.
\begin{align*}
	| \partial_+ \calE \cap H_w^+ |
	=
	| \partial_+ \calE \cap H_w^- |.
\end{align*}
Let $\gamma_w$ is the reflection about $H_w$, then
\begin{align*}
	(\partial_+ \calE \cap H_w^+) \cup \gamma_w (\partial_+ \calE \cap H_w^+),
	\quad
	(\partial_+ \calE \cap H_w^-) \cup \gamma_w (\partial_+ \calE \cap H_w^-)
\end{align*}
also achieve the equality in \eqref{eqn 8.1}.
\end{lemma}
\begin{proof}
The proof is similar to the one of lemma \ref{lem 5.5}, by contradiction and considering the achronal hypersurfaces
\begin{align*}
	&
	(\partial_+ \calE \cap H_w^+) \cup \gamma_w (\partial_+ \calE \cap H_w^+)
	\ \subset\
	\calD((\partial_- \calE \cap H_w^+) \cup \gamma_w (\partial_- \calE \cap H_w^+)),
\\
	&
	(\partial_+ \calE \cap H_w^-) \cup \gamma_w (\partial_+ \calE \cap H_w^-)
	\ \subset\
	\calD((\partial_- \calE \cap H_w^-) \cup \gamma_w (\partial_- \calE \cap H_w^-)).
\end{align*}
\end{proof}

We consider the special class of $(\mathbb{Z}_2)^{\times n}$-reflection symmetric achronal hypersurfaces. Similar to lemma \ref{lem 5.10}, we show that there exists a timelike axis of reflection symmetry for the archronal hypersurface in this special class.
\begin{lemma}\label{lem 8.6}
Let $\calE$ be a $(\mathbb{Z}_2)^{\times n}$-reflection symmetric open causal convex set in $J^+(o)$ with $\partial_- \calE \subset C_0$. Let $\{ w_1,\cdots, w_n\}$ be the orthogonal set of spacelike vectors corresponding to the $(\mathbb{Z}_2)^{\times n}$-reflection symmetry. Let $v$ be the future-directed timelike vector orthogonal to $w_1, \cdots, w_n$. Let $l_v$ be the line in the direction of $v$ through $o$. Then $l_v$ is the axis of reflection symmetry of $\calE$. Any timelike hyperplane is an equal perimeter separation hyperplane of $\calE$ if and only if that $l_v \subset H$.
\end{lemma}

The proof is the same as the one of lemma \ref{lem 5.10}, thus we omit it here.

\subsection{Increase area of achronal hypersurface at edge}

We present a lemma telling how to increase the area of a closed achronal hypersurface at the edge, the intersection of two  tangent hyperplanes. The basic idea is that if the achronal hypersurface is not rotationally symmetric, then we can apply the construction similar in lemma \ref{lem 8.7} to increase the area while preserving the perimeter.
\begin{lemma}\label{lem 8.7}
Introduce the following constructions. Let $\{t,x^1,\cdots, x^n\}$ be the rectangular coordinate system of $\mathbb{M}^{n,1}$.
\begin{enumerate}[label=\alph*.]
\item
Define $H_{a,1} = \{  a x_1 +  t = a\}$ and $H_{a,-1} = \{ -a x_1 + t = a \}$ where $a \in (0,1]$. See figure \ref{fig 21}.
\begin{figure}[h]
\centering
\begin{tikzpicture}[scale=.8, 
	declare function={ h1(\x,\y)=0.5*(1-\x); h2(\x,\y)=0.5*(1+\x); }
	]
	\begin{axis}[view = {50}{25},
		axis lines = left, 
		xlabel = $x^1$, 
		ylabel = $x^2$,
		ticklabel style = {font = \scriptsize},
		unit vector ratio = 1 1 2,
		grid
	]
	\addplot3 [domain=-1.2:1.2, 
		domain y=-2.5:2.5, 
		samples=50, 
		samples y=2, 
		surf, 
		shader = interp, 
	] 
		{min(h1(\x,\y),h2(\x,\y))} ;
	\end{axis}
\end{tikzpicture}
\caption{$H_{a,1}$ and $H_{a,-1}$}
\label{fig 21}
\end{figure}

\item
Consider the set 
\begin{align*}
	&
	L_{[-1,1],b} 
	=
	\{ (x^0, x^1,\cdots x^n): t=0, x^1 \in [-1,1], |x| \leq b  \},
\\
	&
	\begin{aligned}
		S_{[-1,1],b} 
		=
		\partial L_{[-1,1],b} 
		&=
		\{ (x^0, x^1,\cdots x^n): t=0, x^1 \in [-1,1], |x| = b  \},
	\\
		&\phantom{=}
		\cup
		\{ (x^0, x^1,\cdots x^n): t=0, x^1 = 1, |x| < b  \}
	\\
		&\phantom{=}
		\cup
		\{ (x^0, x^1,\cdots x^n): t=0, x^1 = -1, |x| < b  \} .
	\end{aligned}
\end{align*}
Let $C(S_{[-1,1],b})$ be the future outgoing null hypersurface emanating from $S_{[-1,1],b}$, which is the strict past boundary of the future of $L_{[-1,1],b} $, i.e. $C(S_{[-1,1],b}) = \partial_- ( I^+(L_{[-1,1],b}))$. Then
\begin{align*}
	C(S_{[-1,1],b})
	&=
	\big( C_{(t=-b,0)} \cap ( \{ t - x^1 \geq -1 \} \cup \{ t + x^1 \geq -1 \} ) \big)
\\
	&\phantom{=}
	\cup
	\big( I^+_{(t=-b,0)} \cap ( \{ t - x^1 = -1 \} \cup \{ t + x^1 = -1 \} ) \big).
\end{align*}
See figure \ref{fig 22}.
\begin{figure}[h]
\hfill
\begin{subfigure}[h]{0.45\textwidth}
\begin{tikzpicture}[scale=.8,
	declare function={h1(\x)= sqrt( 1.5^2 - \x^2);}
	]	
	\begin{axis}[view = {50}{25},
		axis lines = left, 
		zmin=0, zmax=1.5,
		xlabel = $x^1$, 
		ylabel = $x^2$,
		ticklabel style = {font = \scriptsize},
		unit vector ratio = 1 1 2,
		grid
	]
		\addplot3 [domain=-1:1, 
			domain y=-1:1, 
			samples=10, 
			samples y=10, 
			surf, 
			shader = interp,
		]
			(x, {y*h1(x)}, 0);
	\end{axis}
\end{tikzpicture}
\end{subfigure}
\hfill
\begin{subfigure}[h]{0.45\textwidth}
\centering
\begin{tikzpicture}[scale=.8, 
	declare function={ h4(\x,\y)=-1+\x; h5(\x,\y)=-1-\x); h3(\x,\y)=sqrt(\x^2+\y^2)-1.5; }]
	\begin{axis}[view = {50}{25},
		axis lines = left, 
		xlabel = $x^1$, 
		ylabel = $x^2$,
		ticklabel style = {font = \scriptsize},
		unit vector ratio = 1 1 2,
		grid
	]
	\addplot3 [domain=-2:2, 
		domain y=1:2.5, 
		samples=50, 
		samples y=40, 
		surf, 
		shader = interp,
	]
		{max(0,h3(\x,\y),h4(\x,\y),h5(\x,\y))};
	\addplot3 [domain=-2:2, 
		domain y=-1:1, 
		samples=50, 
		samples y=5, 
		surf, 
		shader = interp,
	]
		{max(0,h3(\x,\y),h4(\x,\y),h5(\x,\y))};
	\addplot3 [domain=-2:2, 
		domain y=-1:-2.5, 
		samples=50, 
		samples y=40, 
		surf, 
		shader = interp,
	]
		{max(0,h3(\x,\y),h4(\x,\y),h5(\x,\y))};
	\end{axis}
\end{tikzpicture}
\end{subfigure}
\caption{$L_{[-1,1],b}$ and $C(S_{[-1,1],b})$.}
\label{fig 22}
\end{figure}

\item
Introduce the causal convex set $I_{a,[-1,1],b} = I^- (H_{a,1}) \cap I^- (H_{a,-1}) \cap I^+(L_{[-1,1],b})$.
\end{enumerate}
Then we have that for the future and strict past boundaries of $I_{a,[-1,1],b}$,
\begin{align*}
	&
	\partial_- (I_{a,[-1,1],b})
	=
	L_{[-1,1],b}
	\cap
	\big( C_{(t=-b,0)}  \cap I^- (H_{a,1}) \cap I^- (H_{a,-1}) \cap \{ t>0 \}\big),
\\
	&
	|\partial_- (I_{a,[-1,1],b})|
	=
	2 \omega_{n-1} b^{n-1} + O ( b^{n-2} ).
\end{align*}
and
\begin{align*}
	&
	\partial_+ (I_{a,[-1,1],b})
	=
	I^+_{(t=-b,0)}
	\cap
	\big( (H_{a,1} \cap \{x^1\geq 0\}) \cup (H_{a,-1} \cap \{x^1\leq 0\}) \big),
\\
	&
	\partial_+ (I_{a,[-1,1],b})
	\subset
	\{ x^1 \in [-1,1], |x| \leq b+a \},
\\
	&
	|\partial_+ (I_{a,[-1,1],b})|
	\leq
	2 \omega_{n-1} \sqrt{1-a^2} \cdot (b+a)^{n-1}.
\end{align*}
See figure \ref{fig 23}.
\begin{figure}[h]
\centering
\begin{tikzpicture}[declare function={ h1(\x,\y)=0.5*(1-\x); h2(\x,\y)=0.5*(1+\x); h3(\x,\y)=sqrt(\x^2+\y^2)-1.5; }]
	\begin{axis}[view = {50}{25},
		axis lines = left, 
		xlabel = $x^1$, 
		ylabel = $x^2$,
		ticklabel style = {font = \scriptsize},
		unit vector ratio = 1 1 2,
		grid
	]
	\addplot3 [domain=-1.2:1.2, 
		domain y=2:2.5, 
		samples=50, 
		samples y=2, 
		surf, 
		shader = interp,
	]
		{min(max(0,h3(\x,\y)),h1(\x,\y),h2(\x,\y))};
	\addplot3 [domain=-1.2:1.2, 
		domain y=1:2, 
		samples=50, 
		samples y=20, 
		surf, 
		shader = interp,
	]
		{min(max(0,h3(\x,\y)),h1(\x,\y),h2(\x,\y))};
	\addplot3 [domain=-1.2:1.2, 
		domain y=-1:1, 
		samples=50, 
		samples y=5, 
		surf, 
		shader = interp,
	]
		{min(max(0,h3(\x,\y)),h1(\x,\y),h2(\x,\y))};
	\addplot3 [domain=-1.2:1.2, 
		domain y=-2:-1, 
		samples=50, 
		samples y=20, 
		surf, 
		shader = interp,
	]
		{min(max(0,h3(\x,\y)),h1(\x,\y),h2(\x,\y))};
	\addplot3 [domain=-1.2:1.2, 
		domain y=-2.5:-2, 
		samples=50, 
		samples y=2, 
		surf, 
		shader = interp,
	]
		{min(max(0,h3(\x,\y)),h1(\x,\y),h2(\x,\y))};
	\end{axis}
\end{tikzpicture}
\caption{$\partial_- (I_{a,[-1,1],b})$.}
\label{fig 23}
\end{figure}
Therefore
\begin{align*}
	|\partial_- (I_{a,[-1,1],b})|
	>
	|\partial_+ (I_{a,[-1,1],b})|
	+ 2 \omega_{n-1} (1-\sqrt{1-a^2}) b^{n-1}
	- O(b^{n-2}).
\end{align*}
Thus for each $a\in (0,1]$, there exists $b_a$ such that for any $b\geq b_a$,
\begin{align*}
	|\partial_- (I_{a,[-1,1],b})|
	>
	|\partial_+ (I_{a,[-1,1],b})|.
\end{align*}
\end{lemma}

\subsection{Reduction from $(\mathbb{Z}_2)^{\times n}$-reflection symmetry to rotational symmetry}

Before verifying the case of equality in theorem \ref{thm 8.1}, we prove the rotational symmetry for the $(\mathbb{Z}_2)^{\times n}$-reflection symmetric achronal hypersurface with the maximal area.
\begin{lemma}\label{lem 8.8}
Let $\calE$ be a $(\mathbb{Z}_2)^{\times n}$-reflection symmetric open causal convex set in $J^+(o)$ with $\partial_- \calE \subset C_0$. Let $\{ w_1,\cdots, w_n\}$ be the orthogonal set of spacelike vectors corresponding to the $(\mathbb{Z}_2)^{\times n}$-reflection symmetry, and $l_v$ be the axis of reflection symmetry. Suppose that the archronal hypersurface $\partial_+ \calE$ and the closed finite lightcone achieve the equality in \eqref{eqn 8.1}. Then $\calE$ is invariant under the spacelike rotations fixing the axis of reflection symmetry $l_v$.
\end{lemma}
\begin{proof}
Without loss of generality, we assume $v=e_0$. Let $p \in (\partial_+ \calE\setminus l_{e_0})$. Suppose that there exists a tangent plane $T_p \partial_+ \calE$ and the corresponding normal vector $n_p$ of $\partial_+ \calE$ at $p$. Let $P_{e_0,p}$ be the $2$-dim plane passing through $p$ and $l_{e_0}$. 

\vspace{5pt}
\noindent{\underline{\textsc{Claim.}} } $n_p \in P_{e_0,p}$.

Assuming the above claim, we prove the lemma. Introduce the parameterisation function $\nu$ of $\partial_+ \calE$
\begin{align*}
	\partial_+ \calE
	=
	\{ (t,x) \in \mathbb{M}^{n,1}: |x| <t = \nu(x) \},
\end{align*}
and we extend $\nu$ by defining $\nu(x) = |x|$ beyond the region of $x$ covered by $\partial_+ \calE$. Then for almost all $x$, the rotational vector derivative of $\nu$ at $x$ vanishes. This follows from that
\begin{align*}
	n_p \parallel (1, \partial_{x^1} \nu, \cdots, \partial_{x^n} \nu),
\end{align*}
and the rotational vector $(0,R_x) \perp P_{e_0,p}$. Then $\nu$ is rotational symmetric since $\nu$ is Lipschitz continuous and $R\nu =0$ for an arbitrary rotational vector field $R$ and almost all $x$. Thus $\partial_+ \calE$ is rotational symmetric about the axis $l_{e_0}$.
\end{proof}

\begin{proof}[Proof of \textsc{Claim}]
We prove the claim by contradiction. Assume that $n_p \not\in P_{e_0,p}$. Then there exists a timelike hyperplane $H$ passing through $l_{e_0}$ that $n_p \not\in H$. By lemmas \ref{lem 8.5} and \ref{lem 8.6}, $H$ is an equal perimeter separation hyperplane of $\calE$ and also separates the area of the strict future boundary $\partial_+ \calE$ equally.  Define $H_+$ be the half space which $n_p$ points into. Let  $\gamma$ be the reflection about $H$, then the archronal hypersurface $\Sigma_p$
\begin{align*}
	\Sigma_p
	=
	(\partial_+ \calE \cap H_-) \cup \gamma (\partial_+ \calE \cap H_-)
\end{align*}
achieves the equality in \eqref{eqn 8.1}. However, with the help of lemma \ref{lem 8.7}, we will show that there exists an interior perturbation of $\partial_+ \calE$ which increase the area, contradicting with \eqref{eqn 8.1}.

We build another rectangular coordinate system near $x$. For any $\epsilon > 0$, construct a coordinate system $\{\bart, \barx^1, \cdots, \barx^n\}$, such that
\begin{align*}
	&
	\eta
	=
	\epsilon^2 (-\ed \bart^2 + (\ed \barx^1)^2 + \cdots +  (\ed \barx^n)^2 ),
\\
	&
	H
	=
	\{\barx^1=0\},
	\quad
	H_-
	=
	\{\barx^1<0\},
\\
	&
	n_p
	\parallel
	\partial_{\bart} + a \partial_{\barx^1},
	\quad
	\bart(p) = a,
	\quad
	\barx^i(p) = 0.
\end{align*}
Let $\eta_{\epsilon} = \epsilon^{-2} \eta$. In $\{ \bart, \barx^i \}$ coordinate system, $\Sigma_p$ is a perturbation of $\partial_+ (I_{a,[-1,1],b})$ introduced in lemma \ref{lem 8.7} near the point $p$.  By lemma \ref{lem 8.7}, there exist $\delta$ sufficiently small and $b$ sufficiently large such that 
\begin{align*}
	|\partial_- (I_{a,[-1,1],b})|
	>
	|\partial_+ (I_{a,[-1,1],b})| + \delta b^{n-1}.
\end{align*}
We consider the open causal convex set $I_{\Sigma_p,b}$, which is a perturbation of $I_{a,[-1,1],b}$, defined as
\begin{align*}
	I_{\Sigma_p,b}
	=
	\calE \cap I^+(L_{\bart=0,\Sigma_p,b}),
	\quad
	L_{\bart=0,\Sigma_p,b}
	=
	\{\bart =0 \} \cap \calE \cap \{ |x| < b \}.
\end{align*}
where $L_{\bart=0,\Sigma_p,b}$ is a perturbation of $L_{[-1,1],b}$. In the following, we show that for sufficiently small $\epsilon$ and large $b$,
\begin{align*}
	|\partial_- (I_{\Sigma_p,b})|_{\eta_{\epsilon}}
	>
	|\partial_+ (I_{\Sigma_p,b})|_{\eta_{\epsilon}}.
\end{align*}
We prove the above in two cases $a\in (0,1)$ and $a=1$.

\vspace{5pt}
\noindent{\underline{\textsc{Case: $a\in (0,1)$.}} }
We shall compare $|\partial_+ (I_{\Sigma_p,b})|_{\eta_{\epsilon}}$ with $|\partial_+ (I_{a,[-1,1],b})|_{\eta_{\epsilon}}$, and $|\partial_+ (I_{\Sigma_p,b})|_{\eta_{\epsilon}}$ with $|\partial_+ (I_{a,[-1,1],b})|_{\eta_{\epsilon}}$.

Consider the orthogonal projection $P^{\perp}_{n_p}$ of the half $\partial_+ (I_{\Sigma_p,b}) \cap H_-$ and $\partial_+ (I_{a,[-1,1],b}) \cap H_-$ to the orthogonal hyperplane $H_{a,-1}$ of $n_p$, where
\begin{align*}
	P^{\perp}_{n_p} (q)
	=
	\{q + s\cdot n_p\} \cap H_{a,-1}.
\end{align*}
Then
\begin{align*}
	P^{\perp}_{n_p} (\partial_+ (I_{\Sigma_p,b}) \cap H_-)
	\subset
	U_{\eta_{\epsilon},o(1)}^{H_{a,-1}} (\partial_+ (I_{a,[-1,1],b}) \cap H_-)
\end{align*}
where $U_{\eta_{\epsilon},o(1)}^{H_{a,-1}} (\partial_+ (I_{a,[-1,1],b}) \cap H_-)$ is the $o(1)$-neighbourhood of $\partial_+ (I_{a,[-1,1],b}) \cap H_-$ in $(H_{a,-1}, \eta_{\epsilon})$. Therefore
\begin{align*}
	|\partial_+ (I_{\Sigma_p,b}) \cap H_-|_{\eta_{\epsilon}}
	&
	\leq
	| P^{\perp}_{n_p} (\partial_+ (I_{\Sigma_p,b}) \cap H_-) |_{\eta_{\epsilon}}
\\
	&
	\leq
	| U_{\eta_{\epsilon},o(1)}^{H_{a,-1}} (\partial_+ (I_{a,[-1,1],b}) \cap H_-) |_{\eta_{\epsilon}}
\\
	&
	\leq 
	| \partial_+ (I_{a,[-1,1],b}) \cap H_- |_{\eta_{\epsilon}}
	+ o(1) b^{n-1}.
\end{align*}
For $|\partial_+ (I_{\Sigma_p,b})|_{\eta_{\epsilon}}$ and $|\partial_+ (I_{a,[-1,1],b})|_{\eta_{\epsilon}}$, note that
\begin{align*}
	L_{[-1,1],b} 
	\subset
	U_{\eta_{\epsilon}, o(1)}^{H_{\bart=0}} (L_{\bart=0, \Sigma_p, b})
\end{align*}
where $U_{\eta_{\epsilon}, o(1)}^{H_{\bart=0}} (L_{\bart=0, \Sigma_p, b})$ is the $o(1)$-neighbourhood of $L_{\bart=0, \Sigma_p, b}$ in $(H_{\bart=0}, \eta_{\epsilon})$. Therefore we have that
\begin{align*}
	|\partial_- (I_{a,[-1,1],b})|_{\eta_{\epsilon}} 
	&=
	| L_{[-1,1],b} |_{\eta_{\epsilon}} 
\\
	&\leq
	| U_{\eta_{\epsilon}, o(1)}^{H_{\bart=0}} (L_{\bart=0, \Sigma_p, b}) |_{\eta_{\epsilon}} 
\\
	&\leq
	| L_{\bart=0, \Sigma_p, b} |_{\eta_{\epsilon}} + o(1) b^{n-1}
\\
	&=
	|\partial_- (I_{\Sigma_p,b})|_{\eta_{\epsilon}} + o(1) b^{n-1}.
\end{align*}
Then by lemma \ref{lem 8.7},
\begin{align*}
	|\partial_+(I_{\Sigma_p,b})|_{\eta_{\epsilon}}
	&\leq
	| \partial_+ (I_{a,[-1,1],b}) \cap H_- |_{\eta_{\epsilon}}
	+ o(1) b^{n-1}
\\
	&\leq
	|\partial_- (I_{a,[-1,1],b})|_{\eta_{\epsilon}} 
	- \delta b^{n-1}
	+ o(1) b^{n-1}
\\
	&\leq
	|\partial_- (I_{\Sigma_p,b})|_{\eta_{\epsilon}} 
	- \delta b^{n-1}
	+ o(1) b^{n-1}.
\end{align*}
Then there exist sufficiently small $\epsilon$ and large $b$,
\begin{align*}
	|\partial_+ (I_{\Sigma_p,b})|_{\eta_{\epsilon}}
	<
	|\partial_- (I_{\Sigma_p,b})|_{\eta_{\epsilon}}.
\end{align*}

\vspace{5pt}
\noindent{\underline{\textsc{Case: $a=1$.}} } 
Comparing $|\partial_+ (I_{\Sigma_p,b})|_{\eta_{\epsilon}}$ with $|\partial_+ (I_{a,[-1,1],b})|_{\eta_{\epsilon}}$, the following inequality obtained in the case $a\in(0,1)$ still holds
\begin{align*}
	2\omega_{n-1} b^{n-1} + O(b^{n-2})
	=
	|\partial_- (I_{a,[-1,1],b})|_{\eta_{\epsilon}} 
	\leq
	|\partial_- (I_{\Sigma_p,b})|_{\eta_{\epsilon}} + o(1) b^{n-1}
\end{align*}
To estimate $|\partial_+ (I_{\Sigma_p,b})|$, we have that
\begin{align*}
	&
	\partial_+  (I_{\Sigma_p,b}) \cap H
	\subset
	H \cap \{ | \bart - 1 | < o(1) \},
\\
	&
	\partial_+  (I_{\Sigma_p,b}) \cap H_- \cap \{\bart =0\}
	\subset
	\{\bart =0\} \cap \{ |\barx^1 + 1 | < o(1) \}.
\end{align*}
Then consider the spacelike hyperplane $H_{1-o(1),1+o(1)}$ that
\begin{align*}
	H_{1-o(1),1+o(1)}
	=
	\{ \frac{\bart}{1-o(1)} - \frac{\barx^1}{1+o(1)} = 1 \}
\end{align*}
Introduce the orthogonal projection $P^{\perp}_{1-o(1),1+o(1)}$ to the hyperplane $H_{1-o(1),1+o(1)}$, then
\begin{align*}
	P^{\perp}_{1-o(1),1+o(1)} ( \partial_+  (I_{\Sigma_p,b}) \cap H_- )
	\subset
	H_{1-o(1),1+o(1)} \cap \{ \bart > 0, \barx^1 < 0 , |x| < b+1\}.
\end{align*}
Therefore
\begin{align*}
	| \partial_+  (I_{\Sigma_p,b}) \cap H_- |_{\eta_{\epsilon}}
	&\leq
	| H_{1-o(1),1+o(1)} \cap \{ \bart > 0, \barx^1 < 0 , |x| < b+1\} |
\\
	&\leq
	o(1) (b+1)^{n-1}.
\end{align*}
Hence there exist sufficiently small $\epsilon$ and large $b$ that
\begin{align*}
	| \partial_+  (I_{\Sigma_p,b}) |_{\eta_{\epsilon}}
	\leq
	o(1)(b+1)^{n-1}
	\leq
	2\omega_{n-1} b^{n-1} - o(1) b^{n-1}
	\leq 
	| \partial_- (I_{\Sigma_p,b}) |_{\eta_{\epsilon}}
\end{align*}

\vspace{5pt}
\noindent{\underline{\textsc{Conclusion.}} } 
We conclude that if $n_p \not\in P_{e_0,p}$, then there exists an interior perturbation of $\Sigma_p$ increasing the area, which leads to a contradiction with the assumption that $\partial_+ \calE$ achieves the equality in \eqref{eqn 8.1}. Therefore the claim $n_p \in P_{e_0,p}$ is true.
\end{proof}

\subsection{Identification of case of equality}\label{sec 7.2.7}

We can identify the case of equality in \eqref{eqn 8.1}. We shall first introduce the following construction of a $(\mathbb{Z}_2)^n$-reflection symmetric achronal hypersurface, similar to construction \ref{con 5.12}.
\begin{construction}\label{con 8.9}
Let $\calE$ be an open causal convex set in $J^+(o)$ with $\partial_- \calE \subset C_0$. Let $p \in \partial_+ \calE$. Construct a $(\mathbb{Z}_2)^{\times n}$-reflection symmetric closed causal convex set $\bar{\calE}$ with $p\in \partial_+\bar{\calE}$.
\begin{enumerate}
\item
Choose a $2$-dimensional timelike plane $P_1$. Choose $w_1 \in P_1$ such that $H_{w_1}$ is the equal perimeter separation hyperplane of $\calE$ and $p \in H_{w_1}^+$. Then define $\calE_1$ as the positive reflection symmetrisation of $\calE$, i.e.
\begin{align*}
	\calE_1 
	=
	(\calE \cap H_{w_1}^+) 
	\cup 
	\gamma_1(\calE \cap H_{w_1}^+).
\end{align*}

\item
Assume that we obtain the orthogonal set of spacelike vectors $\{w_1, \cdots, w_k\}$ and the closed causal convex set $\calE_k$. Choose a $2$-dimensional timelike plane $P_{k+1} \perp \{w_1, \cdots, w_k\}$. Then find $w_{k+1} \in P_{k+1}$ such that $H_{w_{k+1}}$ is the equal perimeter separation hyperplane of $\calE_k$ and $p\in H_{w_{k+1}}^+$. Define $\calE_{k+1}$ as the positive reflection symmetrisation of $\calE_k$, i.e.
\begin{align*}
	\calE_{k+1} 
	=
	(\calE_k \cap H_{w_{k+1}}^+) 
	\cup 
	\gamma_1(\calE_k \cap H_{w_{k+1}}^+).
\end{align*}

\item
Define $\bar{\calE} = \calE_n$. Then $\bar{\calE}$ is $(\mathbb{Z}_2)^{\times n}$-reflection symmetric with the corresponding orthogonal set $\{w_1, \cdots, w_n\}$.
\end{enumerate}
\end{construction}

\begin{proof}[{\bf Proof of the case of equality}]
Let $\Sigma \subset \calD(C_f)$ be a closed achronal hypersurface achieving the equality in \eqref{eqn 8.1}. Consider the open causal convex set $\calE=I(o,\Sigma)$ that
\begin{align*}
	\calE
	=
	I(o,\Sigma)
	=
	I^+(o) \cap I^-(\Sigma),
\end{align*}
then $\partial_+ \calE$ also achieves the equality. Applying the construction \eqref{con 8.9} to $\calE$, we obtain a $(\mathbb{Z}_2)^{\times n}$-reflection symmetric causal convex set $\bar{\calE}$ and the corresponding orthogonal set $\{w_1,\cdots, w_n\}$. By lemma \ref{lem 8.5}, $\bar{\calE}$ also achieves the equality, thus $\bar{\calE}$ is rotational symmetric about the line $l_v$ where $v$ is future-directed unit timelike orthogonal to $\{w_1,\cdots, w_n\}$. Then by lemma \ref{lem 8.3}, we have
\begin{align*}
	\bar{\calE}
	=
	I^+(o) \cap I^-(H_l^v),
\end{align*}
where $n\omega_n l^{n-1} = \int_{\mathbb{S}^{n-1}} f^{n-1} \ed x$. By the construction of $\bar{\calE}$, we have that
\begin{align*}
	\partial_+ \calE \cap H_{w_1}^+ \cap \cdots \cap H_{w_n}^+
	=
	D_l^v \cap H_{w_1}^+ \cap \cdots \cap H_{w_n}^+,
\end{align*}
where $D_l^v = H_l^v \cap I_+(o)$ is the disk of radius $l$ in $H_l^v$. In the following, we shall show that
\begin{align*}
	\partial_+ \calE  \cap H_{w_1}^+ \cap \cdots \cap H_{w_k}^+
	=
	D_l^v \cap H_{w_1}^+ \cap \cdots \cap H_{w_k}^+,
\end{align*}
by induction on $k$ from $n$ to $0$.
\begin{enumerate}[label=\textit{\alph*}.]
\item
It is proved for the case $k=n$.

\item
Suppose the case of $k$ is proved. Consider the orthogonal set $\{w_1,\cdots, w_{k-1}, w'_k = -w_k\}$, which is obtained by reversing the direction of $w_k$. By definition
\begin{align*}
	&
	H_{w'_k}^+ = H_{w_k}^-,
\\
	&
	\partial_+ \calE  \cap (H_{w_1}^+ \cap \cdots \cap H_{w_{k-1}}^+) \cap H_{w_k}^-
	=
	\partial_+ \calE  \cap H_{w_1}^+ \cap \cdots \cap H_{w_{k-1}}^+ \cap H_{w'_k}^+.
\end{align*}
Then we extend $\{w_1, \cdots, w_{k-1}, w'_k\}$ to an orthogonal set $\{w_1, \cdots, w_{k-1}, w'_k, \cdots, w'_n\}$ as in construction \ref{con 8.9} and obtain the corresponding $(\mathbb{Z}_2)^{\times n}$-reflection symmetric open causal convex set $\bar{\calE}'$. Then applying the induction assumption to $\bar{\calE}'$, we have that
\begin{align*}
	\partial_+ \calE  \cap H_{w_1}^+ \cap \cdots \cap H_{w_{k-1}}^+ \cap H_{w'_k}^+
	=
	D_l^{v'} \cap H_{w_1}^+ \cap \cdots \cap H_{w_{k-1}}^+ \cap H_{w'_k}^+,
\end{align*}
where $v'$ is othorgonal to $\{w_1, \cdots, w_{k-1}, w'_k, \cdots, w'_n\}$. In order to show
\begin{align*}
	\partial_+ \calE  \cap H_{w_1}^+ \cap \cdots \cap H_{w_{k-1}}^+
	=
	D_l^v \cap H_{w_1}^+ \cap \cdots \cap H_{w_{k-1}}^+,
\end{align*}
it is sufficient to show that $v'$ is parallel to $v$, which is true by lemma \ref{lem 8.7}, otherwise one can increase the area of $\calE$ by an interior perturbation at the edge $\calE \cap H_{w_k}$ constructed in lemma \ref{lem 8.7}, contradictory with that $\calE$ achieving the equality in \eqref{eqn 8.1}. Thus we prove the case of $k-1$.
\end{enumerate}

Now we can finish the proof of the case of equality. Since $\Sigma \subset \partial_+ \calE = D_l^v$ and $\Sigma$ achieves the equality of \eqref{eqn 8.1}, $\Sigma=D_l^v$.
\end{proof}

\subsection{Implication to a functional inequality on the hyperbolic space}
Consider the hyperbolic model of $I^+(o)$ where
\begin{align*}
	\eta
	=
	-\ed l^2 + l^2 m_{-1}.
\end{align*}
where $m_{-1}$ is the metric of the hyperbolic space of constant curvature $-1$. Let $S_{-1}$ be the hyperboloid
\begin{align*}
	S_{-1}
	=
	\{ (t,x) \in \mathbb{M}^{n,1}: -t^2 + |x|^2 = -1\},
\end{align*}
then we can take $m_{-1} = \eta|_{S_{-1}}$. In geodesic coordinate system and the Poincar\'e disc model, the hyperbolic metric takes the forms
\begin{align*}
	m_{-1}
	=
	\ed s^2 + (\sinh s)^2 \circg
	=
	4\frac{(\ed v^1)^2 + \cdots + (\ed v^n)^2}{(1-|v|^2)^2}.
\end{align*}
Note that
\begin{align*}
	&
	t
	=
	l \cosh s,
	\quad
	r
	=
	l \sinh s,
\\
	&
	s
	=
	\tanh^{-1} |v| 
	= 
	\log \frac{1+|v|}{1-|v|}.
\end{align*} 
Let $\Sigma$ be a closed achronal hypersurface in $I^+(o)$ without boundary in $I^+(o)$. Suppose that $\Sigma$ is the graph of a function $f$ over the hyperbolic space in the hyperbolic coordinate system
\begin{align*}
	\Sigma
	=
	\{ (l, p):  l=f(p)>0, p\in S_{-1} \}.
\end{align*}
Then the metric and the volume form on $\Sigma$ take the form
\begin{align*}
	\eta|_{\Sigma}
	=
	f^2 m_{-1} - \ed f \otimes \ed f,
	\quad
	\dvol_{\Sigma}
	=
	f^n (1- |\ed \log f|_{m_{-1}}^2)^{\frac{1}{2}} \dvol_{m_{-1}}.
\end{align*}

To see the behaviour of $\Sigma$ near $C_0$, we adopt the geodesic coordinate system of $S_{-1}$. The parameterisation function of $\Sigma$ is
\begin{align*}
	l = f(s,\vartheta).
\end{align*}
Then in the spatial polar coordinate system 
\begin{align*}
	t = f(s, \vartheta) \cosh s,
	\quad
	r= f(s, \vartheta) \sinh s,
	\quad
	t+ r = f(s,\vartheta) \exp s.
\end{align*}
Note the function $s \mapsto f(s,\vartheta) \exp s$ is monotonically non-decreasing, then we define the limit function $f_{\infty}$
\begin{align*}
	f_{\infty}(\vartheta) = \lim_{s\rightarrow +\infty} f(s,\vartheta) \exp s.
\end{align*}
Then considering the open finite lightcone $C_{f_{\infty}}$,\footnote{$f_{\infty}$ is lower semicontinuous since it is the limit of a monotonically nondecreasing sequence of continuous functions.} we have that
\begin{align*}
	\Sigma
	\subset
	\calD(C_{f_{\infty}}).
\end{align*}
The perimeter of $C_{f_{\infty}}$ is
\begin{align*}
	| P(C_{f_{\infty}}) |
	=
	\int_{\mathbb{S}^{n-1}} (f_{\infty})^{n-1}(\vartheta) \dvol_{\circg} (\vartheta).
\end{align*}
Applying theorem \ref{thm 8.1} to $\Sigma$ and $C_{f_{\infty}}$, we obtain the following inequality as a corollary.
\begin{corollary}\label{coro 8.10}
Let $f$ be a positive Lipschitz continuous function on the hyperbolic space $(S_{-1}, m_{-1})$, with the Lipschitz constant of $\log f$ being no more than $1$. Let $o$ be a point in $S_{-1}$ and $(s,\vartheta)$ be the geodesic coordinate system at $o$ where the metric takes the form
\begin{align*}
	m_{-1}
	=
	\ed s^2 + (\sinh s)^2 \circg_{\vartheta}.
\end{align*}
Define the function $f_{\infty}$ on the sphere at infinity by
\begin{align*}
	f_{\infty}(\vartheta)
	=
	\lim_{s\rightarrow +\infty} f(s,\vartheta) \exp s.\footnotemark
\end{align*}
\footnotetext{Note that $f_{\infty}$ does not depend on the choice of the origin $o$.}
Then we have
\begin{align}
	\frac{ \int_{S_{-1}} f^n (1- |\ed \log f|_{m_{-1}}^2)^{\frac{1}{2}} \dvol_{m_{-1}} }{\omega_n}
	\leq
	\Big( \frac{\int_{\mathbb{S}^{n-1}} f_{\infty}^{n-1} \dvol_{\circg}}{n \omega_n} \Big)^{\frac{n}{n-1}}.
	\label{eqn 8.3}
\end{align}
The equality is achieved when
\begin{align*}
	&
	f(s,\vartheta) 
	= 
	\frac{c}{\cosh s + k \sinh s \cdot \cos ( d_{\mathbb{S}^{n-1}} (\vartheta, \vartheta_0))},
\\
	&
	f_{\infty}
	=
	\frac{c}{1 + k \cos ( d_{\mathbb{S}^{n-1}} (\vartheta, \vartheta_0))},
	\quad
	c>0,
	|k| < 1.
\end{align*}
where $\vartheta_0 \in \mathbb{S}^{n-1}$ and $d_{\mathbb{S}^{n-1}}$ is the distance in $(\mathbb{S}^{n-1}, \circg)$.
\end{corollary}

\section{Isoperimetric inequality for spacelike hypersurface with boundary in hyperboloid}\label{sec 9}
In this section, we show that the spacelike hyperplane disk has the maximal area among all the achronal hypersurfaces with the boundary of the same perimeter in the hyperboloid.
\begin{theorem}\label{thm 9.1}
Let $\Sigma \subset I^+(o)$ be a closed achronal hypersurface in $I^+(o)$. Suppose that $\Sigma$ is contained in the domain of dependence of a set $E\subset S_{-1}$.
\begin{enumerate}[label=\alph*.]
\item
Assume that $E$ is a set of finite perimeter in the hyperboloid $S_{-1}$. The following isoperimetric inequality holds that
\begin{align}
	\frac{|\Sigma|}{\omega_n}
	\leq
	\Big( \frac{P(E)}{n \omega_n} \Big)^{\frac{n}{n-1}},
	\label{eqn 9.1}
\end{align}
where $P(E)$ is the perimeter of $E$. The equality is achieved when $E$ contains a closed geodesic ball $\overline{B}$ in $(S_{-1}, \eta_{-1})$ and $|E\setminus \overline{B} | =0$.

\item
Assume that $E$ has finite measure in the hyperboloid $S_{-1}$. The following inequality holds that
\begin{align}
	|\Sigma|
	\leq
	| D_{B_{r_{|E|}}} |.
	\label{eqn 9.2}
\end{align}
where $D_{B_{r_{|E|}}}$ is a spacelike hyperplane disk with the same boundary as $B_{r_{|E|}}$ which is the geodesic ball in $(S_{-1}, \eta|_{S_{-1}})$ of the area $|E|$. The equality is achieved when there is a closed geodesic ball $\overline{B_{r_{|E|}}}$ in $(S_{-1}, \eta|_{S_{-1}})$ such that $\overline{B_{r_{|E|}}} \subset E$, $|E \setminus \overline{B_{r_{|E|}}} | = 0$ and $\Sigma = D_{B_{r_{|E|}}}$.
\end{enumerate}
\end{theorem}

We show first that one can increase the area of an achronal hypersurface which has some portion in $I^- (S_{-1})$.
\begin{lemma}\label{lem 9.2}
Let $\Sigma \subset I^+(o)$ be a closed achronal hypersurface in $I^+(o)$ contained in the domain of dependence of a set $E\subset S_{-1}$ with finite measure. Let $\Sigma_- = \Sigma \cap I^-(S_{-1})$. Define
\begin{align*}
	\Sigma'
	=
	S_{-1} \cap I^+( \Sigma_- ),
	\quad
	\bar{\Sigma} = (\Sigma \cap J^+(S_{-1}) ) \cup \Sigma'.
\end{align*}
Then we have that $\bar{\Sigma}$ is achronal and
\begin{align*}
	|\Sigma|
	\leq
	|\bar{\Sigma}|
\end{align*}
and the equality is achieved when $\Sigma_- = \emptyset$.
\end{lemma}
\begin{proof}
To prove that $\bar{\Sigma}$ is achronal, it is sufficient to show that $I^-(p) \cap I^+(q') = \emptyset$ where $p \in \Sigma \cap J^+(S_{-1})$ and $q' \in \Sigma'$. Since $q' \in \Sigma' \subset I^+(\Sigma_-)$, there exists $q \in \Sigma_-$ such that $I^+(q') \subset I^+(q)$. Hence
\begin{align*}
	I^-(p) \cap I^+(q') 
	\subset
	I^-(p) \cap I^+(q)
	=
	\emptyset.
\end{align*}

To prove $|\Sigma| \leq |\bar{\Sigma}|$, it is sufficient to show $|\Sigma_-| \leq |\Sigma'|$. Denote $\partial_-(I^+(\Sigma_-)) \cap I^-(S_{-1})$ by $\Sigma_{-,\partial_- I^+}$. Note that $\Sigma_- \subset \Sigma_{-,\partial_- I^+}$, then it is sufficient to show that  $| \Sigma_{-,\partial_- I^+} | \leq |\Sigma'|$. Parameterising $\Sigma_{-,\partial_- I^+}$ and $\Sigma'$ in the hyperbolic coordinate system of $I^+(o)$,
\begin{align*}
	&
	\Sigma_{-,\partial_- I_+}
	=
	\{ (l,p): l = f_{-,\partial_- I_+}(p), p \in \Omega \subset S_{-1} \},
\\
	&
	\Sigma'
	=
	\{ (l,p): l = 1, p \in \Omega \subset S_{-1} \}.
\end{align*}
Then by the formula of the area, we have
\begin{align*}
	| \Sigma_{-,\partial_-I^+} |
	=
	\int_{\Omega} 
	(f_{-,\partial_- I^+})^n 
	( 1 - \vert \ed \log f_{-,\partial_- I^+} \vert_{m_{-1}}^2 )^{\frac{1}{2}} 
	\dvol_{m_{-1}}
	\leq
	|\Omega|_{m_{-1}}
	=
	| \Sigma' |,
\end{align*}
where the equality is achieved when $\Omega = \emptyset$.
\end{proof}

We apply the Lorentz polarisation to hypersurfaces in $J^+(S_{-1})$.
\begin{lemma}\label{lem 9.3}
Let $\Sigma$ be a closed achronal hypersurface in $\calD(E) \cap J^+(S_{-1})$ where $E \subset S_{-1}$. Let $\calE$ be $J^-(\Sigma) \cap J^+(S_{-1})$. Then $\calE$ is a closed causal convex set. Let $v$ be a future-directed vector, $H$ be a timelike hyperplane and $v\not\in H$. Let $\calE^{\gamma,v}$ and $E^{\gamma,v}$ be the Lorentz polarisations of $\calE$ and $E$ about $(H,v)$. Then $\calE^{\gamma,v} \subset \calD(E^{\gamma,v}) \cap J^+(S_{-1})$ is closed causal convex and
\begin{align*}
	\Sigma 
	\subset 
	\bar{\partial}_+ \calE,
	\quad
	| \Sigma |
	\leq
	| \bar{\partial}_+ \calE |
	=
	| \bar{\partial}_+ (\calE^{\gamma,v}) |.
\end{align*}
\end{lemma}
\begin{proof}
$\calE = J^-(\Sigma) \cap J^+(S^{-1}) \Rightarrow \Sigma \subset \bar{\partial}_+ \calE$. We only need to verify that $| \bar{\partial}_+ \calE | = | \bar{\partial}_+ (\calE^{\gamma,v}) |$. It is similar to the proof of lemma \ref{lem 8.2}, therefore we omit it here.
\end{proof}

We can prove theorem \ref{thm 9.1} now.
\begin{proof}
Part \textit{a.} follows from part \textit{b.} by the isoperimetric inequality in the hyperbolic space, thus it is sufficient to prove part \textit{b.} in the following.

Define 
\begin{align*}
	\bar{\Sigma} = (\Sigma \cap J^+(S_{-1}) \cup (S_{-1} \cap I^+(\Sigma_-))
\end{align*}
as in lemma \ref{lem 9.2}. Introduce
\begin{align*}
	\calE = J^-(\bar{\Sigma}) \cap J^+(S_{-1})
\end{align*}
as in lemma \ref{lem 9.3}. There exists a sequence of families of reflections $ \Gamma_k =\{ \gamma_{1,k}, \cdots, \gamma_{i_k,k} \}$ such that the sequence of sets $\{E_k\}$, the polarisation of $E$ by $\{ \Gamma_k, e_0 \}$
\begin{align*}
	E_k = E^{\gamma_{1,k}, \cdots, \gamma_{i_k,k}, e_0} = E^{\Gamma_k, e_0},
\end{align*}
converges to the geodesic ball $B_{r_{|E|}}(e_0)$ in $S_{-1}$, i.e.
\begin{align*}
	\lim_{k \rightarrow +\infty} | E_k \triangle \overline{B_{r_{|E|}} (e_0)} | = 0.
\end{align*}
Then consider the sequences of closed causal convex sets $\{ \calE_k = \calE^{\Gamma_k,e_0} \}$ and future boundaries $\{\Sigma_k = \bar{\partial}_+ \calE_k \}$. By lemma \ref{lem 9.3}, we have that
\begin{align*}
	\Sigma_k 
	\subset 
	\calD(E_k),
	\quad
	|\Sigma|
	\leq
	|\Sigma_k|.
\end{align*}
We show that $\{ \Sigma_k \}$ converges in the following sense. We extends $\Sigma_k$ by $S_{-1}$ that defining
\begin{align*}
	\bar{\calE}_k 
	= 
	\calE_k \cup S_{-1},
	\quad
	\bar{\Sigma}_k
	=
	\bar{\partial}_+ \bar{\calE}_k,
\end{align*}
then $\Sigma_k \subset \bar{\Sigma}_k$. Then $\{ \bar{\Sigma}_k \}$ converges to $\bar{\Sigma}_{\infty}$ as the graphs of the functions over $S_{-1}$ in the hyperbolic coordinate system. Let $\bar{f}_k$ and $\bar{f}_{\infty}$ be the parameterisation functions of $\bar{\Sigma}_k$ and $\bar{\Sigma}_{\infty}$ respectively in the hyperbolic coordinate system, then
\begin{align*}
	f_{\infty}(p)
	=
	\lim_{k\rightarrow +\infty} f_k(p).
\end{align*}
Then there exists a subsequence $\{k'\}$ such that
\begin{align*}
	\int_{\overline{B_{r_{|E|}}}} f_{\infty}^n (p) \sqrt{1-|\nabla f_{\infty} |^2} \dvol_{m_{-1}}
	\geq
	\lim_{k' \rightarrow +\infty}
	\int_{\overline{B_{r_{|E|}}}} f_{k'}^n (p) \sqrt{1-|\nabla f_{\infty} |^2} \dvol_{m_{-1}}.
\end{align*}
Then
\begin{align*}
	\int_{\overline{B_{r_{|E|}}}} f_{\infty}^n (p) \sqrt{1-|\nabla f_{\infty} |^2} \dvol_{m_{-1}}
	\geq
	\lim_{k' \rightarrow +\infty}
	\int_{E_{k'}} f_{k'}^n (p) \sqrt{1-|\nabla f_{\infty} |^2} \dvol_{m_{-1}},
\end{align*}
since $\lim_{k \rightarrow +\infty} | E_k \triangle B_{r_{|E|}} (e_0) | = 0$ and $\{ f_{k'} \}$ is uniformly bounded. Let
\begin{align*}
	\Sigma_{\infty}
	=
	\{ (l,p): l = f_{\infty}(p), p \in \overline{B_{r_{|E|}}} \},
\end{align*}
then $\Sigma_{\infty} = \bar{\Sigma}_{\infty} \cap \calD( \overline{B_{r_{|E|}}} )$, since $f_k(p) = 1$ in $S_{-1} \setminus E_k$, which implies that $f_{\infty}(p) =1 a.e.$ in $S_{-1} \setminus \overline{B_{r_{|E|}}} $. Hence we obtain that
\begin{align*}
	| E |
	\leq 
	| E_{\infty} |
	\leq
	| D_{B_{r_{|E|}}} |,
\end{align*}
where the last inequality follows from the same argument in the proof of lemma \ref{lem 8.3}. 

The case of equality in \eqref{eqn 9.2} follows the same strategy as in the proof of the case of equality in \eqref{eqn 8.1} in section \ref{sec 7.2.7}, where we construct $(\mathbb{Z}_2)^{\times n}$-reflection symmetric closed achronal hypersurface similarly as in construction \ref{con 8.9} and  use lemma \ref{lem 8.8} to show the spherical symmetry of $(\mathbb{Z}_2)^{\times n}$-reflection symmetric closed achronal hypersurface achieving the equality.
\end{proof}

\section*{Acknowledgements}
\addcontentsline{toc}{section}{Acknowledgements}
The author is supported by the National Natural Science Foundation of China under Grant No. 12201338. The author thanks Demetrios Christodoulou and Pin Yu for valuable discussions.





\Address

\end{document}